\documentclass[12pt,leqno,openany]{book}

\usepackage{makeidx}
\makeindex

\def\diam{\mathop{\rm diam}}
\def\dim{\mathop{\rm dim}}
\def\dist{\mathop{\rm dist}}
\def\Lip{\mathop{\rm Lip}}
\def\supp{\mathop{\rm supp}}

\newtheorem{theorem}{Theorem}
\newtheorem{lemma}[theorem]{Lemma}
\newtheorem{proposition}[theorem]{Proposition}
\newtheorem{sublemma}[theorem]{Sublemma}
\newtheorem{definition}[theorem]{Definition}
\newtheorem{corollary}[theorem]{Corollary}
\newtheorem{problem}[theorem]{Problem}
\newtheorem{remark}[theorem]{Remark}
\newtheorem{claim}[theorem]{Claim}
\newtheorem{assumptions}[theorem]{Assumptions}
\newtheorem{examples}[theorem]{Examples}
\newtheorem{question}[theorem]{Question}
\newtheorem{sassumptions}[theorem]{Standing Assumptions}
\newtheorem{sassumption}[theorem]{Standing Assumption}
\newtheorem{conjecture}[theorem]{Conjecture}

\newcommand{\begintheorem}{\addtocounter{equation}{1}\begin{theorem}}
\newcommand{\beginlemma}{\addtocounter{equation}{1}\begin{lemma}}
\newcommand{\beginproposition}{\addtocounter{equation}{1}\begin{proposition}}
\newcommand{\beginsublemma}{\addtocounter{equation}{1}\begin{sublemma}}
\newcommand{\begindefinition}{\addtocounter{equation}{1}\begin{definition}}
\newcommand{\begincorollary}{\addtocounter{equation}{1}\begin{corollary}}
\newcommand{\beginproblem}{\addtocounter{equation}{1}\begin{problem}}
\newcommand{\beginremark}{\addtocounter{equation}{1}\begin{remark}}
\newcommand{\beginclaim}{\addtocounter{equation}{1}\begin{claim}}
\newcommand{\beginassumptions}{\addtocounter{equation}{1}\begin{assumptions}}
\newcommand{\beginexamples}{\addtocounter{equation}{1}\begin{examples}}
\newcommand{\beginquestion}{\addtocounter{equation}{1}\begin{question}}
\newcommand{\beginsassumptions}{\addtocounter{equation}{1}\begin{sassumptions}}
\newcommand{\beginsassumption}{\addtocounter{equation}{1}\begin{sassumption}}
\newcommand{\beginconjecture}{\addtocounter{equation}{1}\begin{conjecture}}

\begin{document}

\frontmatter

\title{Notes on Metrics, Measures, and Dimensions}

\author{Stephen Semmes \\
	Department of Mathematics \\
	Rice University}

\date{}

\maketitle

\chapter{Preface}

	In this monograph various notions related to metric spaces are
considered, including Hausdorff-type measures and dimensions,
Lipschitz mappings, and the Hausdorff distance between nonempty closed
and bounded subsets of a metric space.  Some familiarity with basic
topics in analysis such as Riemann integrals, open and closed sets,
and continuous functions is assumed, as in \cite{Devinatz, Goldberg,
Rudin1}, for instance.  There is no attempt to be exhaustive in any
way, and more information in a variety of directions can be found in
the bibliography.

\tableofcontents

\mainmatter

\chapter{Metric spaces}
\label{chapter about metric spaces}

\section{Real numbers}
\label{section on real numbers}
\index{real numbers}

	Let us begin by reviewing some aspects of the real numbers.
As usual, we write ${\bf R}$ for the real numbers, and ${\bf R}^n$ for
the $n$-tuples of real numbers, where $n$ is a positive integer.  We
also write ${\bf Z}$ for the integers, and ${\bf Z}_+$ for the
positive integers.

	If $a$, $b$ are real numbers with $a < b$, then the
intervals\index{intervals} $(a, b)$, $[a, b)$, $(a, b]$, and $[a, b]$
can be defined as the sets of real numbers $x$ which satisfy $a < x <
b$, $a \le x < b$, $a < x \le b$, and $a \le x \le b$, respectively.
For $[a, b]$ we also allow $a = b$, in which case $[a, b]$ consists of
the one point.  These intervals are all \emph{bounded}, and their
\emph{lengths}\index{length of an interval} are defined to be $b - a$.
It can sometimes be convenient to consider unbounded intervals of the
form $(-\infty, b)$, $(-\infty, b]$, $(a, \infty)$, $[a, \infty)$, and
$(-infty, \infty) = {\bf R}$, whose lengths are defined to be
$+\infty$.

	If $J$ denotes an interval in the real line, then we may
write $|J|$ for the length of $J$.

	Let $A$ be a subset of ${\bf R}$.  A real number $b$ is said
to be an \emph{upper bound}\index{upper bound for a set of real
numbers} if $a \le b$ for all $a \in A$.  A real number $c$ is said to
be the \emph{least upper bound},\index{least upper bound of a set of
real numbers} or \emph{supremum},\index{supremum of a set of real
numbers} if it is an upper bound for $A$, and if $c \le b$ for every
other upper bound $b$ of $A$.  It is easy to see from the definition
that the supremum is unique if it exists.  A basic property of the
real numbers is completeness\index{completeness of the real numbers in
terms of the existence of suprema} for the ordering, which means that
every nonempty subset $A$ of ${\bf R}$ which has an upper bound has a
least upper bound.  If $A$ is a nonempty subset of ${\bf R}$ which has
no upper bound in the real numbers, then we shall make the convention
of saying that the supremum is equal to $+ \infty$.  The supremum of
$A$ is denoted $\sup A$.

	Similarly, if $A$ is a subset of ${\bf R}$ and $d$ is a real
number, then $d$ is said to be a \emph{lower bound}\index{lower bound
of a set of real numbers} of $A$ if $d \le a$ for all $a$ in $A$.  A
real number $h$ is said to be the \emph{greatest lower
bound},\index{greatest lower bound of a set of real numbers} or
\emph{infimum},\index{infimum of a set of real numbers} of $A$ if it
is a lower bound for $A$, and if $d \le h$ for all other lower bounds
of $A$.  Again, it is easy to see from the definition that the infimum
is unique if it exists.  From the completeness property of the real
numbers described in the previous paragraph, one can show that a
nonempty subset of ${\bf R}$ with a lower bound has an infimum.  One
way to do this is to show that the infimum of $A$ is equal to the
supremum of $-A = \{-a : a \in A\}$, where $-A$ is nonempty and has an
upper bound because $A$ is nonempty and has a lower bound.  Another
approach is to obtain the infimum of $A$ as the supremum of the set
of lower bounds of $A$.  This second argument has the nice feature that
it applies to any linearly-ordered set.  Compare with \cite{Rudin1}.

	If $A$ is a nonempty subset of ${\bf R}$ which does not have a
lower bound, then we shall make the convention that the infimum is equal 
to $- \infty$.  The infimum of $A$ is denoted $\inf A$ in either case.
Often in these notes we shall be concerned with infima and suprema
of nonempty sets of nonnegative real numbers, so that the infimum exists
as a real number, but the supremum might be $+\infty$.

	The real numbers together with $-\infty$, $+\infty$ are called
the \emph{extended real numbers},\index{extended real numbers} as in
\cite{Rudin1}.  The usual ordering on the real numbers extends to the
extended real numbers, and to some extent the arithmetic operations do
too.  For instance, if $x$ is a real number, then $x + \infty$ is
defined to be $+\infty$, $x - \infty$ is defined to be $-\infty$,
while $\infty - \infty$ is normally left undefined.  If $A$ is a
nonempty set of extended real numbers, the infimum and supremum of $A$
can be defined for $A$ in the same way as above, as extended real
numbers.

	Let $I$ be a nonempty set, and suppose that to each element
$i$ of $I$ is associated a nonnegative real number $a_i$.  Consider
the sum
\begin{equation}
\label{sum_{i in I} a_i}
	\sum_{i \in I} a_i.
\end{equation}
If $I$ is finite, then this sum is defined in the usual way.  If $I$
is infinite, then we can define this sum to be the supremum of
\begin{equation}
	\sum_{i \in F} a_i,
\end{equation}
where $F$ runs through all finite subsets of $I$.  Thus the sum may
be equal to $+\infty$.  Sometimes the terms $a_i$ may be allowed to
be nonnegative extended real numbers, and the whole sum is automatically
equal to $+\infty$ if any one of the terms is equal to $+\infty$.

	If the sum (\ref{sum_{i in I} a_i}), then for each $\epsilon >
0$ the number of $i \in I$ such that $a_i \ge \epsilon$ is at most
$1/\epsilon$ times the sum (\ref{sum_{i in I} a_i}).  In particular,
the number of such $i \in I$ is finite.  As a consequence, if the sum
(\ref{sum_{i in I} a_i}) is finite, then the set of $i \in I$ such
that $a_i > 0$ is at most countable.  In these notes, we shall only be
concerned with situations in which $I$ is at most countable anyway.

	The \emph{absolute value}\index{absolute value of a real
number} $|x|$ of a real number $x$ is defined to be equal to $x$ when
$x \ge 0$ and to $-x$ when $x \le 0$.  Thus the absolute value of $x$
is always a nonnegative real number, and one can check that $|x + y|
\le |x| + |y|$ for all real numbers $x$, $y$.  Also, the absolute
value of a product of two real numbers is equal to the product of
their absolute values.

	Let $n$ be a positive integer.  If $x$, $y$ are elements of
${\bf R}^n$ and $t$ is a real number, then we can define $x + y$ and
$t \, x$ as elements of ${\bf R}^n$ by adding the coordinates of $x$
and $y$ in the first case, and multiplying the coordinates of $x$ by
$t$ in the second.  For each $x$ in ${\bf R}^n$, we define $|x|$, the
standard Euclidean norm of $x$, by
\begin{equation}
	|x| = \biggl(\sum_{j=1}^n x_j^2 \biggr)^{1/2}.
\end{equation}
It is easy to see that $|t \, x| = |t| \, |x|$ for all $x \in {\bf R}^n$
and all $t \in {\bf R}$, and it is well-known that
\begin{equation}
	|x + y| \le |x| + |y|
\end{equation}
for all $x, y \in {\bf R}^n$. 

	If $X$ is a set and $E$ is a subset of $X$, then the
\emph{indicator function on $X$ associated to $E$}\index{indicator
function on a set associated to a subset} is denoted ${\bf 1}_E$
and defined by ${\bf 1}_E(x) = 0$ when $x \in X \backslash E$,
${\bf 1}_E(x) = 1$ when $x \in E$.

\section{Some basic notions and results}
\label{section on basic notions and results about metric spaces}
\index{metric spaces}

	Let $(M, d(x,y))$ be a metric space.  For the record, this
means that $M$ is a nonempty set, and that $d(x,y)$ is a nonnegative
real-valued function on $M \times M$ such that $d(x, y) = 0$ if and
only if $x = y$, $d(x, y) = d(y, x)$ for all $x, y \in M$, and
\begin{equation}
\label{the triangle inequality}
\index{triangle inequality}
	d(x, z) \le d(x, y) + d(y, z)
\end{equation}
for all $x, y, z \in M$.  For each positive integer $n$, ${\bf R}^n$
equipped with the standard Euclidean metric\index{standard Euclidean
metric on ${\bf R}^n$} $|x - y|$ defines a metric space.  Of course
every nonempty subset of $M$ can be viewed as a metric space itself,
just by restricting the metric on $M$ to $Y$.

	If $x$ is an element of $M$ and $r$ is a positive real number,
then we write $B(x, r)$ for the open ball in $M$ with center $x$ and
radius $r$, defined by
\begin{equation}
\label{def of B(x, r)}
	B(x, r) = \{y \in M : d(x,y) < r\}.
\end{equation}
Similarly, the closed ball with center $x$ and radius $r$ is denoted
$\overline{B}(x, r)$ and defined by
\begin{equation}
\label{def of overline{B}(x, r)}
	\overline{B}(x, r) = \{y \in M : d(x,y) \le r\}.
\end{equation}

\subsection{Open and closed subsets and normality of metric spaces}

	If $E$ is a subset of $M$ and $p$ is an element of $M$, then
$p$ is said to be a limit point of $E$ if for every $\epsilon > 0$
there is a point $x$ in $E$ such that $x \ne p$ and $d(x, p) <
\epsilon$.  If every limit point of $E$ is also an element of $E$,
then $E$ is said to be a closed\index{closed subset of a metric space}
subset of $M$.  This is the same as saying that a point $p$ in $M$
lies in $E$ if for every $\epsilon > 0$ there is a point $x$ in $E$
such that $d(x, p) < \epsilon$.  A closed ball $\overline{B}(x, r)$ as
in (\ref{def of overline{B}(x, r)}) is always a closed subset of $M$,
and the empty set and $M$ itself are closed subsets of $M$.

	For any subset $E$, the \emph{closure}\index{closure of a
subset of a metric space} is denoted $\overline{E}$ and is defined to
be the union of $E$ and the set of limit points of $E$.  Equivalently,
a point $p$ lies in $\overline{E}$ if for every $\epsilon > 0$ there
is an $x$ in $E$ such that $d(p, x) < \epsilon$.  It is not difficult
to show that the closure of $E$ is always a closed subset of $M$.
Thus $\overline{E} = E$ if and only if $E$ is closed, and the closure
of $E$ can be described as the smallest closed subset of $M$ that
contains $E$.  If $D$ is a subset of $M$, then $D$ is
\emph{dense}\index{dense subsets of a metric space} in $M$ if the
closure of $D$ is equal to $M$.

	A subset $U$ of $M$ is said to be \emph{open}\index{open
subsets of a metric space} if for every $x$ in $U$ there is an $r > 0$
so that $B(x, r)$ is contained in $U$.  An open ball $B(x, r)$ as in
(\ref{def of B(x, r)}) is always an open subset of $M$, as are the
empty set and $M$ itself.  A basic fact is that $U$ is open if and
only if $M \backslash U$ is a closed subset of $M$.  For any subset
$E$ of $M$, the \emph{interior}\index{interior of a subset of a metric
space} $E^\circ$ of $E$ is the set of points $x$ in $E$ for which
there is an $r > 0$ so that $B(x, r)$ is contained in $E$.  One can
check that $E^\circ$ is always an open subset of $E$, $E^\circ = E$ if
and only if $E$ is an open set, and $E^\circ$ is the largest open
subset of $E$.  

	Notice that the interior of $E$ is the same as the complement
in $M$ of the closure of the complement of $E$, i.e., 
\begin{equation}
	E^\circ = M \backslash \overline{(M \backslash E)}.
\end{equation}
If $E_1$, $E_2$ are subsets of $M$, then it is not difficult to show
that
\begin{equation}
	\overline{(E_1 \cup E_2)} = \overline{E_1} \cup \overline{E_2},
		\quad \overline{(E_1 \cap E_2)} 
			\subseteq \overline{E_1} \cap \overline{E_2}
\end{equation}
and
\begin{equation}
	E_1^\circ \cup E_2^\circ \subseteq (E_1 \cup E_2)^\circ,
	    	\quad E_1^\circ \cap E_2^\circ = (E_1 \cap E_2)^\circ.
\end{equation}

	The \emph{boundary}\index{boundary of a subset of a metric
space} of a subset $E$ of $M$ is denoted $\partial E$ and defined by
\begin{equation}
	\partial E = \overline{E} \backslash E^\circ.
\end{equation}
This is equivalent to
\begin{equation}
	\partial E = \overline{E} \cap \overline{(M \backslash E)},
\end{equation}
and in particular
\begin{equation}
	\partial E = \partial (M \backslash E).
\end{equation}
The empty set and $M$ itself have empty boundary in $M$.  Observe that
if $E_1$, $E_2$ are subsets of $M$, then
\begin{equation}
	\partial (E_1 \cup E_2) \subseteq \partial E_1 \cup \partial E_2
\end{equation}
and
\begin{equation}
	\partial (E_1 \cap E_2) \subseteq \partial E_1 \cup \partial E_2.
\end{equation}

	As a metric space, $M$ enjoys the Hausdorff property that for
every pair of distinct elements $x$, $y$ in $M$ there are open sets
$U_1$, $U_2$ such that $x \in U_1$, $y \in U_2$, and $U_1 \cap U_2 =
\emptyset$.  For this one can simply take $U_1$, $U_2$ to be the balls
$B(x, d(x,y)/2)$, $B(y, d(x,y)/2)$.  A stronger version of this is
called ``regularity'' and says that if $x$ is an element of $M$ and
$F$ is a nonempty closed subset of $M$ which does not contain $x$ as
an element, then there are open subsets $V_1$, $V_2$ of $M$ such that
$x \in V_1$, $F \subseteq V_2$, and $V_1 \cap V_2 = \emptyset$.  To
see that this holds, choose $r > 0$ so that $B(x, r) \cap F =
\emptyset$, and then choose $V_1$ to be $B(x, r/2)$, and $V_2$ to be
the set of $z$ in $M$ such that $d(x,z) > r/2$.  An even stronger
separation property is called normality and states that if $F_1$ and
$F_2$ are nonempty closed subsets of $M$ such that $F_1 \cap F_2 =
\emptyset$, then there are open subsets $W_1$, $W_2$ of $M$ such that
$F_1 \subseteq W_1$, $F_2 \subseteq W_2$, and $W_1 \cap W_2 =
\emptyset$.  To find such a pair $W_1$, $W_2$, choose for each $x$ in
$F_1$ a positive real number $r_1(x)$ so that $B(x, r_1(x)) \cap F_2 =
\emptyset$, and choose for each $y$ in $F_2$ a positive real number
$r_2(y)$ so that $B(y, r_2(y)) \cap F_1 = \emptyset$.  Define
$W_1$, $W_2$ by
\begin{equation}
	W_1 = \bigcup_{x \in F_1} B(x, r_1(x)/2), \quad
		W_2 = \bigcup_{y \in F_2} B(y, r_2(y)/2).
\end{equation}
Clearly $W_1$, $W_2$ are open subsets of $M$ which contain $F_1$,
$F_2$, respectively, as subsets.  It is not difficult to show that
$W_1 \cap W_2 = \emptyset$.  One can strengthen this a bit further by
saying that there are open subsets $\widetilde{W}_1$,
$\widetilde{W}_2$ of $M$ such that $F_1 \subseteq \widetilde{W}_1$,
$F_2 \subseteq \widetilde{W}_2$, and the closures of $\widetilde{W}_1$
and $\widetilde{W}_2$ are disjoint.  In fact, this can be obtained as
a consequence of the previous normality property, by first using
normality to find disjoint open subsets $W_1$, $W_2$ of $M$ such that
$F_1 \subseteq W_1$, $F_2 \subseteq W_2$, and then using normality
again to find open subsets $\widetilde{W}_1$, $\widetilde{W}_2$ of $M$
such that $F_1 \subseteq \widetilde{W}_1$, $F_2 \subseteq
\widetilde{W}_2$ and the closures of $\widetilde{W}_1$,
$\widetilde{W}_2$ are contained in $W_1$, $W_2$, respectively.

\subsection{Bounded, compact, and totally bounded sets}

	A subset of $M$ is said to be \emph{bounded}\index{bounded 
subset of a metric space} if it is contained in a ball.  If $A$ is
a bounded subset of $M$, then the \emph{diameter}\index{diameter of
a subset of a metric space} $\diam A$ of $A$ is defined by
\begin{equation}
\label{def of diam A}
	\diam A = \sup \{d(x,y) : x, y \in A\}.
\end{equation}
In this case the diameter of $A$ is a finite real number.  If $A$ is
an unbounded subset of $M$, then the diameter of $A$ is defined to be
$+\infty$.  We make the convention that the diameter of the empty set
is equal to $0$.

	A sequence $\{p_j\}_{j=1}^\infty$ in $M$ is said to converge
to a point $p$ in $M$ if for every $\epsilon > 0$ there is a positive
integer $N$ such that
\begin{equation}
	d(p_j, p) < \epsilon \qquad \hbox{for all } j \ge N.
\end{equation}
In this case we write
\begin{equation}
	\lim_{j \to \infty} p_j = p,
\end{equation}
and one can check that the limit $p$ is unique.  A sequence
$\{p_j\}_{j=1}^\infty$ in $M$ is said to be a \emph{Cauchy
sequence}\index{Cauchy sequence in a metric space} if for every
$\epsilon > 0$ there is a positive integer $N$ so that
\begin{equation}
	d(p_j, p_k) < \epsilon \qquad\hbox{for all } j, k \ge N.
\end{equation}
It is easy to see that a convergent sequence is a Cauchy sequence.  A
metric space in which every Cauchy sequence is a convergent sequence
is said to be \emph{complete}.\index{complete metric spaces} A
well-known fact is that the Euclidean spaces ${\bf R}^n$ with their
standard metrics are complete metric spaces.  Also, any metric
space can be embedded isometrically into a complete metric space,
where the image of the first metric space is dense in the second one,
and this ``completion'' is unique up to isometric equivalence.

	A subset $E$ of $M$ is closed if and only if every sequence of
points in $E$ which converges in $M$ has its limit in $E$.  For any
subset $E$ of $M$, the closure $\overline{E}$ of $E$ is equal to the
set of points in $M$ which occur as limits of sequences of points in
$E$.  This includes constant sequences $\{p_j\}_{j=1}^\infty$, in
which all the $p_j$'s are the same, so that the elements of $E$ are
limits of sequences of points in $E$.  A subset $D$ of $M$ is dense in
$M$ if and only if every element of $M$ occurs as the limit of a
sequence of points in $D$.

	If $\{x_j\}_{j=1}^\infty$ is a sequence of points in some set $X$,
and if $\{j_l\}_{l=1}^\infty$ is a strictly increasing sequence of
positive integers, so that 
\begin{equation}
	j_1 < j_2 < j_3 < \cdots,
\end{equation}
then the sequence $\{x_{j_l}\}_{l=1}^\infty$ obtained by restricting
the $\{x_j\}_{j=1}^\infty$ to $\{j_l\}_{l=1}^\infty$ is called a
\emph{subsequence}\index{subsequence of a sequence} of
$\{x_j\}_{j=1}^\infty$.  In particular, $\{x_j\}_{j=1}^\infty$ is a
subsequence of itself, by taking $j_l = l$ for all $l$.  Note that a
Cauchy sequence in a metric space converges as soon as it has a
convergent subsequence.

	In these notes we shall often be concerned with subsets of
metric spaces which are \emph{compact}.\index{compact subsets of a
metric space} As usual, compactness can be characterized in a number
of equivalent ways.  In brief, a subset $K$ of $M$ is compact if (a)
every for every covering of $K$ by open subsets of $M$ there is a
subcovering of $K$ from this covering with only finitely many
elements, or (b) every infinite subset of $K$ has a limit point in
$K$, or (c) every sequence of points in $K$ has a subsequence which
converges.  If $Y$ is a nonempty subset of $M$ and $K$ is contained in
$Y$, then $K$ is compact as a subset of $M$ if and only if $K$ is
compact as a subset of $Y$, viewed as a metric space itself, with the
restriction of the metric from $M$.  In particular, if $K$ is
nonempty, we can take $Y = K$, so that $K$ is compact as a subset of
$M$ if and only if $K$ is compact as a subset of itself, viewed as a
metric space on its own, with the restriction of the metric from $M$.
This is quite different from the properties of being open or closed.

	In the case of a Euclidean space ${\bf R}^n$, the compact
subsets are exactly the closed and bounded subsets, by the
Heine--Borel theorem.

	A subset $A$ of $M$ is said to be \emph{totally
bounded}\index{totally bounded subsets of a metric space} if for every
$\epsilon > 0$ there is a finite collection of balls in $M$ with
radius $\epsilon$ whose union contains $A$ as a subset.  Compact sets
are totally bounded, and subsets of totally bounded sets are totally
bounded.  As with compactness, a subset $A$ of $M$ is totally bounded
as a subset of $M$ if and only if it is totally bounded as a subset of
itself, viewed as a metric space on its own, with the restriction of
the metric from $M$.  

	A well-known theorem states that a metric space is compact if
and only if it is totally bounded and complete.  This is equivalent to
saying that a subset of a complete metric space is compact if and only
if it is closed and totally bounded.  Another way to look at this is
that a subset $A$ of $M$ is totally bounded if and only if every
sequence of points in $A$ has a subsequence which is a Cauchy
sequence.  Roughly speaking, these results can be treated in much the
same way as classical results about compactness properties of closed
and bounded subsets of ${\bf R}^n$.

\subsection{Separable metric spaces}

	Recall that $M$ is said to be \emph{separable}\index{separable
metric spaces} if there is a subset $D$ of $M$ which is dense and at
most countable.  In this survey we shall be primarily interested in
separable metric spaces.  For example, the Euclidean spaces ${\bf
R}^n$ with their standard metrics are separable, because the set of
points in ${\bf R}^n$ with rational coordinates is countable and
dense.

	A metric space which is totally bounded is also separable.
In fact, a metric space $M$ is separable if and only if for each
$\epsilon > 0$ there is a subset $D_\epsilon$ of $M$ which is at most
countable and such that for each $x \in M$ there is a $y \in D_\epsilon$
that satisfies $d(x, y) < \epsilon$.  The property of being totally
bounded is the same as saying that one can do this with $D_\epsilon$
being finite for each $\epsilon > 0$.

	If $M$ is a metric space, a collection $\mathcal{B}$ of open
subsets of $M$ is said to be a \emph{basis for the
topology}\index{basis for the topology of a metric space} if every
open subset of $M$ can be expressed as a union of subsets of $M$ which
are elements of $\mathcal{B}$.  This is equivalent to saying that for
each $x$ in $M$ and each positive real number $r$ there is an open set
$U$ which is an element of $\mathcal{B}$ such that $x$ is an element
of $U$ and $U$ is a subset of $B(x, r)$.  Separability of $M$ can be
characterized by the existence of a basis $\mathcal{B}$ for the
topology of $M$ which has at most countably-many elements.  For if
$D$ is a dense subset of $M$ which is at most countable, then the
collection of balls of the form $B(y, 1/n)$, where $y$ runs through
the elements of $D$ and $n$ runs through the positive integers,
is a basis for the topology of $M$ which is at most countable.  Conversely,
if $\mathcal{B}$ is a basis for the topology of $M$ which is at most
countable, then one can get a subset $D$ of $M$ which is dense and
at most countable by choosing a point in each element of $\mathcal{B}$
and putting it into $D$.

	A basic fact about separable metric spaces is that every
subset $E$ of a separable metric space $M$ contains a subset which is
dense in $E$ and at most countable.  In other words, every nonempty
subset of $M$ defines a separable metric space itself, with the
restriction of the metric from $M$.  This is not difficult to show,
but notice that a particular dense subset of $M$ may not have any
elements in any particular subset $E$ of $M$.

	A subset of a metric space is said to be \emph{countably
compact}\index{countable compactness} if for every open covering of
the set there is a subcovering with at most countably many elements.
Clearly a compact set is countably compact.  A well-known result
states that a metric space $M$ is countably compact if and only if it
is separable.  Indeed, if $M$ is countably compact, then for every $r
> 0$ there is a family of balls of radius $r$ with at most countably
many elements which covers $M$, and this implies that $M$ is
separable.  Conversely, suppose that $M$ is separable, and that
$\{U_\alpha\}_{\alpha \in A}$ is an open covering of $M$.  Let
$\mathcal{B}$ be a basis for the topology of $M$ which is at most
countable.  Define $\mathcal{B}_1$ to be the subset of $\mathcal{B}$
consisting of those open subsets $V$ of $M$ in $\mathcal{B}$ which are
contained in a $U_\alpha$ for some $\alpha \in A$.  Since
$\mathcal{B}$ is a basis for the topology of $M$, each $U_\alpha$ is a
union of elements of $\mathcal{B}$.  Hence
\begin{equation}
	\bigcup_{V \in \mathcal{B}_1} V = \bigcup_{\alpha \in A} U_\alpha = M.
\end{equation}
For each $V$ in $\mathcal{B}_1$, choose an $\alpha(V) \in A$ such that
$V \subseteq U_{\alpha(V)}$.  Then
\begin{equation}
	\bigcup_{V \in \mathcal{B}_1} U_{\alpha(V)} = M,
\end{equation}
and $\{U_{\alpha(V)} : V \in \mathcal{B}_1\}$ is a subcovering of $M$
from the covering $\{U_\alpha\}_{\alpha \in A}$ which has at most
countably many elements.

\section{The Baire category theorem}
\label{section on Baire category theorem}

	Let $(M, d(x,y))$ be a metric space.  If $U$, $V$ are dense
open subsets of $M$, then the intersection $U \cap V$ is also a dense
open subset of $M$.  This is not difficult to verify.  More precisely,
if $U$ is a dense open subset of $M$ and $D$ is a dense subset of $M$,
then the intersection $U \cap D$ is a dense subset of $M$.  Of course
the intersection of two dense sets in general can be dense.

	Now suppose that $\{U_j\}_{j=1}^\infty$ is a sequence of dense
open subsets of $M$.  In general, the intersection of all of the
$U_j$'s might be the empty set.  This occurs when $M$ is the set of
rational numbers, equipped with the usual metric, for instance.  The
Baire category theorem\index{Baire category theorem} states that if
$M$ is a complete metric space, then in fact $\bigcap_{j=1}^\infty U_j$
is a dense subset of $M$.

	Let us briefly review the proof.  Let $B_0$ be any closed ball
in $M$.  It suffices to show that $B_0$ contains an element of
$\bigcap_{j=1}^\infty U_j$.  Since $U_1$ is a dense open subset of
$M$, there is a closed ball $B_1$ in $M$ such that $B_1 \subseteq
B_0$, $B_1 \subseteq U_1$, and the radius of $B_1$ is less than $1$.
Similarly, because $U_2$ is a dense open subset of $M$, there is a
closed ball $B_2$ in $M$ such that $B_2 \subseteq B_1$, $B_2 \subseteq
U_2$, and the radius of $B_2$ is less than $1/2$.  Proceeding in this
manner, for each positive integer $j$ there is a closed ball $B_n$
such that $B_j \subseteq B_{j-1}$, $B_j \subseteq U_j$, and the radius
of $B_j$ is less then $1/j$.  It follows that the sequence of centers
of the $B_j$'s is a Cauchy sequence, and so the assumption that $M$ is
complete implies that this sequence converges to a point $p$ in $M$.
One can check that $p$ lies in each $B_j$, and hence in each $U_j$,
as well as the given ball $B_0$.

	There is a complementary picture for closed sets with empty
interior.  In any metric space $M$, if $E$ and $F$ are closed sets with
empty interior, then $E \cup F$ is also a closed set with empty interior.
This does not work without the assumption that the sets be closed,
since the real line is the union of the sets of rational and irrational
numbers, each of which has empty interior.  In a complete metric space
the union of sequence of closed sets with empty interior also has empty
interior, although this union is not a closed set in general.

	If $F$ is any nonempty closed subset of $M$ and $j$ is a
positive integer, consider the set $U_j$ which is the union of
the balls $B(x, 1/j)$ with $x \in F$.  This is an open subset of $M$
which contains $F$ as a subset.  It is not difficult to check that
\begin{equation}
	F = \bigcap_{j=1}^\infty U_j.
\end{equation}
Thus every closed subset of a metric space can be expressed as the
intersection of a sequence of open sets.  By passing to complements
one obtains that every open subset of a metric space can be expressed
as the union of a sequence of closed subsets.

\section{Continuous mappings}
\label{section on continuous mappings}

	Let $(M, d(x,y))$ and $(N, \rho(u,v))$ be metric spaces, and
let $E$ be a nonempty subset of $M$.  A mapping $f$ from $E$ to $N$ is
said to be \emph{continuous}\index{continuous mappings between metric
spaces} if for every $x$ in $E$ and every $\epsilon > 0$ there is a
$\delta > 0$ so that
\begin{equation}
	\rho(f(x), f(y)) < \epsilon \quad\hbox{for all } y \in E
					\hbox{ such that } d(y,x) < \delta.
\end{equation}
The mapping $f$ is said to be \emph{uniformly continuous}\index{uniformly
continuous mappings between metric spaces} if for every $\epsilon > 0$
there is a $\delta > 0$ such that 
\begin{equation}
	\rho(f(x), f(y)) < \epsilon \quad\hbox{for all } x, y \in E
					\hbox{ such that } d(y,x) < \delta.
\end{equation}
A basic result states that if $E$ is compact and $f : E \to N$ is
continuous, then $f$ is also uniformly continuous.

	Let us assume for simplicity that $E = M$, to which one can
reduce anyway.  Another basic result states that $f : M \to N$ is
continuous if and only if $f^{-1}(U)$ is an open subset of $M$ for all
open subsets $U$ of $N$.  This is equivalent to saying that $f : M \to
N$ is continuous if and only if $f^{-1}(A)$ is a closed subset of $M$
whenever $A$ is a closed subset of $N$.  An alternate characterization
states that $f : M \to N$ is continuous if and only if for every
sequence $\{x_j\}_{j=1}^\infty$ of points in $M$ which converges to
some point $x$, the sequence $\{f(x_j)\}_{j=1}^\infty$ converges in
$N$ to $f(x)$.

	If $M_1$, $M_2$, $M_3$ are metric spaces and $h_1 : M_1 \to
M_2$, $h_2 : M_2 \to M_3$ are continuous mappings, then the
composition $h_2 \circ h_1$, defined by $(h_2 \circ h_1)(x) =
h_2(h_1(x))$, is a continuous mapping from $M_1$ to $M_3$.  If $h_1 :
M_1 \to M_2$, $h_2 : M_2 \to M_3$ are uniformly continuous, then $h_2
\circ h_1$ is a uniformly continuous mapping from $M_1$ to $M_3$.
These statements are easy to check.

	A mapping $f : M \to N$ is said to be
\emph{bounded}\index{bounded mapping between two metric spaces} if the
image $f(M)$ of $f$ is contained in a bounded subset of $N$.  A basic
result states that if $f : M \to N$ is continuous and $K$ is a compact
subset of $M$, then $f(K)$ is a compact subset of $N$.  In particular,
$f$ is bounded in this case.  Also, if $E$ is a subset of $M$ which is
totally bounded and $f : M \to N$ is uniformly continuous, then $f(E)$
is totally bounded.

	A \emph{homeomorphism}\index{homeomorphism between metric
spaces} from $M$ to $N$ is a one-to-one mapping $f$ of $M$ onto $N$
such that $f : M \to N$ and $f^{-1} : N \to M$ are both continuous.
If $f : M \to N$ is continuous, $M$ is compact, and $f$ is a
one-to-one mapping of $M$ onto $N$, then $f$ is a homeomorphism of $M$
onto $N$.  For instance, $f$ maps closed subsets of $M$ to closed subsets
of $N$ in this case, since closed subsets of $M$ are compact.	

	Suppose that $f_1$, $f_2$ are continuous mappings from $M$ to
$N$, and that $E$ is a dense subset of $M$.  If $f_1$ and $f_2$ are
equal on $E$, then they are equal on all of $M$.  In other words, a
continuous mapping is uniquely determined by its restriction to a
dense subset.  Conversely, suppose that $D$ is a dense subset of $M$,
and that $f$ is a continuous mapping from $D$ into $N$.  In order for
there to be a extension of $f$ to a continuous mapping from $M$ to
$N$, it is sufficient that $f$ be uniformly continuous, and that $N$
be complete.  In this case the extension will also be uniformly
continuous.

\subsection{Real-valued functions}

	Of course constant functions are always continuous.  If
$f_1(x)$, $f_2(x)$ are real-valued continuous functions on $M$, then
the sum $f_1(x) + f_2(x)$ and the product $f_1(x) \cdot f_2(x)$ are
also continuous functions on $M$.  The sum of uniformly continuous
real-valued functions is uniformly continuous, but this does not work
for the product in general.  It does work if at least one of the two
functions is bounded.  If $f(x)$ is a continuous real-valued function
on $M$ and $E$ is a subset of $M$ such that $f(x) \ne 0$ for all $x
\in M$, then $1/f(x)$ is a continuous function on $E$.  If we
assume that $f(x)$ is uniformly continuous, then it may not be true
that $1/f(x)$ is uniformly continuous on $E$, but this is the case if
there is a positive real number $c$ such that $|f(x)| \ge c$ for all
$x$ in $E$. 

	On ${\bf R}^n$, the $n$ coordinate functions
\begin{equation}
	x = (x_1, \ldots, x_n) \mapsto x_j,
\end{equation}
$1 \le j \le n$, are continuous, as one can easily verify.  As a
result, polynomials on ${\bf R}^n$ are continuous functions, and
rational functions, which are quotients of polynomials, are continuous
on the set where the denominator does not vanish.

\beginlemma
\label{d(x,p) is continuous}
For each element $p$ of $M$, $d(x,p)$ is uniformly continuous as a
real-valued function of $x$ on $M$.
\end{lemma}

	This follows from the triangle inequality.  Namely,
\begin{equation}
	d(x, p) \le d(y, p) + d(x,y)
\end{equation}
and
\begin{equation}
	d(y, p) \le d(x, p) + d(x,y),
\end{equation}
so that
\begin{equation}
	|d(x, p) - d(y, p)| \le d(x,y).
\end{equation}

	Let $A$ be a nonempty subset of $M$, and set
\begin{equation}
\label{def of dist(x, A)}
	\dist(x, A) = \inf \{d(x,z) : z \in A\}
\end{equation}
for $x$ in $M$.  Note that 
\begin{equation}
	\dist(x, \overline{A}) = \dist(x, A)
\end{equation}
for all $x$ in $M$, and that $\dist(x, A) = 0$ if and only if
$x$ lies in the closure of $A$.
	
\beginlemma
\label{distance to a subset is continuous}
As a function of $x$ on $M$, $\dist(x, A)$ is uniformly continuous.
\end{lemma}

	For this we use the triangle inequality again.  If $x$ and $y$
are arbitrary elements of $M$ and $z$ is an arbitrary element of $A$,
then
\begin{equation}
	\dist(x, A) \le d(x, z) \le d(x, y) + d(y, z).
\end{equation}
By taking the infimum over $z$ it follows that
\begin{equation}
	\dist(x, A) \le d(x, y) + \dist(y, A).
\end{equation}
Similarly,
\begin{equation}
	\dist(y, A) \le d(x, y) + \dist(x, A),
\end{equation}
so that
\begin{equation}
	|\dist(x, A) - \dist(y, A)| \le d(x, y).
\end{equation}

	Note that for any nonempty subset $A$ of $M$,
\begin{equation}
	\dist(x, A) = \dist(x, \overline{A})
\end{equation}
for all $x$ in $M$, as one can verify from the definitions.

        Suppose that $F_1$, $F_2$ are disjoint nonempty closed subsets
of $M$.  Define a real-valued function $\phi(x)$ on $M$ by
\begin{equation}
\label{def of phi}
        \phi(x) = \frac{\dist(x, F_1)}{\dist(x, F_1) + \dist(x, F_2)}.
\end{equation}
Because $F_1$ and $F_2$ are disjoint and closed, the denominator
in the definition of $\phi(x)$ is never equal to $0$.  Thus
$\phi(x)$ is a continuous function on $M$ which takes values in the
interval $[0,1]$, and it satisfies $\phi(x) = 0$ if and only if
$x \in F_1$ and $\phi(x) = 1$ if and only if $x \in F_2$.

	In general this function $\phi(x)$ may not be uniformly
continuous.  For example, if $M = {\bf R} \backslash \{0\}$, $F_1 =
(-\infty, 0)$, and $F_2 = (0, \infty)$, then $F_1$, $F_2$ are disjoint
nonempty closed subsets of $M$, $\phi(x)$ is defined completely by
$\phi(x) = 0$ when $x \in F_1$ and $\phi(x) = 1$ when $x \in F_2$,
and $\phi$ is clearly not uniformly continuous.

	Let us say that a mapping $f : M \to N$ is \emph{locally
bounded}\index{locally bounded mappings between metric spaces} if for
every $p$ in $M$ there is an $r > 0$ so that $f(B(p, r))$ is a bounded
subset of $N$.  Thus continuous mappings are automatically locally
bounded.  Let us say that $f : M \to N$ is \emph{locally uniformly
continuous}\index{locally uniformly continuous mappings between metric
spaces} if for each $p \in M$ there is an $r > 0$ such that the
restriction of $f$ to $B(p, r)$ is uniformly continuous.  Ordinary
continuity is already a local condition, and so we do not need to
define ``local continuity''.  For real-valued functions, the sum and
product of locally bounded functions is locally bounded, and the sum
and product of locally uniformly continuous functions is locally
uniformly continuous.  If $f(x)$ is a real valued function on $M$
which is locally uniformly continuous and which satisfies $f(x) \ne 0$
for all $x \in M$, then $1/f(x)$ is locally uniformly continuous on
$M$.

	In particular, the function $\phi(x)$ defined in (\ref{def of
phi}) is always locally uniformly continuous, even if it may not be
uniformly continuous.  If there is an $\epsilon > 0$ so that
$d(y, z) \ge \epsilon$ for all $y \in F_1$ and $z \in F_2$, then
$\phi(x)$ is uniformly continuous.

	Recall that a metric space is said to be \emph{locally
compact}\index{locally compact spaces} if for every point $x$ in the
space there is are subsets $K$, $U$ of the space such that $K$ is
compact, $U$ is open, $x \in U$, and $U \subseteq K$.  Thus every
continuous mapping from a locally compact metric space to another
metric space is locally uniformly continuous.

	Let us say that a mapping between two metric spaces is
\emph{countably uniformly continuous}\index{countably uniformly
continuous mappings between metric spaces} if the domain can be
expressed as the union of an at most countable family of subsets, on
each of which the mapping is uniformly continuous.  One can also
assume that the mapping is bounded on each of these subsets, by making
a further decomposition as necessary.  As a basic property of such
mappings, if the domain can be expressed as the countable union of
totally bounded subsets, then so can the image of the mapping.

	The sum or product of two real-valued countably uniformly
continuous functions is also countably uniformly continuous.  This is not
difficult to show, using the observation that if $\{E_i\}_{i \in I}$
and $\{F_j\}_{j \in J}$ are families of subsets of $M$ such that $I$,
$J$ are at most countable and
\begin{equation}
	\bigcup_{i \in I} E_i = \bigcup_{j \in J} F_j = M,
\end{equation}
then $\{E_i \cap F_j \}_{(i, j) \in I \times J}$ is an at most
countable family of subsets of $M$ whose union is also equal to $M$.
Similarly, the reciprocal of a countably uniformly continuous mapping
which does not take the value $0$ is countably uniformly continuous.
A locally uniformly continuous mapping between two metric spaces is
countably uniformly continuous if the domain is separable, since that
implies countable compactness.

	A subset $E$ of $M$ is said to be \emph{$\sigma$-compact} if
it can be expressed as the union of at most countably many compact
sets.  If $M$ is $\sigma$-compact, then every continuous mapping from
$M$ into $N$ is countably uniformly continuous.  If $M$ is separable
and locally compact, then $M$ is $\sigma$-compact.  If $M$ is
$\sigma$-compact, then every open subset of $M$ is $\sigma$-compact as
well, since every open subset of a metric space can be expressed as
the union of a sequence of closed sets.

\subsection{Spaces of continuous mappings}

	We shall write $C(M, N)$ for the space of continuous mappings
from $M$ to $N$.  A mapping $f : M \to N$ is said to be
\emph{bounded}\index{bounded mapping from one metric space to another}
if its image is contained in a bounded subset of $N$, and we write
$C_b(M, N)$ for the space of bounded continuous mappings from $M$ to
$N$.  Note that $C(M, N) = C_b(M, N)$ if $M$ is compact.  If $f_1$, $f_2$
are two bounded continuous mappings from $M$ to $N$, then we set
\begin{equation}
	\theta(f_1, f_2) = \sup \{\rho(f_1(x), f_2(x)) : x \in M \}.
\end{equation}
It is not difficult to verify that this defines a metric on $C_b(M, N)$.
The constant mappings from $M$ to $N$, which take all of $M$ to a single
point in $N$, are obviously bounded and continuous, and in this way
we get a natural embedding of $N$ in $C_b(M, N)$.  This is an isometric
embedding, which is to say that the given metric $\rho(\cdot, \cdot)$
on $N$ agrees with the metric $\theta(\cdot, \cdot)$ applied to the 
constant mappings.

	If $\{f_j\}_{j=1}^\infty$ is a sequence of mappings from $M$
to $N$, and if $f$ is another mapping from $M$ to $N$, then we say
that $\{f_j\}_{j=1}^\infty$ \emph{converges pointwise}\index{pointwise
convergence for a sequence of mappings between two metric spaces} if
\begin{equation}
	\lim_{j \to \infty} f_j(x) = f(x)  \quad\hbox{in } N
\end{equation}
for all $x$ in $M$.  We say that $\{f_j\}_{j=1}^\infty$ converges to
$f$ \emph{uniformly}\index{uniform convergence of a sequence of
mappings between two metric spaces} if for every $\epsilon > 0$
there is an integer $L$ so that
\begin{equation}
	\rho(f_j(x), f(x)) < \epsilon \quad\hbox{for all } j \ge L.
\end{equation}
Thus uniform convergence implies pointwise convergence.

	A well-known result states that if $\{f_j\}_{j=1}^\infty$ is a
sequence of mappings from $M$ to $N$ which converges uniformly to a
mapping $f$ from $M$ to $N$, and if each $f_j$ is continuous, then $f$
is continuous as well.  This does not work in general under the
assumption of pointwise convergence instead of uniform convergence.
Note that a sequence $\{f_j\}_{j=1}^\infty$ in $C_b(M, N)$ converges
to a function $f$ in $C_b(M, N)$ with respect to the metric $\theta$
defined above if and only if $\{f_j\}_{j=1}^\infty$ converges
uniformly to $f$.  If $\{f_j\}_{j=1}^\infty$ converges uniformly to
$f$ and each $f_j$ is uniformly continuous, then $f$ is uniformly
continuous as well.

	If $N$ is complete as a metric space, then $C_b(M, N)$ is
complete as a metric space.  In other words, if $\{f_j\}_{j=1}^\infty$
is a Cauchy sequence in $C_b(M, N)$ with respect to the metric
$\theta$ and $N$ is complete as metric space, then there is a function
$f$ in $C_b(M, N)$ to which $\{f_j\}_{j=1}^\infty$ converges
uniformly.  To show this, one can first use the assumption that
$\{f_j\}_{j=1}^\infty$ is a Cauchy sequence in $C_b(M, N)$ to obtain
that $\{f_j(x)\}_{j=1}^\infty$ is a Cauchy sequence in $N$ for all $x$
in $M$.  Hence $\{f_j\}_{j=1}^\infty$ converges pointwise to some
mapping $f$ from $M$ to $N$, since $N$ is complete.  Using the
assumption that $\{f_j\}_{j=1}^\infty$ is a Cauchy sequence with
respect to the metric $\theta$, one can check that
$\{f_j\}_{j=1}^\infty$ converges to $f$ uniformly.  As a consequence,
$f$ is continuous, and one can also verify that $f$ is bounded in this
case.

	Let us write $C(M)$ for the space of continuous real-valued
functions on $M$, and $C_b(M)$ for the space of bounded continuous
real-valued functions on $M$.  Thus $C(M)$ and $C_b(M)$ are vector
spaces over the real numbers, and on $C_b(M)$ we have the norm
\begin{equation}
	\|f\|_{sup} = \sup \{ |f(x)| : x \in M \}.
\end{equation}
That is, $\|f\|_{sup}$ is a nonnegative real number which is equal to
$0$ exactly when $f$ is identically equal to $0$, $\|t f\|_{sup} = |t|
\, \|f\|_{sup}$ for all real numbers $t$ and all $f$ in $C_b(M)$, and
the triangle inequality for norms holds, namely
\begin{equation}
	\|f_1 + f_2\|_{sup} \le \|f_1\|_{sup} + \|f_2\|_{sup}
\end{equation}
for all $f_1$, $f_2$ in $C_b(M)$.  These properties are not difficult
to verify.  The metric $\theta(\cdot, \cdot)$ defined in general above
is the same as the metric associated to the norm $\|\cdot\|_{sup}$
on $C_b(M)$, which is to say that 
\begin{equation}
	\theta(f_1, f_2) = \|f_1 - f_2\|_{sup}
\end{equation}
for all $f_1$, $f_2$ in $C_b(M)$.  The completeness of $C_b(M)$ as
a metric space implies that $C_b(M)$ equipped with the norm
$\|\cdot\|_{sup}$ is a \emph{Banach space}.\index{Banach space $C_b(M)$}
In fact it is a \emph{Banach algebra} with respect to the additional
operation of multiplication of functions, which means that
\begin{equation}
	\|f_1 \, f_2\|_{sup} \le \|f_1\|_{sup} \, \|f_2\|_{sup}
\end{equation}
for all $f_1$, $f_2$ in $C_b(M)$.

	There is a natural mapping from $M$ to $C(M)$, in which a
point $p$ in $M$ is associated to the continuous function $f_p(x) =
d(x,p)$.  If $M$ is bounded, then this defines a mapping from $M$ into
$C_b(M)$, and one can show that this mapping is an isometric embedding
with respect to the metric $\theta(\cdot, \cdot)$ on $C_b(M)$.  In
other words,
\begin{equation}
	\sup \{|f_p(x) - f_q(x)| : x \in M \} = d(p,q)
\end{equation}
for all $p$, $q$ in $M$, as one can check.  If $M$ is unbounded, then
one can adjust this as follows.  Fix an element $w$ of $M$, which will
serve as a base point.  To each point $p$ in $M$ we can associate the
function $f_p(x) - f_w(x)$.  It is not hard to see that this is a
bounded continuous function on $M$, and that the distance between $f_p
- f_w$ and $f_q - f_w$ is equal to $d(p, q)$, so that we again get an
isometric embedding of $M$ into $C_b(M)$.

	Let us write $UC(M,N)$ for the space of uniformly continuous
mappings from $M$ to $N$, and $UC_b(M,N)$ for the space of bounded
uniformly continuous mappings from $M$ to $N$.  As before, for $N =
{\bf R}$ we may simply write $UC(M)$, $UC_b(M)$, respectively.  Of
course 
\begin{equation}
	C(M,N) = C_b(M,N) = UC(M,N) = UC_b(M,N)
\end{equation}
when $M$ is compact.  If $M$ is totally bounded, then $UC(M,N) =
UC_b(M,N)$, since $f(M)$ is then totally bounded and hence bounded in
$N$ for every uniformly continuous mapping $f$ from $M$ to $N$.

	If $M$ is totally bounded and $N$ is separable, then $UC(M,N)$
is also separable as a metric space equipped with the supremum metric
$\theta(f_1, f_2)$.  In particular, $C(M,N)$ is separable if $M$ is
compact and $N$ is separable.  Let us indicate some of the ingredients
in the proof of this.  Let $\epsilon > 0$ be given.  For each positive
integer $l$, let $\mathcal{UC}_{\epsilon, l}(M,N)$ be the set of
functions $f$ in $UC(M,N)$ such that
\begin{equation}
	\rho(f(x),f(y)) < \epsilon 
		\quad\hbox{when}\quad x, y \in M, d(x, y) < \frac{1}{l}.
\end{equation}
Thus 
\begin{equation}
	UC(M,N) = \bigcup_{l=1}^\infty \mathcal{UC}_{\epsilon, l}(M,N),
\end{equation}
by uniform continuity.  Since $M$ is totally bounded, for each
positive integer $l$ there is a finite subset $F_l$ of $M$ such that
\begin{equation}
	M = \bigcup_{z \in F} B(z, 1/l).
\end{equation}
If $f_1$, $f_2$ are elements of $\mathcal{UC}_{\epsilon, l}(M,N)$
such that 
\begin{equation}
	\rho(f_1(z), f_2(z)) < \epsilon \quad\hbox{when}\quad z \in F_l,
\end{equation}
then
\begin{equation}
	\rho(f_1(w), f_2(w)) < 3 \, \epsilon
					\quad\hbox{when}\quad w \in M.
\end{equation}
Thus an element of $\mathcal{UC}_{\epsilon, l}(M,N)$ is approximately
determined by its values on $F_l$, and one can use this and the
separability of $N$ to get the separability of $UC(M,N)$.

\subsection{Upper and lower semicontinuous functions}

	A real-valued function $f$ on $M$ is said to be \emph{upper
semicontinuous}\index{upper semicontinuous function on a metric space}
if for every $x$ in $M$ and every $\epsilon > 0$ there is a $\delta > 0$
such that
\begin{equation}
\label{def of upper semicontinuity}
	f(y) < f(x) + \epsilon 
	   \quad\hbox{for all } y \in M \hbox{ such that } d(x,y) < \delta.
\end{equation}
A real-valued function $f$ on $M$ is said to be \emph{lower
semicontinuous}\index{lower semicontinuous function on a metric space}
if for every $x$ in $M$ and every $\epsilon > 0$ there is a $\delta >
0$ so that
\begin{equation}
\label{def of lower semicontinuity}
	f(y) > f(x) - \epsilon
      	   \quad\hbox{for all } y \in M \hbox{ such that } d(x,y) < \delta.
\end{equation}
Clearly a real-valued function on $M$ is continuous if and only if it
is both upper and lower semicontinuous.

	Alternatively, $f$ is upper semi-continuous if and only if for
each real number $a$ the set $\{x \in M : f(x) < a\}$ is an open
subset of $M$.  This is equivalent to saying that $\{x \in M : f(x)
\ge a\}$ is a closed subset of $M$ for all real numbers $a$.
Similarly, $f$ is lower semicontinuous if for each real number $a$ the
set $\{x \in M : f(x) > a\}$ is an open subset of $M$, and this is
equivalent to saying that $\{x \in M : f(x) \le a\}$ is a closed
subset of $M$ for each real number $a$.

	In terms of sequences, $f$ is upper semicontinuous if and only
if for every sequence $\{x_j\}_{j=1}^\infty$ in $M$ which converges to
some point $x$ in $M$ one has that
\begin{equation}
	\limsup_{j \to \infty} f(x_j) \le f(x),
\end{equation}
and $f$ is lower semicontinuous if and only if 
\begin{equation}
	f(x) \le \liminf_{j \to \infty} f(x_j)
\end{equation}
for every sequence $\{x_j\}_{j=1}^\infty$ in $M$ which converges to a
point $x$ in $M$.  These statements are not too hard to check, just
from the definitions, in much the same way as for their counterparts
for continuous functions.

	A basic result about semicontinuous functions is that if $M$
is compact and $f$ is an upper semicontinuous function on $M$, then
there is a point $p$ in $M$ such that $f(x) \le f(p)$ for all $x$ in
$M$.  If $f$ is a lower semicontinuous function on $M$ and $M$ is
compact, then there is a $q$ in $M$ such that $f(q) \le f(y)$ for all
$y$ in $M$.  Just as for real-valued continuous functions, these
assertions can be verified either using open coverings of $M$
or sequences and subsequences.

	If $f_1$, $f_2$ are upper semicontinuous functions on $M$,
then the sum $f_1 + f_2$ is upper semicontinuous, while if $f_1$,
$f_2$ are both lower semicontinuous, then $f_1 + f_2$ is lower
semicontinuous.  If $f$ is an upper semicontinuous function on $M$ and
$a$ is a nonengative real number, then $a \, f$ is an upper
semicontinuous function on $M$, and if instead $f$ is a lower
semicontinuous function on $M$ and $a$ is a nonnegative real number,
then $a \, f$ is lower semicontinuous.  A real valued function $f$ on
$M$ is upper semicontinuous if and only if $-f$ is lower
semicontinuous.  

	The product of two nonnegative upper semicontinuous functions
is upper semicontinuous, and the product of two nonnegative lower
semicontinuous functions is lower semicontinuous.  More generally,
the product of a nonnegative continuous function with a function
which is upper or lower semicontinuous is also upper or lower
semicontinuous, respectively.

	Suppose that $\phi(t)$ is a real-valued function on the real
line which is monotone increasing.  As is well known, the left and
right hand limits automatically exist at every point in the real line
in this case, and satisfy
\begin{equation}
\label{lim_{x to t-} phi(x) le phi(t) le lim_{x to t+} phi(x)}
	\lim_{x \to t-} \phi(x) \le \phi(t) \le \lim_{x \to t+} \phi(x).
\end{equation}
The function $\phi$ is continuous on ${\bf R}$ if and only if the
inequalities in (\ref{lim_{x to t-} phi(x) le phi(t) le lim_{x to t+}
phi(x)}) are equalities for all $t \in {\bf R}$.  One can check that
$\phi$ is upper semicontinuous on ${\bf R}$ if and only if the second
inequality in (\ref{lim_{x to t-} phi(x) le phi(t) le lim_{x to t+}
phi(x)}) is an equality for all $t \in {\bf R}$, and $\phi$ is lower
semicontinuous on ${\bf R}$ if and only if the first inequality in
(\ref{lim_{x to t-} phi(x) le phi(t) le lim_{x to t+} phi(x)}) is an
equality for all $t \in {\bf R}$.

	If $f$ is an upper semicontinuous real-valued function on $M$
and $\phi$ is a monotone increase real-valued function on ${\bf R}$
which is upper semicontinuous, then $\phi \circ f$ is an upper
semicontinuous function on $M$.  Similarly, if $f$ is a lower
semicontinuous real-valued function on ${\bf R}$ and $\phi$ is a
monotone increasing real-valued function on ${\bf R}$ which is lower
semicontinuous, then $\phi \circ f$ is a lower semicontinuous function
on $M$.

\section{Connectedness}
\label{section on connectedness}

	Let $(M, d(x,y))$ be a metric space.  Two subsets $A$, $B$
of $M$ are said to be \emph{separated}\index{separated subsets of
a metric space} if
\begin{equation}
\label{def of separatedness for subsets of a metric space}
	\overline{A} \cap B = \emptyset, 
		\quad A \cap \overline{B} = \emptyset.
\end{equation}
A subset $E$ of $M$ is said to be \emph{connected}\index{connected
subsets of a metric space} if it cannot be written as the union of
two separated sets.

	If $E$, $Y$ are subsets of $M$ with $Y$ nonempty and $E
\subseteq Y$, then $E$ is connected as a subset of $M$ if and only if
$E$ is connected as a subset of $Y$, viewed as a metric itself, with
the restriction of the metric from $M$.  This is not difficult to
verify, just from the definition.  In particular, $E$ is connected as
a subset of $M$ if and only if $E$ is connected as a subset of itself,
viewed as a metric space, with the restriction of the metric from $M$.

	The empty set and sets with only one element are automatically
connected.  It is well-known that the connected subsets of the real
line are the empty set and the subsets which are intervals, which may
be open, closed, or mixed, and which may be bounded or unbounded.
This includes the real line itself.

	The connectedness of $M$ itself is equivalent to saying that
$M$ cannot be written as the disjoint union of two nonempty open
subsets, or that $M$ cannot be written as the disjoint union of two
nonempty closed subsets.

	The closure of a connected set is always connected.  The union
of two connected sets is connected as long as they have nonempty
intersection.  More generally, the union of two nonempty connected
sets is connected if and only if they are not separated.

	Let $E$ be a subset of $M$, and suppose that for every pair
of points $p$, $q$ in $E$ there is a subset $E_{p, q}$ of $E$
which is connected and which contains $p$, $q$.  Then $E$ is
connected.

	A set which contains no connected subset with more than one
element is said to be \emph{totally disconnected}.\index{totally
disconnected subsets of a metric space}  This terminology will be
convenient, but it has the feature that the empty set and sets with
one element are called connected and totally disconnected.

	In general, we can define a relation $\sim$ on a metric space
$M$ by saying that $x \sim y$ when $x$ and $y$ are contained in a
connected subset of $M$.  This relation is clearly reflexive and
symmetric, and it is also transitive, because of the earlier remark
about the union of two connected sets being connected when they
contain a common element.

	Thus $\sim$ defines an equivalence relation on $M$.  The
corresponding equivalence classes are called the \emph{connected
components}\index{connected components of a metric space} of $M$.  It
is not hard to see that the connected components of $M$ are indeed
connected subsets of $M$, and in fact they are the maximal connected
subsets of $M$.  The connected components of $M$ are also closed sets,
since the closure of a connected set is always connected.

	If $(M_1, d_1(x,y))$ and $(M_2, d_2(u,v))$ are two metric
spaces and $f$ is a continuous mapping from $M_1$ to $M_2$, and if $E$
is a connected subset of $M_1$, then $f(E)$ is a connected subset of
$M_2$.  As a special case of this, suppose that $a$, $b$ are real
numbers with $a \le b$, and that $p(t)$ is a continuous mapping from
$[a, b]$ into $M$.  We call $p(t)$ a \emph{path}\index{path in a metric
space} in $M$.  The image of a path is always a connected set, since
it is the image of a connected set under a continuous mapping.

	A subset $E$ of $M$ is said to be \emph{pathwise
connected}\index{pathwise connected subsets of a metric space} if for
each $x$, $y$ in $E$ there is a continuous path $p(t)$ defined on an
interval $[a, b]$ and taking values in $E$ such that $p(a) = x$
and $p(b) = y$.  Because the image of the path is a connected set,
it follows that every pair of points in $E$ is contained in a connected
subset of $E$, so that $E$ is itself connected.

	Let us define a relation $\approx$ on $M$ by saying that $x
\approx y$ if there is a path in $M$ from $x$ to $y$, i.e., if there
is a closed interval $[a, b]$ in the real line and a continuous
mapping $p(t)$ from $[a, b]$ into $M$ such that $p(a) = x$, $p(b) =
y$.  Note that $x \approx y$ implies $x \sim y$.  One can check that
$\approx$ defines an equivalence relation on $M$.  The equivalence
classes associated to the relation $\approx$ are pathwise connected
subsets of $M$, and they are called the pathwise-connected
components\index{pathwise-connected components of a metric space} of
$M$.  The pathwise-connected components of $M$ are the maximal
pathwise-connected subsets of $M$, and they are contained in the
ordinary connected components of $M$.

	Connected subsets of the real line are clearly pathwise
connected.  Now fix a positive integer $n$, and suppose that $M$ is an
open subset of ${\bf R}^n$, equipped with the standard metric.  It is
not difficult to show that the pathwise connected components of $M$
are then open subsets of $M$.  In fact the pathwise-connected
components are the same as the connected components, and $M$ is
pathwise-connected if it is connected.  This works more generally
under local connectedness conditions for $M$, which hold in a very
simple way for open subsets of ${\bf R}^n$, just using paths along
line segments locally.

	Fix a positive real number $\epsilon$.  By an
\emph{$\epsilon$-chain}\index{chain@$\epsilon$-chain in a metric
space} in $M$ we mean a finite sequence $z_1, \ldots, z_k$ of points
in $M$ such that $d(z_j, z_{j+1}) < \epsilon$ when $1 \le j \le k-1$.
This condition is considered to hold automatically when $k = 1$, so
that a single point defines an $\epsilon$-chain.  A subset $E$ of $M$
is said to be \emph{$\epsilon$-connected} if for every pair of points
$x$, $y$ in $E$ there is an $\epsilon$-chain $z_1, \ldots, z_k$ in $E$
such that $z_1 = x$ and $z_k = y$.

	Let us again assume that $E = M$ for simplicity.  Define a
relation $\simeq_\epsilon$ on $M$ by saying that $x \simeq_\epsilon y$
if there is an $\epsilon$-chain $z_1, \ldots, z_k$ in $M$ such that
$z_1 = x$ and $z_k = y$.  It is easy to see that $\simeq_\epsilon$
defines an equivalence relation on $M$, and that the corresponding
equivalence classes are the maximal $\epsilon$-connected subsets of
$M$.  By construction, two points in different equivalence classes
associated to $\simeq_\epsilon$ have distance at least $\epsilon$ from
each other.  Using this, one can show that if $M$ is connected in the
ordinary sense, then $M$ is $\epsilon$-connected for all $\epsilon >
0$.

	Conversely, if $M$ is compact and $\epsilon$-connected for
each $\epsilon > 0$, then $M$ is connected.  Note that the set of
rational numbers, equipped with the usual metric, is
$\epsilon$-connected for each $\epsilon > 0$, but not connected in the
usual sense.

\chapter{Some measure functionals and dimensions}
\label{chapter on some measure functionals and dimensions}

	If $M$ is a set, then by a \emph{measure functional on subsets
of $M$}\index{measure functional on subsets of a fixed set} we mean a
function $\mu$ which takes subsets of $M$ to nonnegative real numbers
or $+\infty$ and which satisfies $\mu(\emptyset) = 0$ and $\mu(A) \le
\mu(B)$ whenever $A, B \subseteq M$ and $A \subseteq B$.  In many
cases the measure functional satisfies additional nice properties, at
least for some reaonably-large class of subsets of $M$.

\section{Basic notions}
\label{some basic notions concerning Hausdorff measures and content}

	Let $(M, d(x,y))$ be a metric space, and let $E$ be a subset
of $M$.  Define $\mathcal{U}_{con}(E)$ to be the collection of finite or
countably-infinite families $\{A_i\}_i$ of subsets of $M$ such that
$E$ is contained in the union $\bigcup_i A_i$.  For each nonnegative
real number $\alpha$, define the $\alpha$-dimensional Hausdorff
content\index{Hausdorff content of a subset of a metric space} of $E$,
denoted $H^\alpha_{con}(E)$, to be the infimum of the sums
\begin{equation}
\label{sum_i (diam A_i)^alpha}
	\sum_i (\diam A_i)^\alpha
\end{equation}
over all coverings $\{A_i\}_i$ of $E$ in $\mathcal{U}_{con}(E)$.  If
some $A_i$ is unbounded, so that $\diam A_i = +\infty$, then we
interpret $(\diam A_i)^\alpha$ to be $+\infty$ as well, for all
$\alpha \ge 0$.  If $A_i$ is bounded and $\alpha = 0$, then we
interpret $(\diam A_i)^\alpha$ as being $1$ when $A_i$ is not empty,
and we interpret it as being $0$ if $A_i$ is the empty set.  In
particular, the $\alpha$-dimensional Hausdorff content of the empty
set is taken to be $0$.  Technically, one might also say that for the
empty set the empty covering is admissible, and that the corresponding
empty sum is automatically equal to $0$.

	A simple estimate is 
\begin{equation}
\label{H^alpha_{con}(E) le (diam E)^alpha}
	H^\alpha_{con}(E) \le (\diam E)^\alpha.
\end{equation}
For $\alpha = 0$ we have that $H^0_{con}(E)$ is equal to $0$ when $E$
is the empty set, to $1$ when $E$ is nonempty and bounded, and to
$+\infty$ when $E$ is unbounded.

	Now suppose that $\delta$ is an extended real number such that
$0 < \delta \le \infty$, and that $\alpha$ is a nonnegative real
number.  If $E$ is a subset of $M$, we define $\mathcal{U}_\delta(E)$
to be the collection of finite or countably-infinite families
$\{A_i\}_i$ of subsets of $M$ such that $E$ is contained in the union
$\bigcup_i A_i$ and $\diam A_i < \delta$ for all $i$.  We then define
$H^\alpha_\delta(E)$ to be the infimum of the same sum (\ref{sum_i
(diam A_i)^alpha}) over all $\{A_i\}_i$ in $\mathcal{U}_\delta(E)$
assuming that $\mathcal{U}_\delta(E)$ is nonempty, and otherwise we
define $H^\alpha_\delta(E)$ to be $+\infty$.  If $M$ is separable,
then one can check that $\mathcal{U}_\delta(E)$ is always nonempty.
In that situation $M$ can be covered by a family of at most countably
many balls of radius $r$ for any fixed $r > 0$, for instance.  In
fact, $\mathcal{U}_{\delta}(E)$ is nonempty for all $\delta > 0$ if
and only if $E$ contains a countable dense subset.  By contrast,
$\mathcal{U}_\infty(E)$ is automatically nonempty, because for any
fixed $p$ in $M$, $E$ can be covered by the family of balls $B(p, n)$,
where $n$ runs over all positive integers.

	It is easy to see that
\begin{equation}
\label{H^alpha_{delta_1}(E) ge H^alpha_{delta_2}(E)}
	H^\alpha_{\delta_1}(E) \ge H^\alpha_{\delta_2}(E)
\end{equation}
whenever $0 < \delta_1 \le \delta_2 \le +\infty$.  This is because
\begin{equation}
	\mathcal{U}_{\delta_1}(E) \subseteq \mathcal{U}_{\delta_2}(E),
\end{equation}
while the sums involved in the definitions of
$H^\alpha_{\delta_1}(E)$, $H^\alpha_{\delta_2}(E)$ are identical.  In
other words, $H^\alpha_{\delta_1}(E)$ and $H^\alpha_{\delta_2}(E)$ are
defined as the infima of the same expressions, but with different
collections of coverings, with the coverings for
$H^\alpha_{\delta_1}(E)$ being more restricted than the coverings for
$H^\alpha_{\delta_2}(E)$.  Notice also that
\begin{equation}
\label{H^alpha_infty(E) = H^alpha_{con}(E)}
	H^\alpha_\infty(E) = H^\alpha_{con}(E).
\end{equation}
Again the sums being used are the same, and the coverings of $E$ in
$\mathcal{U}_\infty(E)$ are included in $\mathcal{U}_{con}(E)$.  The
converse does not quite hold, but if $\{A_i\}_i$ lies in
$\mathcal{U}_{con}(E)$ and $\sum_i (\diam A_i)^\alpha < \infty$, then
$\diam A_i < \infty$ for all $i$, and $\{A_i\}_i$ lies in
$\mathcal{U}_\infty(E)$.  It is easy to use this to obtain
(\ref{H^alpha_infty(E) = H^alpha_{con}(E)}).

	When $\alpha = 0$ we are reduced to the following.  If $E$ is
the empty set, then $H^0_\delta(E) = 0$, and otherwise $H^0_\delta(E)$
is equal to the smallest number of sets of diameter less than $\delta$
needed to cover $E$, where this is interpreted as being $+\infty$ if
no finite covering of $E$ by sets with diameter less than $\delta$
exists.

	The \emph{$\alpha$-dimensional Hausdorff
measure}\index{Hausdorff measure of a subset of a metric space} of $E$
is denoted $H^\alpha(E)$ and defined by
\begin{equation}
	H^\alpha(E) = \sup \{H^\alpha_\delta(E) : \delta > 0\}.
\end{equation}
Because of the monotonicity property (\ref{H^alpha_{delta_1}(E) ge
H^alpha_{delta_2}(E)}), this is the same as the limit as $\delta \to
0$ of $H^\alpha_\delta(E)$.  By definition,
\begin{equation}
\label{H^alpha_{con}(E) le H^alpha_delta(E) le H^alpha(E)}
	H^\alpha_{con}(E) \le H^\alpha_\delta(E) \le H^\alpha(E)
\end{equation}
for $0 < \delta \le \infty$.

	If $\alpha = 0$, then one can check that Hausdorff measure
reduces to counting measure.\index{counting measure} That is to say,
$H^0(E)$ is equal to the number of elements of $E$ when $E$ is finite,
and to $+\infty$ when $E$ is infinite.

	Observe that $H^\alpha_{con}$, $H^\alpha_\delta$, and
$H^\alpha$ are all measure functionals on subsets of $M$ in the sense
described at the beginning of the chapter.  Also, for each $\alpha \ge
0$, $H^\alpha_{con}(E) = 0$ implies that $H^\alpha_\delta(E) = 0$ for
all $\delta > 0$, and hence that $H^\alpha(E) = 0$.  Thus for a fixed
$\alpha$, the measure functionals $H^\alpha_{con}$, $H^\alpha_\delta$, and
$H^\alpha$ vanish on the same subsets of $M$.

	We can make analogous definitions using finite coverings.
Namely, for a subset $E$ of $M$, let $\mathcal{U}^f_{con}(E)$
denote the collection of finite families $\{A_i\}_i$ of subsets
of $M$ such that $E$ is contained in $\bigcup_i A_i$.  For each
nonnegative real number $\alpha$, define $HF^\alpha_{con}(E)$
to be the infimum of
\begin{equation}
\label{sum_i (diam A_i)^alpha, rewritten}
	\sum_i (\diam A_i)^\alpha
\end{equation}
over all finite coverings $\{A_i\}_i$ of $E$ in $M$.  Thus
we have the simple inequality
\begin{equation}
\label{HF^alpha_{con}(E) le (diam E)^alpha}
	HF^\alpha_{con}(E) \le (\diam E)^\alpha.
\end{equation}

	For $0 < \delta \le \infty$, define $\mathcal{U}^f_\delta(E)$
to be the collection of finite families $\{A_i\}_i$ of subsets of $M$
such that $E \subseteq \bigcup_i A_i$ and $\diam A_i < \delta$ for
each $i$.  We define $HF^\alpha_\delta(E)$ to be the infimum of
(\ref{sum_i (diam A_i)^alpha, rewritten}) over all $\{A_i\}_i$ in
$\mathcal{U}^f_\delta(E)$ if $\mathcal{U}^f_\delta(E)$ is not empty,
and otherwise to be $+\infty$.  Observe that $E$ is totally bounded if
and only if $\mathcal{U}^f_\delta(E)$ is finite for all $\delta > 0$,
which is also equivalent to saying that $HF^\alpha_\delta(E)$ is
finite for all $\delta > 0$.  If $E$ is unbounded, then
$HF^\alpha_{con}(E)$ and $HF^\alpha_\delta(E)$ are equal to $+\infty$
for all $\alpha$, $\delta$.

	As before, $\mathcal{U}^f_{\delta_1}(E) \subseteq
\mathcal{U}^f_{\delta_2}(E)$ when $\delta_1 \le \delta_2$, and hence
\begin{equation}
	HF^\alpha_{\delta_2}(E) \le HF^\alpha_{\delta_1}(E)
\end{equation}
in this case.  Also,
\begin{equation}
	HF^\alpha_\infty(E) = HF^\alpha_{con}(E),
\end{equation}
almost by definition.  We define $HF^\alpha(E)$ by
\begin{equation}
\label{def of HF^alpha(E)}
	HF^\alpha(E) = \sup \{HF^\alpha_\delta(E) : \delta > 0\},
\end{equation}
which is the same as $\lim_{\delta \to 0} HF^\alpha_\delta(E)$,
by monotonicity.  Thus
\begin{equation}
\label{HF^alpha_{con}(E) le HF^alpha_delta(E) le HF^alpha(E)}
	HF^\alpha_{con}(E) \le HF^\alpha_\delta(E) \le HF^\alpha(E)
\end{equation}
for $0 < \delta \le \infty$.

	When $\alpha = 0$, these quantities are the same as their
counterparts described before, i.e., $HF^0_{con}(E) = H^0_{con}(E)$,
$HF^0_\delta(E) = H^0_\delta(E)$ for $0 < \delta \le \infty$, and
$HF^0(E) = H^0(E)$.  For each $\alpha$ and $\delta$,
$HF^\alpha_{con}$, $HF^\alpha_\delta$, and $HF^\alpha$ define
measure functionals on subsets of $M$.  If $\alpha$ is fixed, then
$HF^\alpha_{con}$, $HF^\alpha_\delta$, and $HF^\alpha$ are equal to
$0$ on the same subsets of $M$.

	Because $\mathcal{U}^f(E) \subseteq \mathcal{U}(E)$ and
$\mathcal{U}^f_\delta(E) \subseteq \mathcal{U}_\delta(E)$ for all
$\delta$, we have that 
\begin{equation}
	H^\alpha_{con}(E) \le HF^\alpha_{con}(E),
	\enspace H^\alpha_\delta(E) \le HF^\alpha_\delta(E),
	\enspace H^\alpha(E) \le HF^\alpha(E).
\end{equation}
We shall see cases when these inequalities are strict, and cases where
equality holds.

	There are a number of variations of these definitions, in
which one might restrict the sets $A_i$ used in the coverings, or use
other measurements of the sizes of the $A_i$'s, etc.  These matters
are related to well-known work of Carath\'eodory, Jordan, and
Minkowski, for instance.

	Note that if $Y$ is a nonempty subset of $M$ and $E$ is a
subset of $Y$, then $H^\alpha_{con}(E)$, $HF^\alpha_{con}(E)$,
$H^\alpha_\delta(E)$, $HF^\alpha_\delta(E)$, $H^\alpha(E)$, and
$HF^\alpha(E)$ are the same for $E$ as a subset of $M$ as they are for
$E$ as a subset of $Y$, viewed as metric space itself, with the
restriction of the metric $d(x, y)$ from $M$.

\section{Coverings by open and closed subsets}
\label{coverings by open and closed subsets}

	Let $(M, d(x,y))$ continue to be some metric space.

\beginlemma
\label{diam overline{A} = diam A}
If $A$ is a subset of $M$, then $\diam \overline{A} = \diam A$.
\end{lemma}

	This is easy to check.  As a result, for the definitions of
the various measures of a set $E$ in Section \ref{some basic notions
concerning Hausdorff measures and content}, one may as well use
coverings by families $\{A_i\}_i$ of \emph{closed} subsets of $M$.  In
other words, if one starts with a covering $\{A_i\}_i$ of some subset
$E$ of $M$ by arbitrary subsets of $M$, then one can replace this with
the family $\{\overline{A}_i\}_i$, and the latter is still be a
covering of $E$.  This new covering also satisfies the same side
conditions as the initial covering used in Section \ref{some basic
notions concerning Hausdorff measures and content}, and does not
change the value of sums of the form $\sum_i (\diam A_i)^\alpha$.

\beginlemma
\label{HF values of E same for overline(E)}
If $E$ is a subset of $M$, then $HF^\alpha_{con}(E) =
HF^\alpha_{con}(\overline{E})$, $HF^\alpha_\delta(E) =
HF^\alpha_\delta(E)$, and $HF^\alpha(E) = HF^\alpha(\overline{E})$ for
all $\alpha \ge 0$ and all $0 < \delta \le \infty$.
\end{lemma}

	This follows easily from the remarks preceding the statement
of the lemma.

\beginlemma
\label{lemma about A_r}
Let $A$ be a subset of $M$, and let $r$ be a positive real number.
Define $A_r$ by
\begin{equation}
\label{def of A_r}
	A_r = \bigcup \{B(a, r) : a \in A\}.
\end{equation}
Then $A_r$ is an open subset of $M$ such that $A \subseteq A_r$ and
$\diam A_r \le \diam A + 2 r$.
\end{lemma}

	This is easy to see.  As a consequence, if one restricts
oneself to coverings by open subsets in Section \ref{some basic
notions concerning Hausdorff measures and content}, then the resulting
values of the various measures of a set $E$ are the same.  

\beginlemma
\label{H and HF the same for compact sets}
If $E$ is a compact subset of $M$, then $H^\alpha_{con}(E) =
HF^\alpha_{con}(E)$, $H^\alpha_\delta(E) = HF^\alpha_\delta(E)$, and
$H^\alpha(E) = HF^\alpha(E)$ for all $\alpha \ge 0$ and $0 < \delta
\le \infty$.
\end{lemma}

	Indeed, one can use coverings of $E$ by open sets, and then
reduce to finite subcoverings by compactness.

\beginlemma
\label{replacing E with an intersection of a sequence of open sets}
Let $E$ be a subset of $M$.  For each $\alpha \ge 0$ there is a subset
$E_1$ of $M$ such that $E \subseteq E_1$, $E_1$ can be written as the
intersection of a sequence of open subsets of $M$, and
$H^\alpha_{con}(E_1) = H^\alpha_{con}(E)$.  The analogous statements
for $H^\alpha_\delta$ and $H^\alpha$ also hold for all $\alpha \ge 0$
and $0 < \delta \le \infty$.
\end{lemma}

	This is not difficult to show.  Fix $\alpha > 0$, and
consider $H^\alpha_{con}(E)$.  If $H^\alpha_{con}(E) = \infty$,
then we can simply take $E_1 = M$.  Otherwise, for each positive
integer $n$, let $\{A_i(n)\}_i$ be a covering of $E$ by open subsets
of $M$ such that 
\begin{equation}
	\sum_i (\diam A_i(n))^\alpha < H^\alpha_{con}(E) + \frac{1}{n}.
\end{equation}
Put
\begin{equation}
	E_1 = \bigcap_{n=1}^\infty \bigcup_i A_i(n).
\end{equation}
Then $E \subseteq E_1$ and $E_1$ is the intersection of a sequence of
open sets by construction, and it is not difficult to see that
$H^\alpha_{con}(E_1) = H^\alpha_{con}(E)$.  Indeed, $H^\alpha_{con}(E)
\le H^\alpha_{con}(E_1)$ because $E \subseteq E_1$, while
$H^\alpha_{con}(E_1) \le H^\alpha_{con}(E)$ because $\{A_i(n)\}_i$ is
a covering of $E_1$ for all $n$.  The argument for
$H^\alpha_\delta(E)$ is quite similar, except that one should also ask
that $\diam A_i(n) < \delta$ for all $i$ and $n$, while for
$H^\alpha(E)$ one should ask that $\diam A_i(n) < 1/n$ for all $i$ and
$n$.

\section{Subadditivity properties}
\label{subadditivity properties}

	Let $M$ be a set.  A measure functional $\mu$ on subsets of $M$ is
said to be \emph{finitely subadditive}\index{finitely subadditive
measure functionals on subsets of a set} if
\begin{equation}
\label{finite subadditivity}
	\mu(E_1 \cup E_2) \le \mu(E_1) + \mu(E_2)
\end{equation}
for all subsets $E_1$, $E_2$ of $M$.  If
\begin{equation}
\label{countable subadditivity}
	\mu\biggl(\sum_j E_j \biggr) \le \sum_j \mu(E_j)
\end{equation}
for any at-most-countable family $\{E_j\}_j$ of subsets of $M$,
then $\mu$ is said to be \emph{countably subadditive}.\index{countably
subadditive measure functionals on subsets of a set}

	For the rest of this section, we assume that $(M, d(x,y))$
is a metric space.

\beginlemma
\label{subaddivity properties of the H and HF measure functionals}
For each $\alpha \ge 0$ and $0 < \delta \le \infty$, $HF^\alpha_{con}$,
$HF^\alpha_\delta$, and $HF^\alpha$ are finitely subadditive, and
$H^\alpha_{con}$, $H^\alpha_\delta$, and $H^\alpha$ are countably
subadditive.
\end{lemma}

	This is an easy consequence of the definitions.  In particular,
if $E$ is a finite subset of $M$ and $\alpha > 0$, then $HF^\alpha(E) = 0$,
and if $E \subseteq M$ is at most countable and $\alpha > 0$ then
$H^\alpha(E) = 0$.

	If $E$ is a bounded subset of ${\bf R}^n$, then it is not
difficult to check that $H^\alpha(E) = HF^\alpha(E) = 0$ for all
$\alpha > 0$.  Using the countable subadditivity of $H^\alpha$,
it follows that $H^\alpha({\bf R}^n) = 0$ for all $\alpha > n$.

	Let us say that two subsets $E_1$, $E_2$ are $\eta$-separated,
where $\eta$ is a positive real number, if $d(x,y) \ge \eta$ for all
$x$ in $E_1$ and all $y$ in $E_2$.

\beginlemma
\label{an additivity property}
Let $E_1$, $E_2$ be a pair of $\eta$-separated subsets of $M$.
For all $\alpha \ge 0$ and $\delta \in (0, \eta]$ we have that
\begin{equation}
	HF^\alpha_\delta(E_1 \cup E_2) 
		\ge HF^\alpha_\delta(E_1) + HF^\alpha_\delta(E_2)
\end{equation}
and
\begin{equation}
	H^\alpha_\delta(E_1 \cup E_2) 
		\ge H^\alpha_\delta(E_1) + H^\alpha_\delta(E_2).
\end{equation}
As a result,
\begin{equation}
	HF^\alpha(E_1 \cup E_2) \ge HF^\alpha(E_1) + HF^\alpha(E_2)
\end{equation}
and
\begin{equation}
	H^\alpha(E_1 \cup E_2) \ge H^\alpha(E_1) + H^\alpha(E_2).
\end{equation}
\end{lemma}

	This is not too difficult to verify.

\section{Dimensions}
\label{section about Hausdorff dimensions, etc.}

	Let $(M, d(x,y))$ be a metric space.

\beginlemma
\label{simple inequalities concerning Hausdorff measures, etc.}
Let $\alpha$, $\beta$ be nonnegative real numbers such that $\alpha
\le \beta$, and let $E$ be a subset of $M$.  For $0 < \delta < \infty$
we have that 
\begin{equation}
	H^\beta_\delta(E) \le \delta^{\beta - \alpha} \, H^\alpha_\delta(E)
\end{equation}
and 
\begin{equation}
	HF^\beta_\delta(E) \le \delta^{\beta - \alpha} \, HF^\alpha_\delta(E).
\end{equation}
In particular, $H^\beta_\delta(E) \le H^\alpha_\delta(E)$ and
$HF^\beta_\delta(E) \le HF^\alpha_\delta(E)$ when $0 < \delta \le 1$.
As a result, $H^\beta(E) \le H^\alpha(E)$, $HF^\beta(E) \le
HF^\alpha(E)$, $H^\alpha_{con}(E) = 0$ implies that $H^\beta_{con}(E)
= 0$, and $HF^\alpha_{con}(E) = 0$ implies that $HF^\beta_{con}(E) =
0$.  Furthermore, if $\alpha < \beta$, then $H^\alpha(E) < \infty$
implies that $H^\beta(E) = 0$, and $HF^\alpha(E) < \infty$ implies
that $HF^\beta(E) = 0$.
\end{lemma}

	This is easy to check, using the fact that $r^\beta \le
\delta^{\beta - \alpha} \, r^\alpha$ when $0 \le r < \delta$.

	The \emph{Hausdorff dimension},\index{Hausdorff dimension of a
subset of a metric space} or \emph{$H$-dimension}, of a subset $E$ of
$M$ is denoted $\dim_H(E)$ and defined to be the infimum of the set
\begin{equation}
	\{\alpha \ge 0 : H^\alpha(E) = 0\}
\end{equation}
when this set is nonempty, and $+\infty$ otherwise.  Similarly, the
\emph{$HF$-dimension} of $E$ is denoted $\dim_{HF}(E)$ and defined to
be the infimum of the set
\begin{equation}
	\{\alpha \ge 0 : HF^\alpha(E) = 0\}
\end{equation}
when this set is nonempty, and $+\infty$ otherwise.  As a variant of
the $HF$-dimension, the \emph{$HF^*$-dimension} of $E$ is denoted
$\dim_{HF^*}(E)$ and defined to be the supremum of the $HF$ dimensions
of $E \cap B(p, n)$ over all positive integers $n$, where $p$ is any
fixed point in $M$.  One can easily check that this does not depend on
the choice of $p$.  Of course the $HF^*$-dimension of $E$ is
automatically equal to the $HF$-dimension of $E$ when $E$ is bounded.
If $E$ is unbounded, then the $HF$-dimension of $E$ is $+\infty$,
while the $HF^*$-dimension of $E$ can be finite.  Note that the empty
set has dimension $0$ in each of these senses.

\beginlemma
\label{dimension inequalities for subsets}
If $E_1$, $E_2$ are subsets of $M$ such that $E_1 \subseteq E_2$,
then the Hausdorff dimension of $E_1$ is less than or equal to the
Hausdorff dimension of $E_2$, the $HF$-dimension of $E_1$ is less than
or equal to the $HF$-dimension of $E_2$, and the $HF^*$-dimension
of $E_1$ is less than or equal to the $HF^*$-dimension of $E_2$.
\end{lemma}

	This is an easy consequence of the definitions.

\beginlemma
\label{H and HF dimensions the same for compact sets}
If $E$ is a compact subset of $M$, then the Hausdorff and
$HF$-dimensions of $E$ are equal.
\end{lemma}

	This follows from the fact that the Hausdorff and $HF$-measures
are equal for compact sets.

\beginlemma
\label{H and HF^* dimensions the same if closed + bounded subsets are compact}
If $E$ is a closed subset of $M$ such that closed and bounded subsets
of $E$ are compact, then the Hausdorff and $HF^*$ dimensions of $E$
are equal.
\end{lemma}

	This is not difficult to verify.

\beginlemma
\label{HF, HF^* dimensions of a union of two sets}
Let $E_1$, $E_2$ be two subsets of $M$.  The $HF$-dimension of $E_1
\cup E_2$ is equal to maximum of the $HF$-dimensions of $E_1$ and
$E_2$, and the $HF^*$-dimension of $E_1 \cup E_2$ is equal to the
maximum of the $HF^*$ dimensions of $E_1$ and $E_2$.
\end{lemma}

	Indeed, one equality can be obtained from Lemma \ref{dimension
inequalities for subsets}, while the other uses the subadditivity of
the $HF$ measures.

\beginlemma
\label{Hausdorff dimension of a countable union}
Let $\{E_j\}_{j=1}^\infty$ be a sequence of subsets of $M$.  The
Hausdorff dimension of $\bigcup_{j = 1}^\infty E_j$ is equal to the
supremum of the Hausdorff dimensions of the $E_j$'s, $j \in {\bf
Z}_+$.
\end{lemma}

	Again, one inequality follows from Lemma \ref{dimension
inequalities for subsets}, and the other can be derived from countable
subaddivity of the Hausdorff measures.

\section{Snowflake transforms}
\label{section about snowflake transforms}

\beginlemma
\label{(r_1 + r_2)^a le r_1^a + r_2^a when 0 < a le 1}
Let $a$ be a real number such that $0 < a \le 1$.  If $r_1$, $r_2$
are nonnegative real numbers, then $(r_1 + r_2)^a \le r_1^a + r_2^a$.
\end{lemma}

	To see this, observe first that
\begin{equation}
	r_1 + r_2 \le \max(r_1, r_2)^{1-a} \, (r_1^a + r_2^a).
\end{equation}
Also, $\max(r_1, r_2)^a \le r_1^a + r_2^a$, and therefore
\begin{equation}
	r_1 + r_2 \le (r_1^a + r_2^a)^{((1-a)/a) + 1} = (r_1^a + r_2^a)^{1/a}.
\end{equation}
This implies the inequality stated in the lemma.

	Let $(M, d(x,y))$ be a metric space, and let $a$ be a real
number such that $0 < a \le 1$.  It is not hard to check that $(M,
d(x,y)^a)$ is then also a metric space, using the lemma for the
triangle inequality.  We call $(M, d(x,y)^a)$ the \emph{snowflake
transform of order $a$}\index{snowflake transform of a metric space}
of the metric space $(M, d(x,y))$.  Observe that the metrics $d(x,y)$
and $d(x,y)^a$ determine the same topology on $M$, which is the same
as saying that the identity mapping on $M$ is a homemorphism as a map
between the metric spaces $(M, d(x,y))$ and $(M, d(x,y)^a)$.  Also,
a subset of $M$ is bounded with respect to $d(x,y)$ if and only if it
is bounded with respect to $d(x,y)^a$.

	If $A$ is a subset of $M$, then the diameter of $A$ with
respect to $d(x,y)^a$ is the same as the diameter of $A$ with respect
to $d(x,y)$ to the $a$th power.  For each nonnegative real number
$\alpha$, one can use this to check that $H^\alpha_{con}(A)$ and
$HF^\alpha_{con}(A)$ with respect to $d(x,y)$ are the same as
$H^{\alpha/ a}_{con}(A)$ and $HF^{\alpha / a}_{con}(A)$ with respect
to $d(x,y)^a$.  Similarly, for each $\delta > 0$, $H^\alpha_\delta(A)$
and $HF^\alpha_\delta(A)$ with respect to $d(x,y)$ are the same as
$H^{\alpha / a}_{\delta^a}(A)$ and $HF^{\alpha / a}_{\delta^a}(A)$
with respect to $d(x,y)^a$.  Taking the supremum over $\delta$, we
obtain that $H^\alpha(A)$ and $HF^\alpha(A)$ with respect to $d(x,y)$
is the same as $H^{\alpha / a}(A)$ and $HF^{\alpha / a}(A)$ with
respect to $d(x,y)^a$.  The Hausdorff, $HF$, and $HF^*$-dimensions of
$A$ with respect to $d(x,y)$ are equal to the Hausdorff, $HF$, and
$HF^*$-dimensions of $A$ with respect to $d(x,y)^a$ divided by $a$.

\chapter{Miscellaneous, 1}
\label{chapter -- miscellaneous, 1}

\section{The real line}
\label{section about the real line}

	If $A$ is any subset of the real line, then there is an
interval $J$ contained in the real line such that $A \subseteq J$ and
$\diam J = \diam A$.  Here one should allow unbounded intervals $J$ to
deal with the case where $A$ is unbounded, and it is convenient to
consider the empty set as a kind of degenerate interval, to be
compatible with some of our conventions.  As a result of this simple
observation, it follows that for subsets of the real line, if one
defines the measure functionals $H^\alpha_{con}$, $H^\alpha_\delta$,
$H^\alpha$, $HF^\alpha_{con}$, $HF^\alpha_\delta$, and $HF^\alpha$ in
the same manner as in Section \ref{some basic notions concerning
Hausdorff measures and content}, except for using coverings by
families of intervals, then one gets the same answers as for the
earlier definitions using coverings by families of general sets.

	Now let us specialize to the case of $\alpha = 1$.

\beginlemma
\label{equalities between measure functionals on R, alpha = 1}
Let $E$ be a subset of the real line.  For $0 < \delta \le \infty$
we have that $H^1_{con}(E) = H^1_\delta(E) = H^1(E)$ and
$HF^1_{con}(E) = HF^1_\delta(E) = HF^1(E)$.
\end{lemma}

	Indeed, if $A$ is any subset of the real line, and $\delta >
0$, then there is a covering of $A$ by intervals with diameter less
than $\delta$ such that the sum of their lengths is equal to the
diameter of $A$.  If $A$ is bounded, then a finite collection of
intervals can be used in the covering.  This permits one to replace
a given covering of $E$ by sets of smaller diameter, while keeping
the sums of the diameters of the sets in the covering fixed.

\beginlemma
\label{H^1 = HF^1 = length for an interval}
If $J = [a, b]$ is a closed and bounded interval in the real line,
then $H^1(J)$ and $HF^1(J)$ are both equal to the length $|b-a|$ of
the interval.
\end{lemma}

	We already know from Lemma \ref{H and HF the same for compact
sets} that $H^1(J) = HF^1(J)$, since $J$ is compact.  It is easy to
see that $HF^1(J) \le |b - a|$, through explicit coverings of $J$.
This follows from (\ref{HF^alpha_{con}(E) le (diam E)^alpha}) and
Lemma \ref{equalities between measure functionals on R, alpha = 1} as
well.  The remaining point is that $HF^1(J) \ge |b-a|$, and we 
leave this as an exercise.

\beginlemma
\label{measure functionals = 0 when alpha > 1, subsets of R}
If $E$ is a bounded subset of ${\bf R}$ and $\alpha > 1$, then
$H^\alpha(E) = HF^\alpha(E) = 0$.  For any subset $E$ of ${\bf R}$ and
$\alpha > 1$, $H^\alpha(E) = 0$.
\end{lemma}

	The first statement can be verified using simple choices
of coverings of $E$, while the second follows from countable
subadditivity of $H^\alpha$.

	As a consequence of these two lemmas, ${\bf R}$ has Hausdorff
and $HF^*$-dimensions equal to $1$, and the Hausdorff, $HF$, and
$HF^*$-dimensions of a closed and bounded interval $J = [a, b]$
with $a < b$ are equal to $1$.

	Let us now explain how a subset of the real line might
be equal to the intersection of a sequence of dense open sets,
and also have Hausdorff dimension $0$.  For each positive integer $n$,
let $\{r_j(n)\}_{j=1}^\infty$ be a sequence of positive real numbers
such that $r_j(n) < 1/(2n)$ for all $n$ and $j$, and
\begin{equation}
	\sum_{j=1}^\infty (2 r_j)^{1/n} < \frac{1}{n}.
\end{equation}
For instance, one can take $r_j(n) = (j + m_n)^{-2n}$ for
sufficiently-large choices of $m_n$.  Next, let $\{x_j\}_{j=1}^\infty$
be an enumeration of the rational numbers.  For each $n$, put
\begin{equation}
	U_n = \bigcup_{j=1}^\infty (x_j - r_j(n), x_j + r_j(n)).
\end{equation}
Thus $U_n$ is a dense open subset of ${\bf R}$, and
\begin{equation}
	H^{1/n}_{1/n}(U_n) \le \sum_{j=1}^\infty (2 r_j)^{1/n} < \frac{1}{n}.
\end{equation}
If $E = \bigcap_{n=1}^\infty U_n$, then $E$ is the intersection of a
sequence of dense open subsets of the real line, and one can check
that $H^\alpha(E) = 0$ for all $\alpha > 0$.

\section{${\bf R}^n$ for general $n$}
\label{section about {bf R}^n for general n}

	Fix a positive integer $n$.  Recall that a subset $E$ of ${\bf
R}^n$ is said to be \emph{convex}\index{convex subsets of ${\bf R}^n$}
if for every pair of points $v$, $w$ in $E$ and every real number $t$
with $0 \le t \le 1$ we have that $t \, v + (1-t) \, w$ also lies in
$E$.  For example, open and closed balls in ${\bf R}^n$ are convex sets.

	If $l$ is a positive integer and $y_1, \ldots, y_l$ are points
in ${\bf R}^n$, then a \emph{convex combination}\index{convex combination
of a finite collection of points in ${\bf R}^n$} is a point in ${\bf R}^n$
of the form
\begin{equation}
	\sum_{j=1}^l \lambda_j \, y_j,
\end{equation}
where $\lambda_1, \ldots, \lambda_l$ are nonnegative real numbers such
that
\begin{equation}
	\sum_{j=1}^l \lambda_j = 1.
\end{equation}
If $A$ is a subset of ${\bf R}^n$, then the \emph{convex
hull}\index{convex hull of a subset of ${\bf R}^n$} is the set
consisting of all points in ${\bf R}^n$ which are convex combinations
of points in $A$.  The convex hull of $A$ is the smallest convex
subset of ${\bf R}^n$ which contains $A$.

	The closure of a convex set is also convex, as one can easily
verify.  The closure of the convex hull of a subset $A$ of ${\bf R}^n$
can be described as the smallest closed convex subset of ${\bf R}^n$
that contains $A$.  Let us mention the classical fact that the
convex hull of a compact subset of ${\bf R}^n$ is compact.  This can
be established by first showing that every element of the convex hull
of a subset $A$ of ${\bf R}^n$ can be written as a convex combination
of less than or equal to $n+1$ elements of $A$.

	If $A$ is a subset of ${\bf R}^n$, then the diameter of the
closure of the convex hull of $A$ is equal to the diameter of $A$.
This simple fact implies that for the measure functionals on ${\bf
R}^n$ described in Section \ref{some basic notions concerning
Hausdorff measures and content}, i.e., the various Hausdorff measures
and contents and so on, one can restrict to coverings by closed convex
sets and get the same values as before.  This works as well in any
normed vector space.  When $n = 1$, this corresponds to the remark in
the previous section that one can restrict to coverings by intervals
and get the same values for the Hausdorff measures and contents, etc.

	Observe that if $A$ is a subset of ${\bf R}^n$, then the
diameter of $A$ is equal to the diameter of the image of $A$ under a
translation or orthogonal transformation on ${\bf R}^n$.  As a result,
the measure functionals described in Section \ref{some basic notions
concerning Hausdorff measures and content} applied to a subset $E$ of
${\bf R}^n$ have the same value as when they are applied to the image
of $E$ under a translation or orthogonal transformation.  A similar
remark works in any metric space, concerning sets which are
isometrically-equivalent; Euclidean spaces happen to have a lot of
such symmetries.

	Now let us specialize to the case where $\alpha = n$ for the
measure functionals in Section \ref{some basic notions concerning
Hausdorff measures and content}.

\beginlemma
\label{estimate for measure in terms of content on R^n, alpha = n}
For each positive integer $n$ there is a positive real number $C(n)$
such that $HF^n(E) \le HF^n_{con}(E)$ and $H^n(E) \le H^n_{con}(E)$
for any subset $E$ of ${\bf R}^n$.
\end{lemma}

	This follows from the geometric fact that if $A$ is a bounded
subset of ${\bf R}^n$ and $\delta$ is a positive real number such that
$\delta < \diam A$, then $A$ can be covered by $\le C(n) \,
\delta^{-n} \, (\diam A)^n$ subsets of ${\bf R}^n$ with diameter
less than $\delta$.

	If $K$ is a compact subset of ${\bf R}^n$ which is the closure
of a region with reasonably-nice boundary, then the classical way to
measure the volume of $K$ is by approximating $K$ by a finite union of
cubes of the same size.  For reasonably-nice compact sets,
approximations of $K$ by unions of cubes contained in $K$ and with
disjoint interiors and approximations by unions of cubes which contain
$K$ as a subset give nearly the same estimate of the volume, and the
volume of $K$ is defined to be the limit of these approximations as
the mesh size tends to $0$.

	Using this, it is not too difficult to show that the
$n$-dimensional Hausdorff measure of $K$ is equal to the classical
volume of $K$ times the $n$-dimensional Hausdorff measure of the unit
cube $Q_0 = \{x \in {\bf R}^n : 0 \le x_j \le 1 \hbox{ for } j = 1,
\ldots, n\}$.  One can also show that $0 < H^n(Q_0) < \infty$.

	More generally, if $m$ is a positive integer less than $n$,
and if $M$ is a compact $C^1$ submanifold of ${\bf R}^n$ perhaps
with boundary, say, then the $m$-dimensional Hausdorff measure of $M$
is the same as the classical $m$-dimensional volume of $M$ times the
$m$-dimensional Hausdorff measure of an $m$-dimensional cube with 
sidelength equal to $1$.

\section{Cantor sets}
\label{section on cantor sets}
\index{Cantor sets}

	Let $\{\alpha_j\}_{j=1}^\infty$ be a sequence of positive real
numbers such that $0 < \alpha_j < 1/2$ for all $j$.  We shall
sometimes denote this sequence simply as $A$.  For each nonnegative
integer $j$ let us define a family $\mathcal{F}_j(A)$ of $2^j$
disjoint closed subintervals of the unit interval $[0, 1]$ as follows.
When $j = 0$ this family consists of only the unit interval $[0, 1]$
itself.  If $\mathcal{F}_j(A)$ has already been defined, and $I$ is an
element of $\mathcal{F}_j(A)$, then we consider the two closed
subintervals of $I$ which have length equal to $\alpha_{j+1}$ times
the length of $I$ and which each contain one of the two endpoints of
$I$.  We let $\mathcal{F}_{j+1}(A)$ be the family of closed
subintervals of intervals in $\mathcal{F}_j(A)$ obtained in this
manner.

	It is easy to see that number of intervals in $\mathcal{F}_j(A)$
is $2^j$, and that they each have length equal to
\begin{equation}
\label{length of interval in mathcal{F}_j(A)}
	\prod_{i=1}^j \alpha_i
\end{equation}
when $j \ge 1$.  Also, the intervals in $\mathcal{F}_j(A)$ are
pairwise disjoint, as one can verify using induction.  Put
\begin{equation}
	E_j(A) = \bigcup_{I \in \mathcal{F}_j(A)} I,
\end{equation}
and observe that
\begin{equation}
	E_{j+1}(A) \subseteq E_j(A)
\end{equation}
for all $j \ge 0$.  Of course each $E_j(A)$ is a nonempty closed
subset of $[0, 1]$.

	Define $E(A)$, the Cantor set associated to the sequence
$A = \{\alpha_j\}_{j=1}^\infty$, by
\begin{equation}
	E(A) = \bigcap_{j=1}^\infty E_j(A).
\end{equation}
Notice that any point in $[0, 1]$ which occurs as an endpoint of an
interval in $\mathcal{F}_j(A)$ for some $j$ also occurs as an endpoint
of an interval in $\mathcal{F}_k(A)$ for all $k \ge j$, as one can
check using induction.  Thus these points all lie in $E(A)$, and it is
not difficult to show that the set of these points is dense in $E(A)$.
In particular, $0$ and $1$ are elements of $E(A)$.

	If $\alpha_j = 1/3$ for all $j$, then $E(A)$ is the classical
middle-thirds Cantor set.  For each $A$, $E(A)$ is a nonempty closed
subset of $[0, 1]$ such that each element of $E(A)$ is a limit point
of $E(A)$, but $E(A)$ does not contain any segment of positive length.
These more general ``Cantor sets'' can be interesting examples for
various questions.

	Suppose now that $A' = \{\alpha'_j\}_{j=1}^\infty$ is another
sequence of positive real numbers such that $\alpha'_j < 1/2$ for all
$j$.  Thus $\mathcal{F}_j(A')$, $E_j(A')$, and $E(A')$ can be defined
just as before.  For each nonnegative integer $j$, there is a unique
increasing homeomorphism $h_j$ from $[0, 1]$ onto itself with the
property that $h_j$ maps each interval in $\mathcal{F}_j(A)$ onto an
interval in $\mathcal{F}_j(A')$, $h_j$ maps each open subinterval of
$[0, 1]$ which is a component of $[0, 1] \backslash E_j(A)$ onto an
open interval in $[0, 1]$ which is a component of $[0, 1] \backslash
E_j(A')$, and the restriction of $h_j$ to any of these intervals in
the domain is an affine function, i.e., a function of the form $a \, x
+ b$ for some positive real number $a$ and real number $b$.  In other
words, one can think of $\mathcal{F}_j(A)$ and $\mathcal{F}_j(A')$ as
corresponding to partitions of $[0, 1]$ with the same number of
pieces, and $h_j$ is the unique increasing piecewise-linear mapping
from $[0, 1]$ onto itself that maps the partition associated to
$\mathcal{F}_j(A)$ to the partition associated to $\mathcal{F}_j(A')$.

	In fact, if $i$ is a nonnegative integer less than or equal to
$j$, then $h_j$ maps each interval in $\mathcal{F}_i(A)$ onto an
interval in $\mathcal{F}_i(A')$, and $h_j(E_i(A)) = E_i(A')$.
However, the restriction of $h_j$ to an interval in $\mathcal{F}_i(A)$
is not affine in general.  We also have that $h_j$ maps each open
subinterval of $[0, 1]$ which is a component of $[0, 1] \backslash
E_i(A)$ to an open subinterval of $[0, 1]$ which is a component of
$[0, 1] \backslash E_i(A')$, and that the restriction of $h_j$ to such
an open interval is affine.  To be more precise, $h_j = h_i$ on $[0,
1] \backslash E_i(A)$ when $i < j$.  Of course $h_l(0) = 0$ and
$h_l(1) = 1$ for all $l$.  It is not difficult to show that
$\{h_j\}_{j=1}^\infty$ converges uniformly on $[0, 1]$ to an
increasing homeomorphism $h$ from $[0, 1]$ onto itself such that $h =
h_j$ on $[0, 1] \backslash E_j(A)$ and $h(E_j(A)) = E_j(A')$ for all
$j$, and $h(E(A)) = E(A')$.

	Let us return now to the situation of a single sequence $A =
\{\alpha_j\}_{j=1}^\infty$ and the corresponding Cantor set $E(A)$.
There is a natural ``Riemann integral'' for continuous functions on
$E(A)$, which is defined in terms of limits of Riemann sums in the
usual way.  Indeed, let $f(x)$ be a real-valued continuous function on
$E(A)$.  Since $E(A)$ is compact, $f(x)$ is uniformly continuous.  For
each $j \ge 0$ and interval $I \in \mathcal{F}_j(A)$, pick a point
$x(I)$ in $I$ which also lies in $E(A)$, such as one of the endpoints
of $I$.  Define the $j$th Riemann sum of $f(x)$ on $E(A)$ to be
\begin{equation}
	2^{-j} \sum_{I \in \mathcal{F}_j(A)} f(x(I)).
\end{equation}
Using the uniform continuity of $f$ on $E(A)$, it is not too difficult
to show that the sequence of these Riemann sums forms a Cauchy
sequence of real numbers.  Hence it converges, and the limit is our
Riemann integral of $f$ on $E(A)$.

\section{Monotone increasing functions on ${\bf R}$}
\label{section on monotone increasing functions on R}

	Let $f(x)$ be a monotone increasing real-valued function on
${\bf R}$.  Thus
\begin{equation}
\label{f(x) a monotone increasing function on R}
	f(x) \le f(y) \ \hbox{ when } \ x \le y.
\end{equation}
Observe that the sum of two monotone increasing functions is monotone
increasing, that the product of a monotone increasing function with a
nonnegative real number is monotone increasing, and that the product
of two nonnegative monotone increasing functions is monotone
increasing.

	If $x$ is a real number, then the left and right limits
of $f$ at $x$ exist, and are given by
\begin{equation}
	\lim_{y \to x-} f(y) = \sup \{f(y) : y < x\}
\end{equation}
and
\begin{equation}
	\lim_{y \to x+} f(y) = \inf \{f(y) : y > x\}.
\end{equation}
Thus
\begin{equation}
	\lim_{y \to x-} f(y) \le f(x) \le \lim_{y \to x+} f(y),
\end{equation}
and $f$ is continuous at $x$ if and only if these three quantities
are equal.  Notice that if $x_1 < x_2$, then
\begin{equation}
	\lim_{y \to x_1+} f(y) \le \lim_{y \to x_2-} f(y).
\end{equation}

	If $D$ denotes the set of points in ${\bf R}$ at which
$f$ is discontinuous, then it is well-known that $D$ is at most
countable.  Indeed, if $x \in D$, then we can associate to $x$
the open interval
\begin{equation}
	J_x = \bigl(\lim_{y \to x-} f(y), \lim_{y \to x+} f(y)\bigr).
\end{equation}
It is not difficult to see that $J_{x_1} \cap J_{x_2} = \emptyset$
when $x_1, x_2 \in D$.  For each $x \in D$, one can choose a rational
number in $J_x$, and this is a one-to-one mapping from $D$ into the
set of rational numbers.  It follows that $D$ is at most countable,
since the set of rational numbers is a countable set.

	If $a$, $b$ are real numbers such that $a < b$, then
\begin{equation}
   \sum_{x \in D \cap (a, b)} 
	\Bigl(\lim_{y \to x+} f(y) - \lim_{y \to x-} f(y)\Bigr)
		\le f(b) - f(a).
\end{equation}
To see this it suffices to show that the analogous inequality holds
when we sum over $x$ in any finite subset of $D \cap (a, b)$.  This is
not difficult to do, just from the definitions.  Notice also that for
$x \in D \cap (a, b)$, the intervals $J_x$ are contained in $(f(a),
f(b))$, and the inequality states that the sum of the lengths of these
$J_x$'s is less than or equal to the length of the interval $(f(a),
f(b))$.  If $f$ is bounded, with $A \le f(x) \le B$ for some real
numbers $A$, $B$ and all $x \in {\bf R}$, then we get that
\begin{equation}
	\sum_{x \in D} 
		\Bigl(\lim_{y \to x+} f(y) - \lim_{y \to x-} f(y)\Bigr) 
			\le B - A.
\end{equation}

	Let $a$, $b$ be real numbers with $a \le b$ again, and suppose
that $h(x)$ is a real-valued function on $[a, b]$.  A partition of
$[a, b]$ can be given by a finite sequence $\{t_j\}_{j=0}^n$ of real
numbers such that
\begin{equation}
	a = t_0 < \cdots < t_n = b.
\end{equation}
If $P = \{t_j\}_{j=1}^n$ is a partition of $[a, b]$, then the variation
of $h$ associated to this partition is denoted $V_P(h)$ and defined by
\begin{equation}
	V_P(h) = \sum_{j=1}^n |h(t_j) - h(t_{j-1})|.
\end{equation}
Notice that
\begin{equation}
	V_P(h) \ge \biggl| \sum_{j=1}^n (h(t_j) - h(t_{j-1})) \biggr|
		= |h(b) - h(a)|.
\end{equation}
The total variation of $h$ on $[a, b]$ is denoted $V_a^b(h)$ and is
defined to be the supremum of $V_P(h)$ over all partitions $P$ of $[a,
b]$.  This may be $+\infty$, and $h$ is said to be of bounded
variation on $[a, b]$ if $V_a^b(h) < \infty$.  

	For any two real-valued functions $h_1$, $h_2$ on $[a, b]$,
we have that
\begin{equation}
	V_a^b(h_1 + h_2) \le V_a^b(h_1) + V_a^b(h_2).
\end{equation}
This is easy to check just from the definitions.  Also, if $h$ is
a real-valued function on $[a, b]$ and $\alpha$ is a real number, 
then
\begin{equation}
	V_a^b(\alpha h) = |\alpha| \, V_a^b(h),
\end{equation}
where the right side is interpreted as being $0$ is $\alpha = 0$,
even if $V_a^b(h)$ is infinite.  Thus linear combinations of functions
of bounded variation are of bounded variation.

	Here are two basic classes of functions of bounded variation.
First, if $h$ is monotone increasing on $[a, b]$, then $h$ is of
bounded variation, and
\begin{equation}
	V_a^b(h) = h(b) - h(a).
\end{equation}
This is simply because $V_P(h) = h(b) - h(a)$ for all partitions $P$
of $[a, b]$.  Second, if $h$ is a $C$-Lipschitz function on $[a, b]$
for some nonnegative real number $C$, then it is easy to see that
\begin{equation}
	V_a^b(h) \le C \, (b - a).
\end{equation}

	Suppose that $c$ is a real number such that $a < c < b$, and
that $h$ is a real-valued function on $[a, b]$.  Let us check that
\begin{equation}
	V_a^b(h) = V_a^c(h) + V_c^b(h).
\end{equation}
If $P$ is any partition of $[a, b]$, then one can add $c$ to $P$ if
$c$ is not alreay in $P$, and this will only increase the variation of
$h$ on the partition, or keep it the same.  Once one has a partition
that contains $c$, one can split it into two partitions, one of $[a,
c]$ and the other of $[c, b]$.  The variation of $h$ associated to the
whole partition is the equal to the sum of the variations associated
to the two smaller partitions.  This permits one to show that
$V_a^b(h) \le V_a^c(h) + V_c^b(h)$.  Conversely, any pair of
partitions of $[a, c]$ and $[c, b]$ can be combined to give a
partition of $[a, b]$, and the other inequality $V_a^b(h) \ge V_a^c(h)
+ V_c^b(h)$ follows from this.

	If $h$ is a function of bounded variation on $[a, b]$, then
one can write $h(x)$ as
\begin{equation}
	h(x) = V_a^x(h) - (V_a^x(h) - h(x)).
\end{equation}
It is not too difficult to show that $V_a^x(h)$ and $V_a^x(h) - h(x)$
are both monotone increasing functions on $[a, b]$.  Thus a function
$h$ on $[a, b]$ has bounded variation if and only if it can be
expressed as the difference of two monotone increasing functions on
$[a, b]$.

	Now suppose that $E$ is a dense subset of ${\bf R}$ and that
$f_1$, $f_2$ are monotone increasing functions on ${\bf R}$ which are
equal on $E$.  Then $f_1$, $f_2$ are equal at any point in ${\bf R}$
at which they are both continuous.  Using monotonicity one can go a
bit further and show that $f_1$, $f_2$ are continuous at the same
points in ${\bf R}$.  In fact, the set of discontinuities is
determined by the common restriction of $f_1$, $f_2$ to $E$.

	If $E$ is any nonempty subset of ${\bf R}$ and $f(x)$ is a
monotone increasing real-valued function on $E$ which is bounded on
bounded subsets of $E$, then we can extend $f$ to a monotone
increasing function on all of ${\bf R}$, e.g., by setting
\begin{equation}
	F(x) = \sup \{f(y) : y \in E, y \le x \}
\end{equation}
when $E \cap (-\infty, x] = \emptyset$.  If $E \cap (-\infty, x] =
\emptyset$ for some $x$'s, then we can simply set $F(x) = \inf \{f(y)
: y \in E\}$ for these $x$'s.

	Let $\{f_j\}_{j=1}^\infty$ be a sequence of monotone
increasing real-valued functions on ${\bf R}$, and assume that $E$ is
a dense subset of ${\bf R}$ such that $\{f_j(x)\}_{j=1}^\infty$
converges as a sequence of real numbers for each $x$ in $E$.
Let us write $f(x)$ for the limit when $x \in E$.  It is easy to
see that $f(x)$ is monotone increasing on $E$, and $f(x)$ is 
bounded on bounded subsets of $E$ because of monotonicity and the
denseness of $E$.

	Let $F(x)$ be a monotone increasing function on ${\bf R}$
which is equal to $f(x)$ when $x \in E$.  Because $E$ is dense in
${\bf R}$, $F(x)$ is uniquely determined by $f(x)$ except at an at
most countable set of discontinuities, and the set of discontinuities
is determined by $f$ on $E$.  It is not difficult to check that
if $x$ is a point in ${\bf R}$ at which $F$ is continuous, then
$\{f_j\}_{j=1}^\infty$ converges in ${\bf R}$ to $F(x)$.

	Suppose that $\{f_j\}_{j=1}^\infty$ is a sequence of monotone
increasing real-valued functions such that $\{f_j(x)\}_{j=1}^\infty$
is a bounded sequence of real numbers for every $x \in {\bf R}$.
Because of monotonicity, this implies that the $f_j$'s are uniformly
bounded on each bounded subset of ${\bf R}$.  The compactness of
closed and bounded subsets of ${\bf R}$ implies that for each $x \in
{\bf R}$ there is a subsequence of $\{f_j\}_{j=1}^\infty$ which
converges in ${\bf R}$.  By standard arguments about subsequences, if
$E$ is any subset of ${\bf R}$ which is at most countable, then there
is a single subsequence $\{f_{j_k}\}_{k=1}^\infty$ of
$\{f_j\}_{j=1}^\infty$ such that $\{f_{j_k}(x)\}_{k=1}^\infty$ is a
convergent sequence of real numbers for all $x \in E$.  We may as well
choose $E$ to be a dense subset of ${\bf R}$, such as the set of
rational numbers.  By the arguments in the preceding paragraphs, it
follows that $\{f_{j_k}(x)\}_{k=1}^\infty$ is a convergent sequence of
real numbers for all but at most a countable set of $x$'s in ${\bf
R}$.  By passing to a subsequence of $\{f_{j_k}\}_{k=1}^\infty$,
we can get a subsequence of $\{f_j\}_{j=1}^\infty$ which converges
pointwise everywhere on ${\bf R}$.

\section{Some mappings between some $C(M)$'s}
\label{section on some mappings between some C(M)'s}

	Let $E$ be a nonempty closed subset of $[0, 1]$, and suppose
that $f(x)$ is a real-valued continuous function on $E$.  There is a
standard way to define an extension $\widehat{f}(x)$ of $f(x)$ to $[0,
1]$, as follows.  One can express $[0, 1] \backslash E$ as a union of
three types of intervals.  The first type, which is the ``generic''
case, are intervals of the form $(a, b)$, where $a$, $b$ are elements
of $E$.  Thus $f(a)$, $f(b)$ are already defined, and one defines
$\widehat{f}(x)$ on $(a, b)$ by $\widehat{f}(x) = A(x)$, where $A(x) =
\alpha \, x + \beta$ is the affine function such that $A(a) = f(a)$,
$A(b) = f(b)$.  The second type of interval is of the form $[0, c)$,
where $c$ is an element of $E$.  On this interval we define
$\widehat{f}(x)$ to be constant and equal to $f(c)$.  The third type
of interval is of the form $(d, 1]$, where $d$ is an element of $E$.
On this interval we again define $\widehat{f}(x)$ to be constant, and
equal to $f(d)$.  Depending on $E$, any of these three types of
intervals may or may not occur.  If the first type occurs, it may
occur finitely or countably many times.  The second and third types
can occur at most once each.

	It is well-known and not difficult to show that
$\widehat{f}(x)$ as defined in this manner is a continuous function on
$[0, 1]$, using the continuity of $f(x)$ on $E$.  Also, by
construction we have that
\begin{equation}
\label{sup {|widehat{f}(x)| : x in [0, 1]} = sup {|f(x)| : x in E}}
	\sup \{|\widehat{f}(x)| : x \in [0, 1]\}
		= \sup \{|f(x)| : x \in E\}.
\end{equation}
Clearly the mapping $f \mapsto \widehat{f}$ is a linear mapping from
$C(E)$ into $C([0, 1])$, and (\ref{sup {|widehat{f}(x)| : x in [0, 1]}
= sup {|f(x)| : x in E}}) says that this is an isometric embedding
with respect to the standard supremum norms and metrics.

	In the other direction, one has the linear mapping from $C([0,
1])$ into $C(E)$ which takes a continuous function on $[0, 1]$ and
simply restricts it to a continuous function on $E$.  The norm of the
restriction is automatically less than or equal to the norm of the
original function.  Of course if one starts with a continuous function
$f(x)$ on $E$, extends it to the continuous function $\widehat{f}(x)$
on $[0, 1]$ as above, and then restricts that back to $E$, then one
gets $f(x)$ back again.  In other words, the composition of the
mapping $f \mapsto \widehat{f}$ from $C(E)$ into $C([0, 1])$ with the
restriction mapping from $C([0, 1])$ into $C(E)$ is equal to the
identity mapping on $C(E)$.

	Let $n$ be a positive integer with $n \ge 2$, and let us write
$[0, 1]^n$ for the Cartesian product of $n$ copies of $[0, 1]$.  In
other words, this is the unit cube in ${\bf R}^n$, consisting of $x$
in ${\bf R}^n$ such that $0 \le x_j \le 1$ for $j = 1, 2, \ldots, n$.
There is a very simple mapping from $C([0, 1])$ into $C([0, 1]^n)$, in
which a function of one variable is viewed as a function of $n$
variables that only depends on the first variable.  One can also look
at this as an extension mapping from functions on an interval to
functions on a cube containing the interval.  This is a linear mapping
that preserves the supremum norm.  In the other direction, one can
define a restriction mapping from $C([0, 1]^n)$ to $C([0, 1])$ by
taking a function on $[0, 1]^n$ and restricting it to a function on
$[0, 1]$ by setting all but the first coordinate equal to $0$.  This
is a linear mapping for which the supremum norm of the intial function
on $[0, 1]^n$ is automatically greater than or equal to the supremum
norm of the restricted function on $[0, 1]$.  The composition of the
extension mapping from $C([0, 1])$ into $C([0, 1]^n)$ with the
restriction mapping from $C([0, 1]^n)$ into $C([0, 1])$ is equal to
the identity mapping on $C([0, 1])$.

	Here is another kind of mapping from $C([0, 1]^n)$ into $C([0,
1])$.  Let $\phi(t)$ be a continuous function from $[0, 1]$ onto $[0,
1]^n$.  Thus $\phi(t)$ should be a ``space-filling curve''.  Recall
that such a mapping cannot be one-to-one, because then it would be a
homeomorphism between $[0, 1]$ and $[0, 1]^n$, and this is not possible
because $[0, 1]$ with a point removed may not be connected, while
$[0, 1]^n$ with a point removed is always connected when $n \ge 2$.
For the mapping from $C([0, 1]^n)$ into $C([0, 1])$, we simply compose
a continuous real-valued function on $[0, 1]^n$ with $\phi(t)$ to get
a continuous real-valued function on $[0, 1]$.  This defines a linear
mapping from $C([0, 1]^n)$ into $C([0, 1])$ which preserves the
supremum norm, because $\phi(t)$ maps $[0, 1]$ onto $[0, 1]^n$.

\chapter{More on metric spaces}
\label{chapter about more on metric spaces}

\section{Lipschitz mappings}
\label{section on Lipschitz mappings}

	Throughout this section we let $(M, d(x,y))$ and $(N,
\rho(u,v))$ be metric spaces.  A mapping $f : M \to N$ is said to be
\emph{Lipschitz}\index{Lipschitz condition for a mapping between two
metric spaces} if there is a nonnegative real number $C$ such that
\begin{equation}
\label{def of Lipschitz condition}
	\rho(f(x), f(y)) \le C \, d(x,y)
\end{equation}
for all $x, y \in M$.  We also say that $f$ is $C$-Lipschitz in this
case, to mention the constant $C$ explicitly.  Lipschitz mappings
are clearly continuous, and even uniformly continuous.  

	The space of Lipschitz mappings from $M$ to $N$ is denoted
$\Lip (M, N)$, and the space of $C$-Lipschitz mappings from $M$ to $N$
is denoted $\Lip_C (M, N)$.  We write $\Lip (M)$ for the space
of Lipschitz real-valued functions on $M$, and $\Lip_C (M)$ for the
space of $C$-Lipschitz real-valued functions on $M$.

	In the case of real-valued functions on Euclidean spaces which
are continuously differentiable, the Lipschitz condition is equivalent
to the norm of the gradient being bounded, and the supremum of the
norm of the gradient is equal to the smallest choice of the constant
in the Lipschitz condition.  As a special case of this, if $v$ is a
vector in ${\bf R}^n$, then the linear function on ${\bf R}^n$ which
maps a point $x$ in ${\bf R}^n$ to the inner product of $x$ with $v$
is Lipschitz with constant equal to the norm of $v$.  The orthogonal
projection of ${\bf R}^n$ onto any linear subspace is Lipschitz.

\beginremark
\label{Lipschitz functions of order alpha}
{\rm If $\alpha$ is a positive real number, then a mapping 
$f : M \to N$ is said to be Lipschitz of order $\alpha$ if there
is a nonnegative real number $C$ such that 
\begin{equation}
\label{def of Lipschitz of order alpha}
	\rho(f(x), f(y)) \le C \, d(x,y)^\alpha
\end{equation}
for all $x, y \in M$.  Thus $f$ is a Lipschitz mapping in the sense
defined before if and only if it is Lipschitz of order $1$.  When
$0 < \alpha < 1$, a mapping $f : M \to N$ is Lipschitz of order $\alpha$
if and only if it is Lipschitz of order $1$ with respect to the
snowflake metric $d(x,y)^\alpha$ on $M$ and still using the metric
$\rho(u,v)$ on $N$.  When $\alpha > 1$, a mapping $f : M \to N$
is Lipschitz of order $\alpha$ if and only if it is Lipschitz
of order $1$ with respect to $d(x,y)$ on $M$ and the snowflake
metric $\rho(u,v)^{1/\alpha}$ on $N$.  If $M$ is ${\bf R}^n$
with the usual metric and $\alpha > 1$, then the only real-valued 
functions which are Lipschitz of order $\alpha$ are the constant
functions, because such a function would have first derivatives
equal to $0$ everywhere.
}
\end{remark}

	The sum of two real-valued Lipschitz functions on $M$ is a
Lipschitz function, as is the product of a real-valued Lipschitz
function and a constant.  In general, the product of two real-valued
Lipschitz functions is not necessarily Lipschitz, but this is the case
if at least one of the two functions is bounded.  If $f(x)$ is a
real-valued Lipschitz function on $M$ and there is a positive real
number $c$ such that $|f(x)| \ge c$ for all $x$ in $M$, then 
$1/f(x)$ is also a Lipschitz function on $M$.

\beginlemma
\label{d(x,p) and dist(x, A) are 1-Lipschitz}
If $p$ is an element of $M$, then $d(x,p)$ is $1$-Lipschitz as a
real-valued function of $x$ on $M$.  If $A$ is a nonempty subset of
$M$, then $\dist(X, A)$, as defined in (\ref{def of dist(x, A)}) in
Section \ref{section on continuous mappings}, is a $1$-Lipschitz
function on $M$.
\end{lemma}

        This follows from the proofs of Lemmas \ref{d(x,p) is
continuous} and \ref{distance to a subset is continuous} in Section
\ref{section on continuous mappings}.

	Let us say that a mapping $f : M \to N$ is \emph{locally
Lipschitz}\index{locally Lipschitz mappings between metric spaces}
if for every $x \in M$ there is an $r > 0$ so that the restriction
of $f$ to $B(x, r)$ is Lipschitz.  One can check that the sum and
product of real-valued locally Lipschitz functions on $M$ are 
also locally Lipschitz, and that if $f(x)$ is a real-valued locally
Lipschitz function on $M$ such that $f(x) \ne 0$ for all $x$ in $M$,
then $1/f(x)$ is a locally Lipschitz function on $M$.

	A mapping $f : M \to N$ is said to be \emph{countably
Lipschitz}\index{countably Lipschitz mappings between metric spaces}
if $M$ can be expressed as the union of at most countably-many subsets
on which $f$ is Lipschitz.  If $f : M \to N$ is locally Lipschitz and
$M$ is separable, then $f$ is countably Lipschitz.  This follows from
the countable compactness of $M$.

	If $f : M \to N$ is countably Lipschitz, then $M$ can be
expressed as the union of at most countably-many subsets on which $f$
is both Lipschitz and bounded.  This is easy to check, by simply
decomposing further the subsets of $M$ on which $f$ is already
Lipschitz.  Furthermore, the sum and product of real-valued
countably Lipschitz functions is countably Lipschitz, as is the
reciprocal of a nonzero real-valued countably Lipschitz function.

\subsection{Diameters, measure functionals, and dimensions}

	Let us write $\diam_M A$ and $\diam_N B$ for the diameters
of subsets $A$ and $B$ of $M$ and $N$, respectively.

\beginlemma
\label{diameters and Lipschitz mappings}
If $A$ is a subset of $M$ and $f : M \to N$ is a $C$-Lipschitz
mapping, then $\diam_N f(A) \le C \, \diam_M A$.
\end{lemma}

	This is easy to check from the definitions.

\beginproposition
\label{Hausdorff-type measure functionals and Lipschitz mappings}
Let $E$ be a subset of $M$, $\alpha$ a nonegative real number, and $f
: M \to N$ a $C$-Lipschitz mapping.  Under these conditions, we have
that 
\begin{equation}
	H^\alpha_{con}(f(E)) \le C^\alpha \, H^\alpha_{con}(E)
\end{equation}
and
\begin{equation}
	HF^\alpha_{con}(f(E)) \le C^\alpha \, HF^\alpha_{con}(E).
\end{equation}
Similarly, for each $\delta$, $0 < \delta \le \infty$, 
\begin{equation}
	H^\alpha_{C \, \delta}(f(E)) \le C^\alpha \, H^\alpha_\delta(f(E))
\end{equation}
and 
\begin{equation}
	HF^\alpha_{C \, \delta}(f(E)) \le C^\alpha \, HF^\alpha_\delta(f(E)).
\end{equation}
Taking the supremum over $\delta$ it follows that 
\begin{equation}
	H^\alpha(f(E)) \le C^\alpha \, H^\alpha(E)
\end{equation}
and
\begin{equation}
	HF^\alpha(f(E)) \le C^\alpha \, HF^\alpha(E).
\end{equation}
In particular, 
\begin{eqnarray}
	{\dim}_H(f(E)) & \le & {\dim}_H(E), 			\\
	{\dim}_{HF}(f(E)) & \le & {\dim}_{HF}(E),		\\
	{\dim}_{HF^*}(f(E)) & \le & {\dim}_{HF^*}(E).	
\end{eqnarray}
\end{proposition}

	Note that the various measure functionals and dimensions of
$f(E)$ are for subsets of $N$, using the metric $\rho(u,v)$, while
the various measure functionals and dimensions for $E$ are for subsets
of $M$, using the metric $d(x,y)$.  

	The proposition can be shown in a straightforward manner using
Lemma \ref{diameters and Lipschitz mappings} and the various
definitions.  When $\alpha = 0$, one can take $C^\alpha = 1$, and the
conditions on $f$ can be weakened.  For instance, the inequalities
$H^0(f(E)) \le H^0(E)$ and $HF^0(f(E)) \le H^0(E)$ just say that the
number of elements of $f(E)$ is less than or equal to the number of
elements of $E$, and this is true for any mapping $f : M \to N$.

\begincorollary
\label{Hausdorff contents and dimensions and countably-Lipschitz mappings}
Suppose that $f : M \to N$ is countably-Lipschitz.  If $E$ is a subset
of $M$ and $\alpha$ is a nonnegative real number such that
$H^\alpha_{con}(E) = 0$, then $H^\alpha_{con}(f(E)) = 0$.  As a consequence,
$\dim_H(f(E)) \le \dim_H(E)$.
\end{corollary}

\subsection{Some approximations and extensions}

\beginproposition
\label{approx bdd unif continuous fcns by Lipschitz fcns}
Let $f(x)$ be a bounded uniformly continuous real-valued function
on $M$.  There exists a sequence $\{f_j\}_{j=1}^\infty$ of
real-valued Lipschitz functions on $M$ which converges to $f$
uniformly.
\end{proposition}

	For this one can take
\begin{equation}
\label{def of Lip approximations}
	f_j(x) = \inf \{f(y) + j \, d(x,y) : y \in M \}.
\end{equation}
By definition, $f_j(x) \le f(x)$ for all $x$ in $M$ and $j \ge 1$, and
it is not difficult to use the boundedness of $f$ to show that the
infimum is always finite.  More precisely, if $a$ is a real number
such that $f(x) \ge a$ for all $x$ in $M$, then $f_j(x) \ge a$ for all
$x$ in $M$ and $j \ge 1$.

	If $x$, $x'$, and $y$ are arbitrary elements of $M$, then
\begin{equation}
	f_j(x) \le f(y) + j \, d(x,y) \le f(y) + j \, d(x', y) + j \, d(x, x').
\end{equation}
As a result,
\begin{equation}
	f_j(x) \le f_j(x') + j \, d(x, x').
\end{equation}
From this it follows that $f_j$ is $j$-Lipschitz.  Let us also point out
that $f_j(x) = f(x)$ for all $x$ in $M$ when $f$ is $j$-Lipschitz.

	Using the boundedness of $f$, it is not difficult to show that
in fact one only needs to consider $y \in M$ with $d(x, y)$ bounded by
a constant times $1/j$ in the infimum in the definition of $f_j(x)$.
One can then check that $f_j(x)$ converges to $f(x)$ uniformly when
$f$ is uniformly continuous.  We leave the details as an exercise.

\beginproposition
\label{extension of Lipschitz functions}
Let $E$ be a nonempty subset of $M$, let $C$ be a nonnegative real
number, and let $h$ be a real-valued $C$-Lipschitz function on $E$.
There exists a $C$-Lipschitz real-valued function $\widehat{h}$ on $M$
such that $\widehat{h}(x) = h(x)$ when $x \in E$.
\end{proposition}

	In this case we set
\begin{equation}
	\widehat{h}(x) = \inf \{h(y) + C \, d(x,y) : y \in M\}.
\end{equation}
Thus
\begin{equation}
\label{widehat{h}(x) le h(z) + C d(x, z)}
	\widehat{h}(x) \le h(z) + C \, d(x, z)
\end{equation}
for all $x$ in $M$ and $z$ in $E$ by definition.  Next, let us
check that
\begin{equation}
\label{widehat{h}(x) ge h(w) - C d(x, w)}
	\widehat{h}(x) \ge h(w) - C \, d(x, w)
\end{equation}
for all $x \in M$ and $w \in E$, and in particular that the infimum in
the definition of $\widehat{h}(x)$ is finite.  For this it is enough
to show that
\begin{equation}
	h(w) - C \, d(x, w) \le h(y) + C \, d(x,y)
\end{equation}
for all $y$ in $E$, which is the same as
\begin{equation}
	h(w) \le h(y) + C \, d(x, w) + C \, d(x,y).
\end{equation}
This holds because
\begin{equation}
	h(w) \le h(y) + C \, d(w, y) \le h(y) + C \, d(x, w) + C \, d(x, y),
\end{equation}
by the hypothesis that $h$ is $C$-Lipschitz on $E$ and the triangle
inequality.

	It follows in particular from (\ref{widehat{h}(x) le h(z) + C
d(x, z)}) and (\ref{widehat{h}(x) ge h(w) - C d(x, w)}) that
$\widehat{h}(x) = h(x)$ when $x \in E$.  If $x$, $x'$ are arbitrary
elements of $M$ and $y$ is any element of $E$, then
\begin{equation}
	\widehat{h}(x) \le h(y) + C \, d(x, y) 
				\le h(y) + C \, d(x', y) + C \, d(x, x'),
\end{equation}
and hence
\begin{equation}
	\widehat{h}(x) \le \widehat{h}(x') + C \, d(x, x').
\end{equation}
This shows that $\widehat{h}$ is $C$-Lipschitz, and the lemma follows.

\subsection{Spaces of Lipschitz mappings}

	For the rest of the section we assume that $C$ is a fixed
nonnegative real number, and we consider the collection $\Lip_C(M,N)$
of all $C$-Lipschitz mappings from $M$ to $N$.  In particular, if $M$
is bounded, then each Lipschitz mapping from $M$ to $N$ is bounded,
and we can consider $\Lip_C(M,N)$ as a subset of the space $C_b(M,N)$
of bounded continuous mappings from $M$ to $N$.  As in Section
\ref{section on continuous mappings}, we view this space as a metric
space equipped with the metric $\theta(f_1, f_2) = \sup \{\rho(f_1(x),
f_2(x)) : x \in M\}$.  Recall that $C_b(M, N)$ is a complete metric
space with respect to this metric when $N$ is complete, and notice
that $C_b(M, N)$ is bounded if $N$ is bounded.

\beginlemma
\label{Lip_C(M, N) is a closed in various ways}
If $\{f_j\}_{j=1}^\infty$ is a sequence in $\Lip_C(M,N)$ which
converges pointwise to a mapping $f$ from $M$ to $N$, then $f$
is also an element of $\Lip_C(M,N)$.
\end{lemma}

	This is easy to verify from the definitions.

\beginlemma
\label{M, N totally bounded implies that Lip_C(M, N) is too}
If $M$ and $N$ are totally bounded, then $\Lip_C (M, N)$ is
totally bounded as a subset of $C_b(M, N)$.
\end{lemma}

	Exercise.

\beginlemma
\label{pointwise convergence on separable spaces with nice ranges}
Suppose that $M$ is separable and $N$ has the property that closed and
bounded subsets of it are compact.  Let $\{f_j\}_{j=1}^\infty$ be a
sequence of mappings in $\Lip_C(M,N)$ such that
$\{f_j(p)\}_{j=1}^\infty$ is a bounded subset of $N$ for some $p$ in
$M$. Then there is a subsequence $\{f_{j_k}\}_{k=1}^\infty$ of
$\{f_j\}_{j=1}^\infty$ which converges pointwise on $M$, i.e.,
$\{f_{j_k}(x)\}_{k=1}^\infty$ is a convergent sequence in $M$ for
every $x$ in $M$.
\end{lemma}

	Because the $f_j$'s are Lipschitz with bounded constant, the
assumption that $\{f_j(p)\}_{j=1}^\infty$ is a bounded sequence in $N$
for some $p$ in $M$ implies that this holds for all $p$ in $M$.  The
assumption that closed and bounded subsets of $N$ are compact then
implies that for each $x$ in $M$ there is a subsequence of
$\{f_j(x)\}_{j=1}^\infty$ that converges in $M$.  By the separability
of $M$ there is a sequence $\{x_i\}_{i=1}^\infty$ of points in $M$
which is dense in $M$.  Standard diagonalization arguments about
subsequences imply that there is a single subsequence
$\{f_{j_k}\}_{k=1}^\infty$ of $\{f_j\}_{j=1}^\infty$ such that
$\{f_{j_k}(x_i)\}_{k=1}^\infty$ is a convergent sequence in $N$ for
each $i$.  Once one has that $\{f_{j_k}\}_{k=1}^\infty$ converges
pointwise on a dense subset of $M$, it is not difficult to show that
it converges pointwise on all of $M$.

\beginlemma
\label{from pointwise convergence to uniform convergence}
Suppose that $M$ is totally bounded, $N$ is complete, $E$ is a dense
subset of $M$, and $\{f_j\}_{j=1}^\infty$ is a sequence of mappings in
$\Lip_C(M, N)$ which converges pointwise on $E$.  Under these
conditions, $\{f_j\}_{j=1}^\infty$ converges in $C_b(M, N)$ to an
element of $\Lip_C(M, N)$.
\end{lemma}

	This is not difficult to verify.

\beginproposition
\label{compactness of Lip_C(M, N) when M totally bounded and N compact}
Suppose that $M$ is totally bounded and $N$ is compact.  As a subset
of $C_b(M, N)$, $\Lip_C(M, N)$ is then compact.
\end{proposition}

	This is a version of the famous Arzela--Ascoli theorem.  One
way to look at it is to use Lemmas \ref{Lip_C(M, N) is a closed in
various ways} and \ref{M, N totally bounded implies that Lip_C(M, N)
is too} to say that $\Lip_C(M,N)$ is closed and totally bounded as a
subset of $C_b(M,N)$, which is complete as a metric space since $N$
is.  Alternatively, Lemmas \ref{pointwise convergence on separable
spaces with nice ranges} and \ref{from pointwise convergence to
uniform convergence} can be applied to show that every sequence in
$\Lip_C(M,N)$ has a subsequence which converges to an element of
$\Lip_C(M,N)$.

\section{The Hausdorff metric}
\label{section on the Hausdorff metric}

	Let $(M, d(x,y))$ be a metric space.  As in Lemma \ref{lemma
about A_r}, if $A$ is a subset of $M$ and $r$ is a positive real
number, then
\begin{equation}
\label{def of A_r, 2}
	A_r = \bigcup \{B(a, r) : a \in A\}.
\end{equation}
is an open subset of $M$ that contains $A$ and has diameter $\le \diam
A + 2 r$.

	Let $\alpha$ be a nonnegative real number, and let $E$ be a
subset of $M$ such that such that $H^\alpha_{con}(E) < \infty$.  For
every $\epsilon > 0$ there is an $\eta > 0$ so that
\begin{equation}
	HF^\alpha_{con}(E_\eta) < HF^\alpha_{con}(E) + \epsilon.
\end{equation}
More generally, if $0 < \delta \le \infty$, $H^\alpha_\delta(E) < \infty$,
and $\epsilon > 0$, then there is an $\eta > 0$ so that
\begin{equation}
	HF^\alpha_\delta(E_\eta) < HF^\alpha_\delta(E) + \epsilon.
\end{equation}
That is, any covering of $E$ by finitely many subsets of $M$ can be
enlarged slightly to get a covering of $E_\eta$, and one can use this
estimate $HF^\alpha_{con}(E_\eta)$, $HF^\alpha_\delta(E_\eta)$.  This
type of argument does not work in general for $HF^\alpha$ itself, as
one can see by considering the case where $E$ contains a single point.

	If $E_1$, $E_2$ are nonempty subsets of $M$ and $t$ is a
positive real number, then we say that $E_1$ and $E_2$ are
\emph{$t$-close}\index{close subsets of a metric space} if
\begin{equation}
\label{definition of t-closeness, 1}
	\hbox{for every $x \in E_1$ there is a $y \in E_2$ such that }
		d(x,y) < t
\end{equation}
and
\begin{equation}
\label{definition of t-closeness, 2}
	\hbox{for every $y \in E_2$ there is an $x \in E_1$ such that }
		d(x,y) < t.
\end{equation}
This is equivalent to saying that 
\begin{equation}
	E_1 \subseteq (E_2)_t, \ E_2 \subseteq (E_1)_t.
\end{equation}
Notice that if $E_1$, $E_2$ are bounded nonempty subsets of $M$ which
are $t$-close for some $t > 0$, then
\begin{equation}
	|\diam E_1 - \diam E_2| \le 2 t.
\end{equation}

\beginlemma
\label{triangle inequality for closeness}
Suppose that $E_1$, $E_2$, and $E_3$ are nonempty subsets of $M$
and $t$, $u$ are positive real numbers such that $E_1$, $E_2$ are
$t$-close and $E_2$, $E_3$ are $u$-close.  Then $E_1$, $E_3$ are
$(t + u)$-close.
\end{lemma}

	This is easy to verify.

\beginlemma
\label{total boundedness, separability, and t-closeness}
A subset $A$ of $M$ is totally bounded if and only if for
every $t > 0$, $A$ is $t$-close to a finite subset of $M$.
Moreover, $M$ is separable if and only if for every $t > 0$,
$M$ is $t$-close to a subset of $M$ which is at most countable.
\end{lemma}

	This is a straightforward consequence of the definitions.

\beginlemma
\label{t-closeness and epsilon-connectedness}
Suppose that $E$, $F$ are nonempty subsets of $M$ which are $t$-close
for some $t > 0$, and that $E$ is $\epsilon$-connected for some
$\epsilon > 0$.  Then $F$ is $(\epsilon + 2 t)$-connected.
\end{lemma}

	In other words, any two elements of $E$ can be connected
by an $\epsilon$-chain, and because $E$ and $F$ are $t$-close, 
it is easy to see that any two elements of $F$ are connected
by an $(\epsilon + 2t)$-chain.

\subsection{The space of closed and bounded subsets of $M$}

	Let us write $\mathcal{S}(M)$ for the collection of all
nonempty closed and bounded subsets of $M$.  If $E_1$, $E_2$ are two
elements of $\mathcal{S}(M)$, define the \emph{Hausdorff
distance}\index{Hausdorff distance between subsets of a metric space}
$D(E_1, E_2)$ between $E_1$ and $E_2$ by
\begin{equation}
\label{def of Hausdorff distance}
	D(E_1, E_2) = \inf \{t > 0 : E_1, E_2 \hbox{ are $t$-close}\}.
\end{equation}
The boundedness of $E_1$ and $E_2$ ensures that the set of $t > 0$
such that $E_1$, $E_2$ are $t$-close is nonempty.  The restriction to
closed sets implies that $D(E_1, E_2) = 0$ if and only if $E_1 = E_2$.

\beginlemma
\label{Hausdorff distance defines a metric}
The Hausdorff distance $D(E_1, E_2)$ defines a metric on
$\mathcal{S}(M)$.
\end{lemma}

	This is easy to check, using Lemma \ref{triangle inequality
for closeness} to show that the Hausdorff distance satisfies the
triangle inequality.

\beginlemma
\label{totally bounded subsets of M, mathcal{S}(M)}
If $A$ is a totally bounded subset of $M$, then $\mathcal{A} = \{E \in
\mathcal{S}(M) : E \subseteq A\}$ is a totally bounded subset of
$\mathcal{S}(M)$.  Hence if $M$ is totally bounded, then
$\mathcal{S}(M)$ is totally bounded, and separable in particular.
\end{lemma}

	The main point for this lemma is that total boundedness of $A$
implies that $A$ can be approximated by finite sets, and as a result
subsets of $A$ can be approximated by finite subsets of a fixed finite
set for each fixed degree of approximation.

\beginlemma
\label{total boundedness and sequences of subsets of M}
Suppose that $\{E_j\}_{j=1}^\infty$ is a sequence of nonempty closed
and bounded subsets of $M$ which converges to a nonempty closed
and bounded subset $E$ of $M$ in the Hausdorff metric.  If each
$E_j$ it totally bounded, then so is $E$.
\end{lemma}

	This can be derived fairly directly from the definitions.

\beginlemma
\label{a description of the limit of a sequence of c+b subsets}
Suppose that $\{E_j\}_{j=1}^\infty$ is a sequence of nonempty closed
and bounded subsets of $M$ which converges to a nonempty closed and
bounded subset $E$ of $M$.  For each $x$ in $E$, there is a sequence
$\{x_j\}_{j=1}^\infty$ of points in $M$ such that $x_j \in E_j$ for
all $j$ and $\lim_{j \to \infty} x_j = x$.  Conversely, if
$\{j_l\}_{l=1}^\infty$ is a strictly increasing sequence of positive
integers, $x_{j_l} \in E_{j_l}$ for each $l$, and
$\{x_{j_l}\}_{l=1}^\infty$ converges in $M$ to a point $x$, then $x
\in E$.
\end{lemma}

	This characterization of limits of sequences in
$\mathcal{S}(M)$ is a straightforward consequence of the definitions.

\subsection{Subsets, functions, and embeddings}

	There is a natural embedding of $M$ into $\mathcal{S}(M)$, in
which a point $p$ in $M$ is associated to the subset $\{p\}$ of $M$.
This is an isometric embedding, which is to say that
\begin{equation}
	D(\{p\}, \{q\}) = d(p,q)
\end{equation}
for all $p$, $q$ in $M$.

	If $E$ is a nonempty subset of $M$, let us write $f_E(x)$ for
the function $\dist(x, E)$.  Thus $f_E(x)$ is a continuous real-valued
function on $M$ which is $1$-Lipschitz.  We may as well restrict our
attention to nonempty closed subsets $E$ of $M$ here, since $f_E(x) =
f_{\overline{E}}(x)$ for all $x$ in $M$.  

	Assume for the moment that $M$ is bounded.  Then the mapping
$E \mapsto f_E$ defines a mapping from $\mathcal{S}(M)$ into the space
$C_b(M)$ of bounded continuous real-valued functions on $M$.  In fact
this is an isometric embedding of $\mathcal{S}(M)$ into $C_b(M)$, with
respect to the supremum distance on $C_b(M)$.  Explicitly, if $E_1$,
$E_2$ are elements of $\mathcal{S}(M)$, then
\begin{equation}
	D(E_1, E_2) = \sup \{|f_{E_1}(x) - f_{E_2}(x)| : x \in M\}.
\end{equation}
This is not too difficult to verify.

	If $M$ is not bounded, then we can use a basepoint to get an
embedding of $\mathcal{S}(M)$ into $C_b(M)$.  That is, we fix an element
$w$ of $M$, and associate to each nonempty subset $E$ of $M$ the
continuous real-valued function $f_E(x) - f_{\{w\}}(x)$ on $M$.
It is easy to see that this is a bounded function on $M$ when $E$
is a bounded nonempty subset of $M$.  One can also check that this
leads to an isometric embedding of $\mathcal{S}(M)$ into $C_b(M)$,
whether or not $M$ is bounded.

\subsection{Compactness}

\beginlemma
\label{intersection of a decreasing sequence of compact sets}
Suppose that $\{K_j\}_{j=1}^\infty$ is a sequence of nonempty
compact subsets of $M$ such that $K_{j+1} \subseteq K_j$ for
all $j \ge 1$, and put $K = \bigcap_{j=1}^\infty K_j$.  Under
these conditions, $\{K_j\}_{j=1}^\infty$ converges to $K$ with
respect to the Hausdorff metric.
\end{lemma}

	Note that $K$ is a nonempty compact subset of $M$ in this
case, by a well-known result.  To prove the lemma, we use essentially
the same type of argument as for that result.  It suffices to show
that for each $\epsilon > 0$ there is a positive integer $N$ such that
$K_j \subseteq K(\epsilon)$ for all $j \ge N$, where $K(\epsilon)$
denotes the set of $x \in M$ such that $\dist(x, K) < \epsilon$.  This
is the same as $K_\epsilon$ in the sense of (\ref{def of A_r, 2}), but
we write $K(\epsilon)$ now to avoid confusion with $K_j$.  Let
$\epsilon > 0$ be given.  Because $K(\epsilon)$ is an open subset of
$M$, and the union of $K(\epsilon)$ with the complements of the
$K_j$'s is equal to all of $M$, it follows from the compactness of
$K_1$ that a finite collection of these open sets are enough to cover
$K_1$.  The monotonicity property of the $K_j$'s implies more
precisely that $K_1$ is contained in $(M \backslash K_N) \cup
K(\epsilon)$ for some positive integer $N$.  Hence $K_j \subseteq
K(\epsilon)$ when $j \ge N$, as desired.

\beginlemma
\label{completeness of mathcal{S}(M) when M is compact}
If $M$ is compact, then $\mathcal{S}(M)$ is complete as a
metric space with respect to the Hausdorff metric.
\end{lemma}

	In this case $\mathcal{S}(M)$ consists simply of the nonempty
compact subsets of $M$.  Let $\{E_j\}_{j=1}^\infty$ be a Cauchy
sequence in $\mathcal{S}(M)$, which we would like to show converges
in $\mathcal{S}(M)$.  By passing to a subsequence we may assume that
\begin{equation}
\label{D(E_j, E_{j+1}) < 2^{-j}}
	D(E_j, E_{j+1}) < 2^{-j}
\end{equation}
for all $j \ge 1$.  This uses the fact that a Cauchy sequence converges
if any of its subsequences converge.

	For each $j$, set
\begin{equation}
	\widetilde{E}_j = \{x \in M : \dist(x, E_j) \le 2^{-j}\}.
\end{equation}
Each $\widetilde{E}_j$ is a closed subset of $M$, and hence compact.
Also, $D(\widetilde{E}_j, E_j) \le 2^{-j}$, and so if
$\{\widetilde{E}_j\}_{j=1}^\infty$ converges, then so does
$\{E_j\}_{j=1}^\infty$, and they have the same limit.  Thus it
suffices to show that $\{\widetilde{E}_j\}_{j=1}^\infty$ converges.

	A key point now is that
\begin{equation}
	\widetilde{E}_{j+1} \subseteq \widetilde{E}_j
\end{equation}
for all $j$.  This follows from (\ref{D(E_j, E_{j+1}) < 2^{-j}}) and
the definition of $\widetilde{E}_j$.  Therefore
$\{\widetilde{E}_j\}_{j=1}^\infty$ converges in $\mathcal{S}(M)$ to
$\bigcap_{j=1}^\infty \widetilde{E}_j$, as in Lemma \ref{intersection
of a decreasing sequence of compact sets}.

	Note that Lemma \ref{completeness of mathcal{S}(M) when M is
compact} also works when closed and bounded subsets of $M$ are
compact, even if $M$ itself might be unbounded.  This is because any
Cauchy sequence $\{E_j\}_{j=1}^\infty$ in $\mathcal{S}(M)$ has all of
the $E_j$'s contained in a compact subset of $M$.

\beginproposition
\label{M compact implies that mathcal{S}(M) is compact}
If $M$ is compact, then $\mathcal{S}(M)$ is compact, with
respect to the Hausdorff metric.
\end{proposition}

	This follows from Lemmas \ref{totally bounded subsets of M,
mathcal{S}(M)} and \ref{completeness of mathcal{S}(M) when M is
compact}, which imply that $\mathcal{S}(M)$ is totally bounded and
complete when $M$ is compact.

\beginlemma
\label{Hausdorff limits and connectedness}
If $M$ is compact, and $\{A_j\}_{j=1}^\infty$ is a sequence in
$\mathcal{S}(M)$ which converges in the Hausdorff metric to a nonempty
compact subset $A$ of $M$, and each $A_j$ is connected, then $A$ is
connected.
\end{lemma}

	Assume for the sake of a contradiction that $A$ is not
connected.  This implies that $A$ can be expressed as the union of two
nonempty disjoint closed subsets $A^1$, $A^2$.  Because of
compactness, there is an $\epsilon > 0$ so that $d(x, y) \ge \epsilon$
for all $x \in A^1$ and $y \in A^2$.  The convergence of
$\{A_j\}_{j=1}^\infty$ to $A$ in the Hausdorff metric would then imply
that $A_j$ admits a similar decomposition for $j$ large enough, in
contradiction to hypothesis.  This proves the lemma.

\beginlemma
\label{compact and totally disconnected spaces}
Suppose that $M$ is compact and totally disconnected, i.e., that
$M$ does not contain a connected subset with at least two elements.
Then for every $\eta > 0$ there is a $\delta > 0$ so that a
$\delta$-connected subset of $M$ has diameter less than or equal to
$\eta$.
\end{lemma}

	Indeed, suppose for the sake of a contradiction that
there is an $\eta > 0$ so that for each positive integer $n$
there is a $(1/n)$-connected subset $E_n$ of $M$ with diameter
at least $\eta$.  We may as well assume that each $E_n$ is a closed
subset of $M$, since this can be arranged by taking the closure.
Since $M$ is compact, there is a subsequence of the $E_n$'s
which converges in the Hausdorff metric to a nonempty closed
subset $E$.  One can check that $E$ is connected and has diameter
at least $\eta$, contradicting the assumption that $M$ is totally
disconnected.

\section{General forms of integration}
\label{section about general forms of integration}

	Let $(M, d(x,y))$ be a metric space.  We make the standing
assumption in this section that closed and bounded subsets of $M$ are
compact.  If $f$ is a real-valued continuous function on $M$, then we
define the support of $f$ to be the closure of the set of $x \in M$
such that $f(x) \ne 0$, and we denote this by $\supp f$.  Thus the
support of $f$ is a closed subset of $M$ by definition, which is empty
exactly when $f$ is identically equal to $0$.

	The vector space of real-valued continuous functions on $M$
with bounded support is denoted $C_{00}(M)$.  We write $C_0(M)$ for
the vector space of real-valued continuous functions $f$ on $M$ which
``tend to $0$ at infinity'', in the sense that for each $\epsilon > 0$
there is a ball $B$ in $M$ such that $|f(x)| < \epsilon$ for all $x
\in M \backslash B$.  Thus
\begin{equation}
	C_{00}(M) \subseteq C_0(M),
\end{equation}
and of course $C_{00}(M) = C_0(M) = C(M)$, the space of all
real-valued continuous functions on $M$, when $M$ is compact.  Because
of our standing assumptions on $M$, functions in $C_{00}(M)$ are
bounded and uniformly continuous, and indeed the same is true for
functions in $C_0(M)$.

	Recall that $C_b(M)$ denotes the vector space of bounded
real-valued continuous functions on $M$.  If $f$ is an element
of $C_b(M)$, then the supremum norm of $f$ is defined by
\begin{equation}
\label{recalling the supremum norm}
	\|f\|_{sup} = \sup_{x \in M} |f(x)|.
\end{equation}
The supremum metric on $C_b(M)$ is defined by $\|f_1 - f_2\|_{sup}$,
and $C_b(M)$ is complete as a metric space when equipped with this
metric.

	It is not difficult to show that $C_0(M)$ is a closed subset
of $C_b(M)$ with respect to the supremum metric, so that $C_0(M)$ is
also a complete metric space when equipped with this metric.  In fact,
$C_0(M)$ is equal to the closure of $C_{00}(M)$ in $C_b(M)$.
We leave it as an exercise to the reader to check that $C_0(M)$
is separable as a metric space.  Of course this is very similar to
the separability of $C(N)$ when $N$ is a compact metric space.

	If $\{\phi_j\}_{j=1}^\infty$ is a sequence of functions in
$C_{00}(M)$ and $\phi$ is another function in $C_{00}(M)$, then we
define \emph{restricted convergence of $\{\phi_j\}_{j=1}^\infty$ to
$\phi$}\index{restricted convergence of a sequence in $C_{00}(M)$} to
mean that there is a compact subset $K$ of $M$ such that the support
of $\phi_j$ is contained in $K$ for all $j$ and
$\{\phi_j\}_{j=1}^\infty$ converges to $\phi$ uniformly.  In this case
the support of $\phi$ is also contained in $K$.  A subset
$\mathcal{E}$ is $C_{00}(M)$ is said to be \emph{dense in $C_{00}(M)$
with respect to restricted convergence}\index{dense subsets of
$C_{00}(M)$ with respect to restricted convergence} if for every
function $\phi$ in $C_{00}(M)$ there is a sequence of functions
$\{\phi_j\}_{j=1}^\infty$ in $\mathcal{E}$ with restricted convergence
to $\phi$.  One can check that $C_{00}(M)$ has a subset which is at
most countable and dense with respect to restricted convergence.

\subsection{Nonnegative linear functionals on $C_{00}(M)$}
\index{nonnegative linear functionals on $C_{00}(M)$}

	Let $\lambda$ be a nonnegative linear functional on
$C_{00}(M)$.  In other words, $\lambda$ is a linear mapping from the
vector space $C_{00}(M)$ to the real numbers such that $\lambda(f) \ge
0$ whenever $f$ is a function on $M$ in $C_{00}(M)$ which is
nonnegative, which is to say that $f(x) \ge 0$ for all $x \in M$.

\beginexamples
\label{examples of nonnegative linear functionals on C_{00}(M)}
{\rm (i)  If $M = {\bf R}^n$, then $\lambda(f)$ equal to the ordinary
Riemann integral of $f$ defines a nonnegative linear functional
on $C_{00}({\bf R}^n)$.

(ii) The ``Riemann integrals'' associated to Cantor sets $E(A)$
as in Section \ref{section on cantor sets} define nonnegative linear
functionals on $C_{00}({\bf R})$. 

(iii) In ${\bf R}^n$ when $n \ge 2$ one can take the integral of a
function over a compact submanifold of dimension $k < n$, with respect
to the element of $k$-dimensional volume.

(iv) For any $M$ one has the nonnegative linear functionals
corresponding to ``Dirac masses'',\index{Dirac masses} in which one
chooses a point $p \in M$ and defines $\lambda(f)$ to simply be
$f(p)$.  }
\end{examples}

	Notice that if $\lambda_1$, $\lambda_2$ are nonnegative linear
functionals on $C_{00}(M)$ and $a_1$, $a_2$ are nonnegative real
numbers, then $a_1 \, \lambda_1 + a_2 \, \lambda_2$ is a nonnegative
linear functional on $C_{00}(M)$.

\beginremark
\label{remark about f --> f_+, f_-}
{\rm If $f(x)$ is a real-valued function on $M$, define $f_+(x)$
and $f_-(x)$ by
\begin{equation}
\label{def of f_+(x), f_-(x)}
	f_+(x) = \max(f(x), 0), \quad f_-(x) = \max(-f(x), 0).
\end{equation}
Thus $f_+(x)$, $f_-(x)$ are nonnegative functions on $M$ and $f(x) =
f_+(x) - f_-(x)$.  If $f(x)$ is continuous, then so are $f_+(x)$,
$f_-(x)$.  Of course the supports of $f_+(x)$, $f_-(x)$ are contained
in the support of $f(x)$, and $f_+(x)$, $f_-(x)$ are bounded functions
on $M$ if $f(x)$ is.  }
\end{remark}

	Suppose that $f$, $h$ are functions in $C_{00}(M)$ such that
$h$ is nonnegative and $|f(x)| \le h(x)$ for all $x$ in $M$.  Then
$h + f$ and $h - f$ are nonnegative functions, so that $\lambda(h + f)$
and $\lambda(h - f)$ are nonnegative real numbers.  This implies that
\begin{equation}
\label{|lambda(f)| le lambda(h)}
	|\lambda(f)| \le \lambda(h).
\end{equation}

	As an application of this, suppose that $B$ is a nonempty
closed and bounded subset of $M$, and let $\phi_B(x)$ be a function in
$C_{00}(M)$ such that $0 \le \phi_B(x) \le 1$ for all $x \in M$ and
$\phi_B(x) = 1$ when $x \in B$.  For instance, one can take
$\phi_B(x)$ to be a simple piecewise-linear function of $\dist(x, B)$.
If $f$ is a function in $C_{00}(M)$ with support contained in $B$,
then
\begin{equation}
\label{|lambda(f)| le lambda(phi_B) ||f||_{sup}}
	|\lambda(f)| \le \lambda(\phi_B) \, \|f\|_{sup}.
\end{equation}
Indeed, it is enough to check this when $|f(x)| \le 1$ for all $x \in
B$, since one can reduce to this case by multiplying $f$ by a positive
real number.  In this case, $|f(x)| \le \phi_B(x)$ for all $x \in M$,
and (\ref{|lambda(f)| le lambda(phi_B) ||f||_{sup}}) follows from
(\ref{|lambda(f)| le lambda(h)}).  

	The inequality (\ref{|lambda(f)| le lambda(phi_B)
||f||_{sup}}) implies that $\lambda$ is \emph{continuous with respect
to restricted convergence},\index{continuity of a nonnegative linear
functional on $C_{00}(M)$ with respect to restrcied convergence} in
the sense that if $\{f_j\}_{j=1}^\infty$ is a sequence in $C_{00}(M)$
and $f$ is a function in $C_{00}(M)$ with restricted convergence of
$\{f_j\}_{j=1}^\infty$ to $f$, then
\begin{equation}
\label{continuity of lambda with respect to restricted convergence}
	\lim_{j \to \infty} \lambda(f_j) = f.
\end{equation}
As a consequence of this, if $\lambda_1$, $\lambda_2$ are two
nonnegative linear functionals on $C_{00}(M)$ which are equal on a
dense subset $\mathcal{E}$ of $C_{00}(M)$ with respect to restricted
convergence, then $\lambda_1$, $\lambda_2$ are equal on all of
$C_{00}(M)$.

	We say that $\lambda$ is a \emph{bounded nonnegative linear
functional on $C_{00}(M)$}\index{bounded nonnegative linear functional
on $C_{00}(M)$} if there is a nonnegative real number $C$ such that
\begin{equation}
\label{def of lambda being a bounded linear functional}
	|\lambda(f)| \le C \, \|f\|_{sup}
\end{equation}
for all $f$ in $C_{00}(M)$.  If $M$ is compact, then this holds
automatically with $C = \lambda({\bf 1})$, where ${\bf 1}$ denotes the
constant function on $M$ equal to $1$.

	If $\lambda$ is bounded, then $\lambda$ defines a uniformly
continuous and in fact Lipschitz mapping from $C_{00}(M)$ to the real
numbers with respect to the supremum metric on $C_{00}(M)$, and hence
extends to a continuous mapping from $C_0(M)$ to the real numbers
which is also linear and satisfies (\ref{def of lambda being a bounded
linear functional}).  Using linearity and nonnegativity, we can extend
$\lambda$ in a natural manner to all of $C_b(M)$, as follows.

	Fix a basepoint $p$ in $M$.  For each positive integer $j$,
define $\phi_j \in C_{00}(M)$ by $\phi_j(x) = 1$ when $d(x, p) \le j$,
$\phi_j(x) = d(x,p) - j$ when $j \le d(x,p) \le j+1$, and $\phi_j(x) =
0$ when $d(x,p) \ge j+1$.  Thus $\phi_j(x)$ is a continuous
piecewise-linear function of $d(x,p)$, $0 \le \phi_j(x) \le 1$ for
all $j \ge 1$ and $x \in M$, and $\phi_j(x) \le \phi_{j+1}(x)$
for all $j$ and $x$.  As a result, 
\begin{equation}
	0 \le \lambda(\phi_j) \le \lambda(\phi_{j+1}) \le C
\end{equation}
for all $j$, and therefore the sequence
$\{\lambda(\phi_j)\}_{j=1}^\infty$ converges.  In particular, for
each $\epsilon > 0$ there is a positive integer $L$ such that
\begin{equation}
	0 \le \lambda(\phi_k - \phi_j) < \epsilon
				\quad\hbox{when } k \ge j \ge L.
\end{equation}
Now let $f$ be any bounded real-valued continuous function on $M$.
Clearly
\begin{equation}
	|f(x) (\phi_k(x) - \phi_j(x))| 
		\le \|f\|_{sup} \, (\phi_k(x) - \phi_j(x))
\end{equation}
for all $x \in M$ when $k \ge j$, and hence
\begin{equation}
	|\lambda(f (\phi_k - \phi_j))| 
		\le \|f\|_{sup} \, \lambda(\phi_k - \phi_j).
\end{equation}
Thus $\{\lambda(\phi_j \, f)\}_{j=1}^\infty$ is a Cauchy sequence of
real numbers, and therefore a convergent sequence.  One can define
$\lambda(f)$ to be the limit of this sequence, and this extends
$\lambda$ to a nonnegative linear functional on $C_b(M)$.

	It is easy to see that this extension of $\lambda$ to $C_b(M)$
is also bounded in the sense that (\ref{def of lambda being a bounded
linear functional}) for all $f$ in $C_b(M)$, and with the same
constant $C$ as on $C_{00}(M)$.  This implies that the extension of
$\lambda$ to $C_b(M)$ is continuous as a mapping from $C_b(M)$ to
${\bf R}$, with respect to the supremum metric.  There is a stronger
continuity condition that this extension of $\lambda$ to $C_b(M)$
satisfies, which is that if $\{f_j\}_{j=1}^\infty$ is a sequence of
functions in $C_b(M)$ and $f$ is another function in $C_b(M)$, if the
supremum norms $\|f_j\|_{sup}$, $j \in {\bf Z}_+$, are bounded, and if
for each bounded subset $B$ of $M$ the restriction of the $f_j$'s to
$B$ converges uniformly to the restriction of $f$ to $B$, then
\begin{equation}
\label{lim_{j to infty} lambda(f_j) = lambda(f)}
	\lim_{j \to \infty} \lambda(f_j) = \lambda(f).
\end{equation}
To put it a bit differently, the boundedness of $\lambda$ implies that
(\ref{lim_{j to infty} lambda(f_j) = lambda(f)}) holds when
$\{f_j\}_{j=1}^\infty$ converges to $f$ uniformly, and in fact we have
that (\ref{lim_{j to infty} lambda(f_j) = lambda(f)}) also holds under
the more general convergence conditions described above.  This is not
too hard to check, using the same kind of considerations as before.

	Another nice feature of this extension of $\lambda$ to
$C_b(M)$ is that if $f(x)$ is a function in $C_b(M)$ which is
nonnegative, then
\begin{equation}
	\lambda(f) = \sup \{\lambda(\phi) : \phi \in C_{00}(M), 
					0 \le \phi(x) \le f(x)
					\hbox{ for all } x \in M\}.
\end{equation}
Again, this is not too difficult to check, using the same kind of
considerations as in the definition of the extension of $\lambda$ to
$C_b(M)$.

\beginremark
\label{pushing forward linear functionals on C_{00}(M)}
{\rm Let $M_1$, $M_2$ be metric spaces as in this section, with
closed and bounded subsets being compact.  Suppose that
$\rho$ is a continuous mapping from $M_1$ to $M_2$ which
is \emph{proper}, which means that the inverse image of a compact
subset of $M_2$ under $\rho$ is a compact subset of $M_1$.
If $\lambda_1$ is a nonnegative linear functional on $C_{00}(M_1)$,
then we can get a nonnegative linear functional $\lambda_2$ on
$C_{00}(M_2)$ using $\rho$ by the formula
\begin{equation}
	\lambda_2(f) = \lambda_1(f \circ \rho)
\end{equation}
for $f$ in $C_{00}(M_2)$.  That is, $f \circ \rho$ is an element
of $C_{00}(M_1)$ when $f$ is an element of $C_{00}(M_2)$, under the
assumption that $\rho : M_1 \to M_2$ is continuous and proper.
If $\lambda_1$ is also bounded, then the assumption that $\rho$ is
proper is not really needed, because $f \circ \rho$ is a bounded
continuous function on $M_1$ whenever $f$ is a bounded continuous
function on $M_2$ and $\rho : M_1 \to M_2$ is continuous.  Of
course $\lambda_2$ is then bounded as well.
}
\end{remark}

\subsection{Supports of nonnegative linear functionals}

	Let us return to the situation where $\lambda$ is a
nonnegative linear functional on $C_{00}(M)$, which may or may not be
bounded.  If $U$ is an open subset of $M$, then we say that $\lambda$
vanishes on $U$ if $\lambda(f) = 0$ whenever $f \in C_{00}(M)$ has
support contained in $U$.  For example, if $\phi$ is a nonnegative
function in $C_{00}(M)$ and $\lambda(\phi) = 0$, then $\lambda$
vanishes on the open set $\{x \in M : \phi(x) > 0\}$.

	Suppose that $U_1$, $U_2$ are open subsets of $M$ such that
$\lambda$ vanishes on $U_1$, $U_2$, and let us show that $\lambda$
vanishes on the union $U_1 \cup U_2$ as well.  Let $f$ be any function
in $C_{00}(M)$ such that the support of $f$ is contained in $U_1 \cup
U_2$.  Define functions $f_1$, $f_2$ by
\begin{equation}
	f_i(x) = f(x) \, 
   \frac{\dist(x, M \backslash U_{3-i})}{\dist(x, M \backslash U_1) 
					+ \dist(x, M \backslash U_2)},
						\qquad i = 1, 2,
\end{equation}
where this is interpreted as being $0$ when $x \in M \backslash (U_1
\cup U_2)$.  It is not difficult to see that $f_1$, $f_2$ are elements
of $C_{00}(M)$ with support contained in $U_1$, $U_2$, respectively,
and which satisfy $f_1 + f_2 = f$.  Thus
\begin{equation}
	\lambda(f) = \lambda(f_1) + \lambda(f_2) = 0,
\end{equation}
because $\lambda$ vanishes on $U_1$, $U_2$.  This shows that $\lambda$
vanishes on $U_1 \cup U_2$.

	More generally, suppose that $\{U_i\}_{i \in I}$ is an
arbitrary family of open subsets of $M$ on which $\lambda$ vanishes.
Then $\lambda$ vanishes on the union of all the $U_i$'s.  Indeed, let
$f$ be a function in $C_{00}(M)$ with support contained in the union
of the $U_i$'s.  Since the support of $f$ is compact, it is contained
in the union of finitely many $U_i$'s.  This permits one to apply the
result of the preceding paragraph.

	As a consequence, there is an open subset $U$ of $M$ on which
$\lambda$ vanishes, and which contains any open subset of $M$ on
which $\lambda$ vanishes.  Put $E = M \backslash U$.  This is a
closed subset of $M$, and it is called the support of $\lambda$.

	Suppose that $f$ is a function in $C_{00}(M)$ such that $f(x)
= 0$ for all $x \in E$.  For each positive integer $j$, let $\psi_j$
be the real-valued continuous function on $M$ such that $\psi_j(x) =
0$ when $\dist(x, F) \le 1/j$, $\psi_j(x) = j \dist(x, F) - 1$ when
$1/j \le \dist(x, F) \le 2/j$, and $\psi_j(x) = 1$ when $\dist(x, F)
\ge 2/j$.  Thus $\psi_j$ is a continuous piecewise-linear function of
$\dist(x, F)$.  It is not difficult to see that the functions
$\psi_j(x) \, f(x)$ converge uniformly to $f(x)$ as $j \to \infty$,
and that the support of $\psi_j \, f$ is contained in the support of
$f$ and in $M \backslash U$ for each $j$.  Hence $\lambda(\psi_j \, f)
= 0$ for all $j$, and $\lambda(\psi_j \, f)$ converges to $\lambda(f)$
as $j \to \infty$.  Thus $\lambda(f) = 0$.

	If the support of $\lambda$ is a compact subset of $M$, and if
$\phi$ is a nonnegative function in $C_{00}(M)$ such that $\phi(x) =
1$ for all $x$ in the support of $\lambda$, then we have that
\begin{equation}
	\lambda(f) = \lambda(\phi \, f)
\end{equation}
for all $f$ in $C_{00}(M)$.  Because $|\phi(x) \, f(x)| \le
\|f\|_{sup} \, \phi(x)$ for all $x \in M$, we get that
\begin{equation}
	|\lambda(f)| \le \lambda(\phi) \, \|f\|_{sup}
\end{equation}
for all $f$ in $C_{00}(M)$.  Thus $\lambda$ is bounded when it has
compact support.

\subsection{Weak convergence}

	Let $\{\lambda_j\}_{j=1}^\infty$ be a sequence of nonnegative
linear functionals on $C_{00}(M)$.  We say that $\{\lambda_j\}_{j=1}^\infty$
is \emph{uniformly bounded on bounded subsets of $M$}\index{uniformly
bounded on bounded subsets on a metric space $M$, for a sequence of
nonnegative linear functionals on $C_{00}(M)$} if for each bounded
subset $B$ of $M$ there is a nonnegative real number $C_B$ such that
\begin{equation}
\label{uniformly bounded on bounded subsets of M condition}
	|\lambda_j(f)| \le C_B \, \|f\|_{sup}
\end{equation}
for all $f$ in $C_{00}(M)$ with support contained in $B$ and for all
$j$.  This condition clearly implies that
$\{\lambda_j(f)\}_{j=1}^\infty$ is a bounded sequence of real numbers
for each $f$ in $C_{00}(M)$.  This converse also works, as one can
show using (\ref{|lambda(f)| le lambda(phi_B) ||f||_{sup}}).

	We say that $\{\lambda_j\}_{j=1}^\infty$ is \emph{uniformly
bounded}\index{uniformly bounded on a metric space $M$, for a sequence
of nonnegative linear functionals on $C_{00}(M)$} if there is a
nonnegative real number $C$ such that
\begin{equation}
\label{uniformly bounded condition}
	|\lambda_j(f)| \le C \, \|f\|_{sup}
\end{equation}
for all $f$ in $C_{00}(M)$ and all $j$.  In particular, each
$\lambda_j$ should be a bounded nonnegative linear functional on
$C_{00}(M)$ for this condition to hold.

	Now suppose that $\{\lambda_j\}_{j=1}^\infty$ is a sequence of
nonnegative linear functionals on $C_{00}(M)$, and that $\lambda$ is
another nonnegative linear functional on $C_{00}(M)$.  We say that
$\{\lambda_j\}_{j=1}^\infty$ \emph{converges weakly to
$\lambda$}\index{weak convergence of a sequence of nonnegative linear
functionals on $C_{00}(M)$ to another nonnegative linear functional on
$C_{00}(M)$} if
\begin{equation}
\label{def of weak convergence for a seq. of nonnegative linear functionals}
	\lim_{j \to \infty} \lambda_j(f) = \lambda(f)
\end{equation}
for all $f$ in $C_{00}(M)$.  Sometimes this is referred to as
\emph{$weak^*$ convergence}, but we shall use ``weak convergence'' for
simplicity.  If $\{\lambda_j\}_{j=1}^\infty$ converges weakly to
$\lambda$, then $\{\lambda_j(f)\}_{j=1}^\infty$ is a bounded sequence
of real numbers for every $f$ in $C_{00}(M)$, and hence
$\{\lambda_j\}_{j=1}^\infty$ is uniformly bounded on bounded subsets
of $M$.  

	Notice that if $\{\lambda_j\}_{j=1}^\infty$ is a uniformly
bounded sequence of nonnegative linear functionals on $C_{00}(M)$
which converges weakly to $\lambda$, then $\lambda$ is a bounded
nonnegative linear functional.

	As a basic example, suppose that $\{p_j\}_{j=1}^\infty$ is a
sequence of points in $M$ and that $p$ is another point in $M$, and
define linear functionals $\lambda_j$, $\lambda$ on $C_{00}(M)$ by
$\lambda_j(f) = f(p_j)$, $\lambda(f) = f(p)$.  Then
$\{\lambda_j\}_{j=1}^\infty$ converges weakly to $\lambda$ if and
only if $\{p_j\}_{j=1}^\infty$ converges to $p$ in $M$.

	Let $\{\lambda_j\}_{j=1}^\infty$ be a sequence of nonnegative
linear functionals on $C_{00}(M)$ which is uniformly bounded on
bounded subsets of $M$, and let $\mathcal{E}$ be a dense subset of
$C_{00}(M)$ with respect to restricted convergence.  Assume that
for every $f$ in $\mathcal{E}$, $\{\lambda_j(f)\}_{j=1}^\infty$ is
a convergent sequence of real numbers.  Using the uniform boundedness
of $\{\lambda_j\}_{j=1}^\infty$ on bounded subsets of $M$, one can
check that $\{\lambda_j(f)\}_{j=1}^\infty$ is a convergent sequence
of real numbers for every $f$ in $\mathcal{E}$.  Denoting the limit
as $\lambda(f)$, one can also check that $\lambda(f)$ is a nonnegative
linear functional on $C_{00}(M)$.  Thus $\{\lambda_j\}_{j=1}^\infty$
converges weakly to $\lambda$ in the sense described above.

	Suppose that $\{\lambda_j\}_{j=1}^\infty$ is a sequence of
nonnegative linear functionals on $C_{00}(M)$ which is uniformly
bounded on bounded subsets of $M$, and that $\mathcal{E}$ is a dense
subset of $C_{00}(M)$ with respect to restricted convergence which is
at most countable as well.  For any function $f$ in $C_{00}(M)$, the
sequence $\{\lambda_j(f)\}_{j=1}^\infty$ is a bounded sequence of real
numbers, and hence has a subsequence that converges.  By standard
arguments, there is a subsequence $\{\lambda_{j_l}\}_{l=1}^\infty$ of
$\{\lambda_j\}_{j=1}^\infty$ such that
$\{\lambda_{j_l}(f)\}_{l=1}^\infty$ is a convergent sequence of real
numbers for every $f$ in $\mathcal{E}$.  As in the remarks in the
preceding paragraph, $\{\lambda_{j_l}\}_{l=1}^\infty$ converges
weakly to a nonnegative linear functional on $C_{00}(M)$.

\section{A few basic inequalities}
\label{section on a few basic inequalities}

	Let $J_1, \ldots, J_n$ be a collection of bounded intervals
in the real line, and let $a_1, \ldots, a_m$ be positive real numbers.
Consider the function 
\begin{equation}
	\phi(x) = \sum_{p=1}^m a_p \, {\bf 1}_{J_p}(x)
\end{equation}
on ${\bf R}$, and the associated quantity
\begin{equation}
\label{sum_{p=1}^m a_p |J_p|}
	\sum_{p=1}^m a_p \, |J_p|,
\end{equation}
which is the total integral of $\phi$.

	If $L_1, \ldots, L_n$ are disjoint intervals in ${\bf R}$,
then
\begin{equation}
	\sum_{q=1}^n |J \cap L_q| \le |J|
\end{equation}
for every interval $J$ in ${\bf R}$.  As a result,
\begin{equation}
	\sum_{q=1}^n \sum_{p=1}^m a_p \, |J_p \cap L_q| 
			\le \sum_{p=1}^m a_p \, |J_p|.
\end{equation}

	If $t$ is a positive real number and $L$ is an
interval in ${\bf R}$ such that $\phi(x) > t$ for all $x \in L$,
then
\begin{equation}
\label{|L| le t^{-1} sum_{p=1}^m a_p |J_p cap L|, 1}
	|L| \le t^{-1} \sum_{p=1}^m a_p \, |J_p \cap L|.
\end{equation}
Using this one can check that
\begin{equation}
\label{chebychev's inequality, 1}
	HF^1_{con}(\{x \in {\bf R} : \phi(x) > t\}) 
			\le t^{-1} \sum_{p=1}^m a_p \, |J_p|.
\end{equation}
It may be that $\phi(x) \le t$ for all $x \in {\bf R}$, so that the
left side of (\ref{chebychev's inequality, 1}) is equal to $0$, and
the inequality holds automatically.  If $\{x \in {\bf R} : \phi(x) >
t\} \ne \emptyset$, then $\{x \in {\bf R}: \phi(x) > t \}$ can be
expressed as a finite union of bounded disjoint intervals, and the
$1$-dimensional $HF$-content of this set is less than or equal to the
sum of the lengths of these intervals.  This permits one to derive
(\ref{chebychev's inequality, 1}) from (\ref{|L| le t^{-1} sum_{p=1}^m
a_p |J_p cap L|, 1}).

	Now suppose that $A$ is a nonempty set which is at most
countable, and that $\{J_p\}_{p \in A}$ is a family of bounded open
intervals in ${\bf R}$ indexed by $A$.  Also let $\{a_p\}_{p \in A}$
be a family of positive real numbers indexed by $p$.  As before, set
\begin{equation}
	\phi(x) = \sum_{p \in A} a_p \, {\bf 1}_{J_p}(x)
\end{equation}
for $x$ in ${\bf R}$, and consider the associated quantity
\begin{equation}
	\sum_{p \in A} a_p \, |J_p|.
\end{equation}
The restriction to open intervals here is not too serious, and one
can always enlarge intervals slightly to get open intervals, without
increasing the sum above too much.

	It may be that $\phi(x) = +\infty$ for some $x$'s, so that
$\phi(x)$ may be an extended real-valued function in this case.  We
shall be concerned with the situation where $\sum_{i \in I} a_i \,
|J_i|$ is finite, and we shall see that the set where $\phi(x)$ is
infinite is small in a reasonable sense.  

	Observe that for each positive real number $t$ and each
$x \in {\bf R}$ such that $\phi(x) > t$, there is a finite
subset $A_1$ of $A$ such that 
\begin{equation}
	\sum_{p \in A_1} a_p \, {\bf J}_p(x) > t.
\end{equation}
Because the $J_p$'s are assumed to be open intervals, it is
not difficult to check that
\begin{equation}
	\{x \in {\bf R} : \phi(x) > t \}
\end{equation}
is an open subset of ${\bf R}$.  In other words, $\phi(x)$ is lower
semicontinuous.

	As before, notice that if $L_1, \ldots, L_n$ are disjoint
intervals in ${\bf R}$, then
\begin{equation}
	\sum_{q=1}^n \sum_{p \in A} a_p \, |J_p \cap L_q|
		\le \sum_{p \in A} a_p \, |J_p|.
\end{equation}
This is again basically just a matter of rearranging sums.

	If $t$ is a positive real number and $L$ is a closed and
bounded interval in ${\bf R}$ such that $\phi(x) > t$ for all $x$ in
$L$, then
\begin{equation}
\label{|L| le t^{-1} sum_{p in A} a_p |J_p cap L|, 2}
	|L| \le t^{-1} \sum_{p \in A} a_p \, |J_p \cap L|.
\end{equation}
Indeed, because $L$ is compact, one can verify that there is a finite
subset $A_1$ of $A$ such that 
\begin{equation}
	\sum_{p \in A_1} a_p \, {\bf 1}_{J_p}(x) > t 
			\quad\hbox{for all } x \in L,
\end{equation}
so that (\ref{|L| le t^{-1} sum_{p=1}^m a_p |J_p cap L|, 1}) may be
applied.  Once one has (\ref{|L| le t^{-1} sum_{p in A} a_p |J_p cap
L|, 2}) for closed and bounded intervals, the same inequality follows
for arbitrary intervals, by approximating an arbitrary interval by
closed and bounded subintervals.

	If $t$ is a positive real number, then let us check that
\begin{equation}
\label{chebychev's inequality, 2}
	H^1_{con}(\{x \in {\bf R} : \phi(x) > t \}) 
			\le t^{-1} \sum_{p \in A} a_p \, |J_p|.
\end{equation}
We may as well assume that there is at least one $x \in {\bf R}$
such that $\phi(x) > t$, since otherwise the left side of
the inequality is equal to $0$, and the inequality holds automatically.
Because $\{x \in {\bf R} : \phi(x) > t \}$ is an open subset
of ${\bf R}$, it can be expressed as the union of a family
$\{L_q\}_{q \in B}$ of disjoint open intervals in ${\bf R}$,
where $B$ is a nonempty set which is at most countable.  Thus
it suffices to show that
\begin{equation}
	\sum_{q \in B} |L_q| \le t^{-1} \sum_{p \in A} a_p \, |J_p|.
\end{equation}
To get this it is enough to show that for each finite subset $B'$ of
$B$ and closed and bounded intervals $L'_q \subseteq L_q$, $q \in B'$,
we have that
\begin{equation}
	\sum_{q \in B'} |L'_q| \le t^{-1} \sum_{p \in A} a_p \, |J_p|.
\end{equation}
By compactness, there is a finite subset $A_1$ of $A$ such that
\begin{equation}
	\sum_{p \in A_1} a_p \, {\bf 1}_{J_p}(x) > t
			\quad\hbox{for all } x \in \bigcup_{q \in B'} L'_q,
\end{equation}
and it is enough to show that
\begin{equation}
	\sum_{q \in B'} |L'_q| \le t^{-1} \sum_{p \in A_1} a_p \, |J_p|.
\end{equation}
This inequality can be derived from (\ref{|L| le t^{-1} sum_{p=1}^m
a_p |J_p cap L|, 1}).

	As a consequence of (\ref{chebychev's inequality, 2}), notice
that 
\begin{equation}
	H^1_{con}(\{x \in {\bf R} : \phi(x) = + \infty\}) = 0
\end{equation}
when $\sum_{p \in A} a_p \, |J_p| < \infty$.

	Now suppose that $(M, d(x,y))$ is a metric space, $\alpha$ is
a real number such that $\alpha \ge 1$, and $f$ is a $C$-Lipschitz
real-valued function on $M$ for some $C > 0$.  Also let $E$ be a
nonempty bounded subset of $M$, and let $\epsilon > 0$ be given.
Because $E$ is assumed to be bounded, $HF^\alpha_{con}(E)$ is finite,
and there is a finite collection of nonempty subsets $V_1, \ldots,
V_k$ of $M$ such that $E \subseteq \bigcup_{i=1}^k V_i$ and
\begin{equation}
\label{sum_{i=1}^k (diam V_i)^alpha < HF^alpha_{con}(E) + epsilon}
	\sum_{i=1}^k (\diam V_i)^\alpha < HF^\alpha_{con}(E) + \epsilon.
\end{equation}
For $i = 1, \ldots, k$, let $J_i$ be a closed and bounded interval
in the real line such that $f(V_i) \subseteq J_i$ and 
\begin{equation}
	\diam J_i = \diam f(V_i) \le C \diam V_i.
\end{equation}
Consider the function
\begin{equation}
	\phi(x) = \sum_{i=1}^k (\diam V_i)^{\alpha - 1} \, {\bf 1}_{J_i}(x)
\end{equation}
for $x \in {\bf R}$.  For each $x \in {\bf R}$, $E \cap f^{-1}(x)$ is
a subset of $M$ which is contained in the union of the $V_i$'s such
that $x \in J_i$.  Hence
\begin{equation}
	HF^{\alpha - 1}_{con}(E \cap f^{-1}(x)) \le \phi(x)
\end{equation}
for each $x \in {\bf R}$.  The analogue of (\ref{sum_{p=1}^m a_p |J_p|})
in this case is 
\begin{eqnarray}
	\sum_{i=1}^k (\diam V_i)^{\alpha - 1} \, |J_i|		
		& \le & C \, \sum_{i=1}^k (\diam V_i)^\alpha	\\ 
		& < & C \cdot HF^\alpha_{con}(E) + C \cdot \epsilon.
								 \nonumber
\end{eqnarray}
From (\ref{chebychev's inequality, 1}) it follows that for each $t > 0$,
\begin{eqnarray}
\label{HF^1_{con}({HF^{alpha - 1}_{con}(E cap f^{-1}(x)) > t }) le ..., 1}
\lefteqn{HF^1_{con}(\{x \in {\bf R} : 
		HF^{\alpha - 1}_{con}(E \cap f^{-1}(x)) > t \})} \\
	& & \le t^{-1} (C \cdot HF^\alpha_{con}(E) + C \cdot \epsilon).
								\nonumber
\end{eqnarray}
Since this is true for each $\epsilon > 0$, we get that
\begin{eqnarray}
\label{HF^1_{con}({HF^{alpha - 1}_{con}(E cap f^{-1}(x)) > t }) le ..., 2}
\lefteqn{HF^1_{con}(\{x \in {\bf R} : 
		HF^{\alpha - 1}_{con}(E \cap f^{-1}(x)) > t \})} \\
	& & \le C \, t^{-1} \, HF^\alpha_{con}(E).
								\nonumber
\end{eqnarray}

	Similarly, suppose that $E$ is a nonempty subset of $M$ such
that $H^\alpha_{con}(E) < \infty$, and let $\epsilon > 0$ be given.
There is an at most countable collection $\{V_i\}_{i \in I}$ of
nonempty subsets of $M$ such that $E \subseteq \bigcup_{i \in I} V_i$
and
\begin{equation}
\label{sum_{i in I} (diam V_i)^alpha < H^alpha_{con}(E) + epsilon}
	\sum_{i \in I} (\diam V_i)^\alpha < H^\alpha_{con}(E) + \epsilon.
\end{equation}
Of course 
\begin{equation}
	\diam f(V_i) \le C \diam V_i
\end{equation}
for each $i \in I$.  Choose a bounded open interval $J_i$ in ${\bf R}$
for each $i \in I$ so that $f(V_i) \subseteq J_i$ and
\begin{equation}
\label{sum (diam V_i)^{alpha - 1} |J_i| < C H^alpha_{con}(E) + C epsilon}
	\sum_{i \in I} (\diam V_i)^{\alpha - 1} \, |J_i|
		< C \cdot H^\alpha_{con}(E) + C \cdot \epsilon.
\end{equation}
In other words, the length of $J_i$ may have to be strictly larger
than $\diam f(V_i)$ in order for $J_i$ to be an open interval which
contains $f(V_i)$, and we ask that it be only a little bit larger in
this way.  Consider the function
\begin{equation}
	\phi(x) = \sum_{i \in I} (\diam V_i)^{\alpha - 1} \, {\bf 1}_{J_i}(x)
\end{equation}
for $x \in {\bf R}$.  As before, if $x \in {\bf R}$, then $E \cap
f^{-1}(x)$ is contained in the union of the $V_i$'s such that $x \in
J_i$, and thus
\begin{equation}
	H^{\alpha - 1}_{con}(E \cap f^{-1}(x)) \le \phi(x).
\end{equation}
Using (\ref{chebychev's inequality, 2}) and (\ref{sum (diam
V_i)^{alpha - 1} |J_i| < C H^alpha_{con}(E) + C epsilon}) we get that
\begin{eqnarray}
\label{H^1_{con}({H^{alpha - 1}_{con}(E cap f^{-1}(x)) > t }) le ..., 1}
\lefteqn{H^1_{con}(\{x \in {\bf R} : 
			H^{\alpha - 1}_{con}(E \cap f^{-1}(x)) > t \})}  \\
	&& \le t^{-1} (C \cdot H^\alpha_{con}(E) + C \cdot \epsilon).
								\nonumber
\end{eqnarray}
Because this holds for all $\epsilon > 0$, we obtain
\begin{eqnarray}
\label{H^1_{con}({H^{alpha - 1}_{con}(E cap f^{-1}(x)) > t }) le ..., 2}
\lefteqn{H^1_{con}(\{x \in {\bf R} : 
			H^{\alpha - 1}_{con}(E \cap f^{-1}(x)) > t \})}  \\
	&& \le C \, t^{-1} \, H^\alpha_{con}(E).
								\nonumber
\end{eqnarray}

\section{Topological and Hausdorff dimensions}
\label{section about topological and Hausdorff dimensions}

	Two general references are \cite{H-W, Nagata}.

	Let $(M, d(x,y))$ be a metric space.  A subset $W$ of $M$ is
said to be \emph{clopen}\index{clopen subset of a metric space} if
$W$ is both open and closed.  Thus the empty set and $M$ itself are
both clopen subsets of $M$, and $M$ has a nonempty proper subset which
is clopen if and only if $M$ is not connected.  Observe that finite
unions and intersections of clopen sets are clopen.

	We say that $M$ has topological dimension $0$ if for each $x$
in $M$ and $\epsilon > 0$ there is a clopen subset $W$ of $M$ such
that $x$ is an element of $W$ and $W$ is a subset of $B(x, \epsilon)$.
Clearly, if two metric spaces are homeomorphic to each other and one
of them has topological dimension $0$, then so does the other.

	If $E$ is a nonempty subset of $M$, then we say that $E$ has
topological dimension $0$ if it has topological dimension $0$ when
viewed as a metric space in its own right, using the restriction of
the metric $d(x,y)$ to $E$.  It is convenient to say that the empty
set has topological dimension equal to $-1$.  Notice that if a
nonempty subset $E$ of $M$ has topological dimension $0$ and $E_1$ is
a nonempty subset of $E$, then $E_1$ also has topological dimension
$0$.  As a basic class of examples, a nonempty subset of the real line
has topological dimension $0$ if and only if it does not contain an
interval of positive length.

	It is easy to see that if $M$ has topological dimension $0$,
then $M$ is totally disconnected, i.e., $M$ does not contain a
connected subset with at least two elements.  The converse holds when
$M$ is compact, or merely locally compact.  This can be derived from
Lemma \ref{compact and totally disconnected spaces}.  It is known that
there are nonempty bounded subsets of the plane which are totally
disconnected and which become connected after the addition of a single
point.  Such a set does not have topological dimension $0$, as can
be shown using results described later in this section.

	Suppose that $E$ is a nonempty closed subset of $M$.  In this
case, $E$ has topological dimension $0$ if and only if for every $x$
in $E$ and every $\epsilon > 0$ there are closed subsets $F_1$, $F_2$
of $E$ such that $x$ is an element of $F_1$, $F_1$ is a subset of
$B(x, \epsilon)$, $F_1 \cap F_2 = \emptyset$, and $F_1 \cup F_2 = E$.

\beginlemma
\label{clopen neighborhoods of compact subsets of spaces with top dim 0}
Let $(M, d(x,y))$ be a metric space with topological dimension $0$.
Suppose that $K$ is a compact subset of $M$ and that $U$ is an open
subset of $M$ such that $K \subseteq U$.  Then there is a clopen subset
$W$ of $M$ such that $K \subseteq W \subseteq U$.
\end{lemma}

	Indeed, the assumption that $M$ has topological dimension $0$
implies that for each point $p$ in $K$ there is a clopen subset $W(p)$
of $M$ such that $p \in W(p)$ and $W(p) \subseteq U$.  Because $K$ is
compact, there exist a finite collection $p_1, \ldots, p_n$ of points
in $K$ such that $K \subseteq \bigcup_{i=1}^n W(p_i)$.  If we take $W
= \bigcup_{i=1}^n W(p_i)$, then $W$ is clopen and $K \subseteq W
\subseteq U$, as desired.

\beginproposition
\label{separating disjoint closed sets in a separable space with top dim 0}
Let $(M, d(x,y))$ be a separable metric space with topological dimension $0$.
If $F_1$, $F_2$ are disjoint closed subsets of $M$, then there are
clopen subsets $V_1$, $V_2$ of $M$ such that $F_1 \subseteq V_1$,
$F_2 \subseteq V_2$, $V_1 \cap V_2 = \emptyset$, and $V_1 \cup V_2 = M$.
\end{proposition}

	Notice that if $F_1$ or $F_2$ is compact, then this follows
from Lemma \ref{clopen neighborhoods of compact subsets of spaces with
top dim 0}.  

	Because $M$ has topological dimension $0$, one can choose for
each $p$ in $M$ a clopen subset $U(p)$ of $M$ such that $p$ is an
element of $U(p)$ and either $U(p) \cap F_1 = \emptyset$ or $U(p) \cap
F_2 = \emptyset$.  The assumption that $M$ be separable implies that
$M$ is countably compact, and hence there is a sequence
$\{p_j\}_{j=1}^\infty$ of points in $M$ such that $M =
\bigcup_{j=1}^\infty U(p_j)$.

	Define subsets $W_j$ of $M$ for $j \ge 1$ by $W_1 = U(p_1)$,
$W_j$ is equal to $U(p_j)$ minus the union of the $U(p_i)$'s for $1
\le i < j$.  Thus for each $j$ we have that $W_j$ is a clopen subset
of $M$, $W_j \subseteq U(p_j)$,
\begin{equation}
	\bigcup_{i=1}^j W_i = \bigcup_{i=1}^j U(p_i),
\end{equation}
and therefore
\begin{equation}
	\bigcup_{i=1}^\infty W_i = M.
\end{equation}
The $W_j$'s are also pairwise disjoint, by construction.  Let $V_1$ be
the union of the $W_i$'s which intersect $F_1$, and let $V_2$ be the
union of all the other $W_i$'s.  Then $V_1 \cap V_2 = \emptyset$, $V_1
\cup V_2 = M$, $F_1 \subseteq V_1$, and $F_2 \subseteq V_2$.  Because
the $W_i$'s are all clopen, $V_1$ and $V_2$ are both open sets, and
hence they are closed sets as well.  This proves the proposition.

	Let us reformulate the proposition for closed subsets of a
metric space as follows.  If $(M, d(x,y))$ is a separable metric space
and $E$ is a nonempty closed subset of $M$ with topological dimension
$0$, and if $F_1$, $F_2$ are disjoint closed subsets of $E$, then
there are disjoint closed subsets $\widehat{F}_1$, $\widehat{F}_2$ of
$M$ such that $F_1 \subseteq \widehat{F}_1$, $F_2 \subseteq
\widehat{F}_2$, and $\widehat{F}_1 \cup \widehat{F}_2 = E$. 

	Here is another formulation which is sometimes more
convenient.  Suppose that $(M, d(x,y))$ is a separable metric space,
$E$ is a nonempty closed subset of $M$ with topological dimension $0$,
and that $F_1$, $F_2$ are disjoint closed subsets of $M$.  Then there
are disjoint closed subsets $\widetilde{F}_1$, $\widetilde{F}_2$ of
$M$ such that $F_1 \subseteq \widetilde{F}_1$, $F_2 \subseteq
\widetilde{F}_2$, and $E \subseteq \widetilde{F}_1 \cup
\widetilde{F}_2$.

	Now suppose that $(M, d(x,y))$ is a separable metric space,
and that $E_1$, $E_2$ are two nonempty closed subsets of $M$ with
topological dimension $0$.  Let us check that the union $E_1 \cup E_2$
has topological dimension $0$ too.  Let $A$, $B$ be disjoint closed
subsets of $E_1 \cup E_2$.  Using the result mentioned in the previous
paragraph, one can get disjoint closed subsets of $M$ which contain
$A$, $B$, respectively, and whose union contains $E_1 \cup E_2$.  More
precisely, one first enlarges $A$, $B$ to disjoint closed subsets
$A_1$, $B_1$ of $M$ such that $A \subseteq A_1$, $B \subseteq B_1$,
and $E_1 \subseteq A_1 \cup B_1$.  Then one enlarges $A_1$, $B_1$ to
get disjoint closed subsets $A_2$, $B_2$ of $M$ such that $A_1
\subseteq A_2$, $B_1 \subseteq B_2$, and $E_2 \subseteq A_2 \cup B_2$.
Thus $A_2$, $B_2$ are disjoint closed subsets of $M$ which satisfy $A
\subseteq A_2$, $B \subseteq B_2$, and $E_1 \cup E_2 \subseteq A_2
\cup B_2$.  From this it follows that $E_1 \cup E_2$ has topological
dimension $0$.

	In fact there is a version of this for sequences of closed
sets with topological dimension $0$.  Specifically, let $(M, d(x,y))$
be a separable metric space, and let $\{E_j\}_{j=1}^\infty$ be a sequence
of nonempty closed subsets of $M$ with topological dimension $0$.
Under these conditions, $\bigcup_{j=1}^\infty E_j$ has topological
dimension $0$ as well.  Without loss of generality, we may assume
that $M = \bigcup_{j=1}^\infty E_j$.  Let $A$, $B$ be disjoint
closed subsets of $M$.  We would like to show that there are disjoint
closed subsets of $M$ which contain $A$, $B$, respectively, and whose
union is equal to $M$.

	Basically, we would like to repeat the process as before, but
now it is not clear that the union over all $j \ge 1$ leads to closed
subsets of $M$.  This can be fixed using normality of $M$.  Because
$E_1$ has topological dimension $0$, there are disjoint closed subsets
$A_1$, $B_1$ of $M$ such that $A \subseteq A_1$, $B \subseteq B_1$,
and $E_1 \subseteq A_1 \cup B_1$.  By the normality of $M$, there are
disjoint closed subsets $\widetilde{A}_1$, $\widetilde{B}_1$ of $M$
such that $A_1$ is contained in the interior of $\widetilde{A}_1$ and
$B_1$ is contained in the interior of $\widetilde{B}_1$.  Next, one
uses the hypothesis that $E_2$ has topological dimension $0$ to obtain
disjoint closed subsets $A_2$, $B_2$ of $M$ such that $\widetilde{A}_1
\subseteq A_2$, $\widetilde{B}_1 \subseteq B_2$, and $E_2 \subseteq
A_2 \cup B_2$.  Again normality can be applied to get disjoint closed
subsets of $M$ which contain $A_2$, $B_2$ in their interiors.
Repeating the process we obtain for each $j$ closed subsets $A_j$,
$B_j$, $\widetilde{A}_j$, $\widetilde{B}_j$ such that $A_j \cap B_j =
\emptyset$, $\widetilde{A}_j \cap \widetilde{B}_j = \emptyset$, and
$E_j \subseteq A_j \cup B_j$ for all $j$, $A_j$, $B_j$ are contained in
the interiors of $\widetilde{A}_j$, $\widetilde{B}_j$, respectively,
$\widetilde{A}_j \subseteq A_{j+1}$ and $\widetilde{B}_j \subseteq
B_{j+1}$, and $A \subseteq A_1$, $B \subseteq B_1$.  Thus
\begin{equation}
	\bigcup_{j=1}^\infty A_j = \bigcup_{j=1}^\infty \widetilde{A}_j,
		\qquad
	\bigcup_{j=1}^\infty B_j = \bigcup_{j=1}^\infty \widetilde{B}_j.
\end{equation}
If we set $U = \bigcup_{j=1}^\infty A_j$ and $V = \bigcup_{j=1}^\infty
B_j$, then we have that $A \subseteq U$, $B \subseteq V$, $U \cap V =
\emptyset$, $U \cup V = M$, and $U$, $V$ are both open subsets of $M$.
Hence they are also closed, and it follows that $M$ has topological
dimension $0$.

	A nice consequence of this is that if $E_1$, $E_2$ are
nonempty subsets of $M$ with topological dimension $0$ and $E_1 \cup
E_2 = M$, and if $E_1$ is closed, then $M$ has topological dimension
$0$.  We already saw an argument for this when both $E_1$ and $E_2$
are closed, and it does not work for arbitrary subsets of $M$, since
the real line is the union of the rationals and irrationals, each of
which have topological dimension $0$.  If $E_1$ is closed, then $M
\backslash E_1$ is open, and it also has topological dimension $0$,
since it is a subset of $E_2$.  An open set can be expressed as the
union of a sequence of closed sets, and it follows that $M$ can be
expressed as the union of a sequence of closed subsets with
topological dimension $0$, and hence has topological dimension $0$
itself.

	Now suppose that $M$ has $1$-dimensional Hausdorff content
equal to $0$.  By countable subadditivity, this holds if $M$ can be
expressed as the union of a sequence of subsets with $1$-dimensional
Hausdorff content equal to $0$.  If $f$ is a real-valued function on
$M$ which is Lipschitz, or countably-Lipschitz, then $f(M)$ is a
subset of the real line with $1$-dimensional Hausdorff content equal
to $0$, and in particular $f(M)$ does not contain an interval of
positive length.  By choosing $f(x)$ to be of the form $d(p, x)$,
where $p$ is any element of $M$, it follows that there are plenty of
balls in $M$ which are clopen sets, and hence that $M$ has topological
dimension $0$.  One can also use this to find plenty of disjoint
closed sets whose union is equal to $M$ and which contain a given pair
of disjoint closed subsets of $M$.

	If $n$ is a positive integer, then $M$ is said to have
topological dimension $\le n$ if for every point $x$ in $M$ and every
$\epsilon > 0$ there is an open subset $U$ of $M$ such that $x \in U$,
$U \subseteq B(x, \epsilon)$, and $\partial U$ has topological dimension
$\le n-1$, as a metric space itself.  In other words, the property of
having topological dimension $\le n$ is defined inductively on $n$.

	Suppose that $M$ has $(n+1)$-dimensional Hausdorff content
equal to $0$, and that $f : M \to N$ is Lipschitz, or countably
Lipschitz.  Then $f^{-1}(t)$ has $n$-dimensional Hausdorff content
equal to $0$ for all $t \in {\bf R}$ except for a set of
$1$-dimensional Hausdorff content equal to $0$.  In particular,
this holds for a dense set of $t$'s.  Using this it is easy to see
that $M$ has topological dimension $\le n$ in this case.

\chapter{Miscellaneous, 2}
\label{chapter -- miscellaneous, 2}

\section{Some geometry and analysis on the real line}
\label{section on some geometry and analysis on the real line}

	Let $\mathcal{E}$ be a nonempty finite collection of intervals
in the real line.  Now, if $J_1$, $J_2$, and $J_3$ are intervals
in ${\bf R}$ which contain at least one point in common, then it is
easy to see that one of these three intervals is contained in the other
two.  Using this simple observation, it follows that there is a subcollection
$\mathcal{E}_1$ of $\mathcal{E}$ such that
\begin{equation}
	\bigcup_{J \in \mathcal{E}_1} J = \bigcup_{J \in \mathcal{E}} J
\end{equation}
and any element of ${\bf R}$ is contained in at most $2$ intervals
in $\mathcal{E}$.

	Now suppose that $\mu(x)$ is a monotone increasing function
on the real line which is bounded, and set
\begin{equation}
	A = \inf_{x \in {\bf R}} \mu(x), \quad B = \sup_{x \in {\bf R}} \mu(x).
\end{equation}
Because $\mu$ is monotone increasing, we have that $\mu(x) \to A$ as
$x \to -\infty$ and $\mu(x) \to B$ as $x \to +\infty$.  Thus, for each
$\epsilon > 0$, there are real numbers $L$, $N$ such that
\begin{equation}
	\mu(x) \le A + \epsilon
\end{equation}
when $x \le L$ and
\begin{equation}
	\mu(x) \ge B - \epsilon
\end{equation}
when $x \ge N$.  For each $x \in {\bf R}$, put
\begin{equation}
	\mu(x-) = \lim_{y \to x-} \mu(y), \qquad 
		\mu(x+) = \lim_{y \to x+} \mu(x+).
\end{equation}
It will be convenient to set $\mu(-\infty) = A$ and
$\mu(+\infty) = B$.

	Let $J$ be an interval in ${\bf R}$, and let us define the
$\mu$-length of $J$, denoted $\mu(J)$, as follows.  First, if $J$ is
of the form $(u, v)$, then we set
\begin{equation}
	\mu(J) = \mu(v-) - \mu(u+).
\end{equation}
Second, if $J$ is of the form $[u, v)$, then we set
\begin{equation}
	\mu(J) = \mu(v-) - \mu(u-).
\end{equation}
Third, if $J$ is of the form $(u, v]$, then we set
\begin{equation}
	\mu(J) = \mu(v+) - \mu(u+).
\end{equation}
Fourth, if $J$ is of the form $[u, v]$, then we set
\begin{equation}
	\mu(J) = \mu(v+) - \mu(u-).
\end{equation}
In this case it may be that $u = v$, so that $J$ consists of a single
point.  When $J$ consists of a single point, $\mu(J) = 0$ if and only
if $\mu$ is continuous at that point.  

	Notice that
\begin{equation}
	\mu({\bf R}) = B - A.
\end{equation}
If $J_1, \ldots, J_l$ are disjoint intervals in ${\bf R}$, then it is
not hard to check that
\begin{equation}
	\sum_{i=1}^l \mu(J_i) \le B - A.
\end{equation}
If instead $J_1, \ldots, J_l$ are intervals in ${\bf R}$ such that
each point in ${\bf R}$ is contained in at most two of the $J_i$'s, then
\begin{equation}
	\sum_{i=1}^l \mu(J_i) \le 2 (B - A).
\end{equation}
Again, this is not too difficult to verify.

	Fix a real number $\alpha$ such that $0 < \alpha \le 1$, and
define a function $\mu^*_\alpha(x)$ on ${\bf R}$ by
\begin{eqnarray}
\label{def of mu^*_alpha}
\lefteqn{\qquad \mu^*_\alpha(x) = }					\\
	& & \sup \biggl\{ \frac{\mu(J)}{|J|^\alpha} :
				\ J \hbox{ is a bounded open interval in }
					{\bf R} \hbox{ such that } x \in J\},
							\nonumber
\end{eqnarray}
where $|J|$ denotes the ordinary length of $J$, as usual.  Thus
$\mu^*_\alpha(x) \ge 0$, and $\mu^*_\alpha(x) = +\infty$ is possible,
although it does not happen too much, as we shall see.

	If $\mu(x)$ is constant, then $\mu(J) = 0$ for all intervals
$J$ in ${\bf R}$, and $\mu^*_\alpha(x) = 0$ for all $x$.  Let us assume
that $\mu(x)$ is not constant, which implies that $\mu^*_\alpha(x) > 0$
for all $x \in {\bf R}$.

	Fix a positive real number $t$, and consider the set
\begin{equation}
\label{{x in {bf R} : mu^*_alpha(x) > t}}
	\{x \in {\bf R} : \mu^*_\alpha(x) > t\}.
\end{equation}
If $\mu^*_\alpha(x) > t$, then there is a bounded open interval $J$
in ${\bf R}$ such that $x \in J$ and 
\begin{equation}
\label{frac{mu(J)}{|J|^alpha} > t}
	\frac{\mu(J)}{|J|^\alpha} > t.
\end{equation}
Thus $\mu^*_\alpha(y) > t$ for all $y \in J$ in these circumstances,
so that $J$ is contained in the set (\ref{{x in {bf R} :
mu^*_alpha(x) > t}}).  This shows that (\ref{{x in {bf R} :
mu^*_alpha(x) > t}}) is an open subset of ${\bf R}$.  Also,
\begin{equation}
\label{frac{mu(J)}{|J|^alpha} > t implies |J| le t^{-1/alpha} (B - A)^alpha}
	\frac{\mu(J)}{|J|^\alpha} > t  \ \hbox{ implies that } \ 
		|J| \le t^{-1/\alpha} \, (B - A)^\alpha.
\end{equation}

	Suppose that $K$ is a compact subset of 
(\ref{{x in {bf R} : mu^*_alpha(x) > t}}).  Thus each element of
$K$ is contained in a bounded open interval $J$ such that
$|J|^\alpha < t^{-1} \, \mu(J)$.  Because of compactness, $K$
is contained in the union of finitely many such intervals.
Using the remark at the beginning of the section, there
are finitely many bounded open intervals $J_1, \ldots, J_l$
in ${\bf R}$ such that 
\begin{equation}
	K \subseteq \bigcup_{i=1}^l J_i,
\end{equation}
each point in ${\bf R}$ is contained in at most two of the $J_i$'s,
and 
\begin{equation}
	|J_i|^\alpha < t^{-1} \, \mu(J_i)
\end{equation}
for each $i$.  Thus
\begin{equation}
\label{HF^alpha_{con}(K) le ... le 2 t^{-1} (B - A)}
	HF^\alpha_{con}(K) \le \sum_{i=1}^l |J_i|^\alpha
		< \sum_{i=1}^l t^{-1} \, \mu(J_i)
		\le 2 \, t^{-1} \, (B - A).
\end{equation}

	Suppose that $\alpha = 1$, and let us check that
\begin{equation}
\label{H^1({x in {bf R} : mu^*_1(x) > t}) le 2 t^{-1} (B - A)}
	H^1(\{x \in {\bf R} : \mu^*_1(x) > t\}) \le 2 \, t^{-1} \, (B - A).
\end{equation}
We know that (\ref{{x in {bf R} : mu^*_alpha(x) > t}}) is an open
subset of the real line, so that it is either empty or it can be
expressed as the union of an at most countable family of disjoint open
intervals in ${\bf R}$, and to show (\ref{H^1({x in {bf R} : mu^*_1(x)
> t}) le 2 t^{-1} (B - A)}) it is enough to show that the sum of the
lengths of these intervals is less than or equal to $2 \, t^{-1} \, (B
- A)$.  It suffices to show that the sum of the lengths of any finite
collection of these open intervals is less than or equal to $2 \,
t^{-1} \, (B - A)$, and to get this it is enough to show that the sum
of the lengths of any finite collection of pairwise disjoint closed
and bounded intervals contained in (\ref{{x in {bf R} : mu^*_alpha(x)
> t}}) is less than or equal to $2 \, t^{-1} \, (B - A)$.  This last
assertion follows from (\ref{HF^alpha_{con}(K) le ... le 2 t^{-1} (B -
A)}) with $\alpha = 1$.

	Let $E$ be the set of $x \in {\bf R}$ such that
\begin{eqnarray}
	&& \liminf_{y \to x-} \frac{\mu(x) - \mu(y)}{x - y}, \quad
		\limsup_{y \to x-} \frac{\mu(x) - \mu(y)}{x - y}, \\
	&& \liminf_{y \to x+} \frac{\mu(y) - \mu(x)}{y - x}, \quad
		\limsup_{y \to x+} \frac{\mu(y) - \mu(x)}{y - x}
							\nonumber
\end{eqnarray}
are finite and equal to each other.  One can show that
\begin{equation}
	H^1_{con}(E) = 0.
\end{equation}
In other words, $\mu$ is differentiable almost everywhere on ${\bf R}$,
in the sense that the set of exceptions has $1$-dimensional Hausdorff
content equal to $0$.

\section{Products of nonnegative linear functionals}
\label{section on products of nonnegative linear functionals}
\index{products of nonnegative linear functionals}

	Let $(M_1, d_1(\cdot, \cdot)$ and $(M_2, d_2(\cdot, \cdot))$
be metric spaces.  Consider the metric space which is the Cartesian
product $M_1 \times M_2$, consisting of ordered pairs $(x, y)$
with $x \in M_1$ and $y \in M_2$, equipped with the metric
\begin{equation}
	\rho((x,y), (x',y')) = d_1(x,x') + d_2(y,y').
\end{equation}
There are other natural choices for the metric on $M_1 \times M_2$,
such as the maximum of $d_1(x,x')$, $d_2(y,y')$, and which are
equivalent in terms of the topology that they determine.  

	Notice that if $U_1$, $U_2$ are open subsets of $M_1$, $M_2$,
respectively, then $U_1 \times U_2$ is an open subset of $M_1 \times
M_2$.  In fact these subsets of $M_1 \times M_2$ form a basis for the
topology of $M_1 \times M_2$.  Similarly, the Cartesian product of two
closed subsets of $M_1$, $M_2$ is a closed subset of $M_1 \times M_2$.
If $\{(x_j, y_j)\}_{j=1}^\infty$ is a sequence of points in $M_1
\times M_2$ and $(x, y)$ is another point in $M_1 \times M_2$, then
$\{(x_j, y_j)\}_{j=1}^\infty$ converges to $(x, y)$ in $M_1 \times
M_2$ if and only if $\{x_j\}_{j=1}^\infty$ converges to $x$ in $M_1$
and $\{y_j\}_{j=1}^\infty$ converges to $y$ in $M_2$.  If $K_1$, $K_2$
are compact subsets of $M_1$, $M_2$, respectively, then $K_1 \times
K_2$ is a compact subset of $M_1 \times M_2$.

	As in Section \ref{section about general forms of
integration}, let us assume that closed and bounded subsets of $M_1$
and $M_2$ are compact.  This implies that closed and bounded subsets
of $M_1 \times M_2$ are compact as well.  Suppose that $\lambda_1$,
$\lambda_2$ are nonnegative linear functionals on $C_{00}(M_1)$,
$C_{00}(M_2)$, respectively.  We want to define a corresponding
nonnegative linear functional $\lambda_1 \times \lambda_2$ on
$C_{00}(M_1 \times M_2)$.

	What should $\lambda_1 \times \lambda_2$ be like?  A basic
property is that if $\phi_1$, $\phi_2$ are functions in $C_{00}(M_1)$,
$C_{00}(M_2)$, respectively, then $\phi(x, y) = \phi_1(x) \,
\phi_2(y)$ is a function in $C_{00}(M_1 \times M_2)$, and we would
like to have
\begin{equation}
\label{basic property of lambda_1 times lambda_2}
	(\lambda_1 \times \lambda_2)(\phi) 
		= \lambda_1(\phi_1) \, \lambda_2(\phi_2).
\end{equation}
In fact this property characterizes $\lambda_1 \times \lambda_2$ uniquely,
because linear combinations of functions of the form
\begin{equation}
\label{phi(x, y) = phi_1(x) phi_2(y), phi_1 in C_{00}(M_1), phi_2 in C_{00}(M_2)}
	\phi(x, y) = \phi_1(x) \, \phi_2(y), \quad
			   \phi_1 \in C_{00}(M_1), \ \phi_2 \in C_{00}(M_2),
\end{equation}
are dense in $C_{00}(M_1 \times M_2)$ with respect to restricted
convergence.  

	We can define $\lambda_1 \times \lambda_2$ in a more direct
manner as follows.  Suppose that $\phi(x, y)$ is any function in
$C_{00}(M_1 \times M_2)$.  For each $y \in M_2$, we can consider 
$\phi(x, y)$ as a function of $x$ on $M_1$, which is an element of
$C_{00}(M_1)$.  Let us write this as $\phi_{1, y}(x)$.  Thus
\begin{equation}
	\lambda_1(\phi_{1, y})
\end{equation}
is a real number for each $y \in M_2$, and it is not hard to see
that this defines a function on $M_2$ which lies in $C_{00}(M_2)$.
Therefore, we can apply $\lambda_2$ to this function to get a real
number, and we can define $(\lambda_1 \times \lambda_2)(\phi)$ to
be this number.  One can easily verify that this defines a nonnegative
linear functional on $C_{00}(M_1 \times M_2)$ which satisfies the
condition (\ref{basic property of lambda_1 times lambda_2}).

	We can just as well do this in the other order.  If $\phi(x,
y)$ is a function in $C_{00}(M_1 \times M_2)$, then $\phi_{2, x}(y) =
\phi(x, y)$ is a function in $C_{00}(M_2)$ for each $x \in M_1$.  As a
result, we can apply $\lambda_2$ to this function for each $x \in M_1$
to obtain a real number $\lambda_2(\phi_{2, x})$, and this defines a
function in $C_{00}(M_1)$.  We can apply $\lambda_1$ to this function
to get a real number, which provides an alternate approach to defining
$(\lambda_1 \times \lambda_2)(\phi)$.  

	The approaches to defining $\lambda_1 \times \lambda_2$
described in the preceding two paragraphs give the same answer,
because they both define nonnegative linear functionals on $C_{00}(M_1
\times M_2)$ which agree on linear combinations of functions of the
form (\ref{phi(x, y) = phi_1(x) phi_2(y), phi_1 in C_{00}(M_1), phi_2
in C_{00}(M_2)}).  This statement is analogous to the familiar
fact in advanced calculus that double integrals can be given in terms
of iterated integrals in either order.

	Note that if $\lambda_1$, $\lambda_2$ are bounded nonnegative
linear functionals on $C_{00}(M_1)$, $C_{00}(M_2)$, respectively, then
the product $\lambda_1 \times \lambda_2$ is a bounded nonnegative
linear functional on $C_{00}(M_1 \times M_2)$.  More precisely,
if $C_1$, $C_2$ are nonnegative real numbers such that
\begin{equation}
	|\lambda_1(\phi_1)| \le C_1 \, \|\phi_1\|_{sup, M_1}
\end{equation}
for all $\phi_1$ in $C_{00}(M_1)$, and
\begin{equation}
	|\lambda_2(\phi_2)| \le C_2 \, \|\phi_2\|_{sup, M_2}
\end{equation}
for all $\phi_2$ in $C_{00}(M_2)$, then
\begin{equation}
	|(\lambda_1 \times \lambda_2)(\phi)| 
		\le C_1 \, C_2 \, \|\phi\|_{sup, M_1 \times M_2}
\end{equation}
for all $\phi$ in $C_{00}(M_1 \times M_2)$.  Here we write
$\|f\|_{sup, M}$ for the supremum norm of a real-valued function $f$
on a metric space $M$.

	If $\lambda_1$, $\lambda_2$ are bounded, so that $\lambda_1
\times \lambda_2$ is also bounded, then $\lambda_1$, $\lambda_2$, and
$\lambda_1 \times \lambda_2$ have natural extensions to bounded
nonnegative linear functionals on $C_b(M_1)$, $C_b(M_2)$, and $C_b(M_1
\times M_2)$, respectively, as in Section \ref{section about general
forms of integration}.  For these extensions we continue to have that
\begin{equation}
	(\lambda_1 \times \lambda_2)(\phi) 
		= \lambda_1(\phi_1) \, \lambda_2(\phi_2)
\end{equation}
for $\phi_1$ in $C_b(M_1)$ and $\phi_2$ in $C_b(M_2)$, where
$\phi(x, y) = \phi_1(x) \, \phi_2(y)$.

	Suppose that $\lambda_1$, $\lambda_2$ are nonnegative linear
functionals on $C_{00}(M_1)$, $C_{00}(M_2)$, with supports equal to
$F_1 \subseteq M_1$, $F_2 \subseteq M_2$, respectively.  One can check
that the support of $\lambda_1 \times \lambda_2$ in $M_1 \times M_2$
is equal to $F_1 \times F_2$.  In particular, if $\lambda_1$,
$\lambda_2$ have compact supports in $M_1$, $M_2$, then $\lambda_1
\times \lambda_2$ has compact support in $M_1 \times M_2$.

\section{Some harmonic analysis}
\label{section on some harmonic analysis}

	Fix a positive integer $n$, and let us consider ${\bf R}^n$
with its usual metric.  

	Let $\lambda$ be a nonnegative linear functional on
$C_{00}({\bf R}^n)$, and let $f$ be a function in $C_{00}({\bf R}^n)$.
For each $y \in {\bf R}^n$, put
\begin{equation}
	f_y(x) = f(y - x).
\end{equation}
Consider $\lambda(f_y)$, as a real-valued function of $y$ on ${\bf
R}^n$.  This function is called the \emph{convolution of $\lambda$ and
$f$},\index{convolutions} and we shall denote it $(\lambda * f)(y)$.
It is not difficult to show that $\lambda(f_y)$ is continuous as a
function of $y$.

	If $\lambda$ has compact support in ${\bf R}^n$, then the
convolution $(\lambda * f)(y)$ can be defined for any continuous
real-valued function $f$ on ${\bf R}^n$, and one can check that
$(\lambda * f)(y)$ is again continuous.  If both $\lambda$ and $f$
have compact support in ${\bf R}^n$, then the convolution $\lambda *
f$ has compact support in ${\bf R}^n$.

	When $\lambda$ is a bounded nonnegative linear functional
on $C_{00}({\bf R}^n)$, we have a natural extension of $\lambda$ to
a bounded nonnegative linear functional on $C_b({\bf R}^n)$, as in
Section \ref{section about general forms of integration}.  In this
case the convolution $(\lambda * f)(y)$ can be defined for all
bounded continuous real-valued functions $f$ on ${\bf R}^n$, and
one can check that the convolution is again a bounded continuous
function on ${\bf R}^n$.  More precisely, if $\lambda$ is bounded
with constant $C \ge 0$, then 
\begin{equation}
	\|\lambda * f \|_{sup} \le C \, \|f\|_{sup}.
\end{equation}

	Let $\lambda$ be a bounded nonnegative linear functional on
$C_{00}({\bf R}^n)$, with its natural extension to a bounded
nonnegative linear functional on $C_b({\bf R}^n)$.  We can extend
$\lambda$ to a functional defined on complex-valued functions on $M$,
by setting
\begin{equation}
	\lambda(f_1 + i \, f_2) = \lambda(f_1) + i \, \lambda(f_2)
\end{equation}
when $f_1$, $f_2$ are real-valued functions on $M$.  For each $w \in
{\bf R}^n$, let $E^w(x)$ be the complex exponential function on ${\bf
R}^n$ defined by
\begin{equation}
	E^w(x) = \exp (2 \pi i x \cdot w),
\end{equation}
where $x \cdot w$ is the usual inner product on ${\bf R}^n$,
\begin{equation}
	x \cdot w = \sum_{j=1}^n x_j \, w_j.
\end{equation}
Recall that 
\begin{equation}
	|E^w(x)| = 1 \ \hbox{ and } \ \overline{E^w(x)} = E^w(-x)
\end{equation}
for all $x \in {\bf R}^n$, where $\overline{a}$ denotes the
complex conjugate of a complex number $a$, and
\begin{equation}
	E^w(x + y) = E^w(x) \, E^w(y)
\end{equation}
for all $x, y \in {\bf R}^n$.  Define the \emph{Fourier
transform}\index{Fourier transform} of $\lambda$ to be the
function $\widehat{\lambda}(w)$ on ${\bf R}^n$ given by
\begin{equation}
	\widehat{\lambda}(w) = \lambda(E^{-w}).
\end{equation}
One can check that $\widehat{\lambda}$ is a bounded continuous function
on ${\bf R}^n$, and more precisely that if $\lambda$ is a bounded
nonnegative linear functional with constant $C$, then 
\begin{equation}
	\|\widehat{\lambda}\|_{sup} \le C.
\end{equation}
It is easy to see from the definitions that
\begin{equation}
	\lambda * E^w(y) = \widehat{\lambda}(w) \, E^w(y).
\end{equation}
In other words, $E^w$ is an eigenfunction for the linear operator
\begin{equation}
	f \mapsto \lambda * f
\end{equation}
on $C_b({\bf R}^n)$, with eigenvalue $\widehat{\lambda}(w)$.

	If $\lambda_1$, $\lambda_2$ are nonnegative linear functionals
on $C_{00}({\bf R}^n)$, then we can try to define the
convolution\index{convolutions} $\lambda_1 * \lambda_2$, as a
nonnegative linear functional on $C_{00}({\bf R}^n)$, by
\begin{equation}
	(\lambda_1 * \lambda_2)(f) = (\lambda_1 \times \lambda_2)(F),
			\quad F(x,y) = f(x + y).
\end{equation}
This makes sense if $\lambda_1$, $\lambda_2$ are both bounded, in
which case $\lambda_1 * \lambda_2$ is also a bounded nonnegative
linear functional on $C_{00}({\bf R}^n)$, and in fact one can let $f$
be a bounded continuous real-valued function on ${\bf R}^n$.  To be
more precise, if $\lambda_1$, $\lambda_2$ are bounded with constants
$C_1$, $C_2$, then $\lambda_1 * \lambda_2$ is bounded with constant
$C_1 \, C_2$.  If at least one of $\lambda_1$, $\lambda_2$ has compact
support, and $f$ is a continuous real-valued function on ${\bf R}^n$
with compact support, then one can show that the definition of
$(\lambda_1 * \lambda_2)(f)$ also makes sense.  Observe that when both
$\lambda_1$, $\lambda_2$ have compact support, the convolution
$\lambda_1 * \lambda_2$ has compact support as well, and one can use
any continuous real-valued function $f$ on ${\bf R}^n$ in the formula
above.

	This notion of convolution is compatible with the previous
one in the sense that
\begin{equation}
	(\lambda_1 * \lambda_2) * f = \lambda_1 * (\lambda_2 * f)
\end{equation}
in the situations where the convolutions are defined, as one can
verify.  It is clear from the definition of $\lambda_1 * \lambda_2$
that
\begin{equation}
	\lambda_1 * \lambda_2 = \lambda_2 * \lambda_1.
\end{equation}
If $\lambda_1$, $\lambda_2$ are bounded nonnegative linear functionals
on $C_{00}({\bf R}^n)$, then
\begin{equation}
	(\lambda_1 * \lambda_2)^{\widehat{}} (w)
		= \widehat{\lambda_1}(w) \, \widehat{\lambda_2}(w)
\end{equation}
for all $w \in {\bf R}^n$.

	These notions also apply to the case of the $n$-dimensional
torus ${\bf T}^n = {\bf R}^n / {\bf Z}^n$.  Because ${\bf T}^n$ is
compact, $C({\bf T}^n) = C_{00}({\bf T}^n)$, and a nonnegative linear
functional $\lambda$ on $C({\bf Z}^n)$ is automatically bounded with
constant $C = \lambda(1)$.  If $f$ is any continuous function on ${\bf
T}^n$, then $\lambda * f$ can be defined as before, and is a
continuous function on ${\bf T}^n$.  When $w \in {\bf Z}^n$, the
complex exponential function $E^w$ can be viewed as a function on
${\bf T}^n$, since it is periodic, and we define the Fourier transform
$\widehat{\lambda}$ of $\lambda$ to be the function on ${\bf Z}^n$
given by $\widehat{\lambda}(w) = \lambda(E^{-w})$.  As before,
$\lambda * E^w(y) = \widehat{\lambda}(w) \, E^w(y)$.  If $\lambda_1$,
$\lambda_2$ are two nonnegative linear functionals on $C({\bf T}^n)$,
then the convolution $\lambda_1 * \lambda_2$ can be defined in the
same way as before, to give another nonnegative linear functional on
$C({\bf T}^n)$.  Again we have that $\lambda_1 * \lambda_2 = \lambda_2
* \lambda_1$, $(\lambda_1 * \lambda_2) * f = \lambda_1 * (\lambda_2 *
f)$ for any continuous function $f$ on ${\bf T}^n$, and $(\lambda_1 *
\lambda_2)^{\widehat{}} (w) = \widehat{\lambda}_1(w) \,
\widehat{\lambda}_2(w)$ for all $w \in {\bf Z}^n$.

	Similarly, these notions can be applied to ${\bf Z}^n$ as
well.  Every function on ${\bf Z}^n$ is continuous, and $C_{00}({\bf
Z}^n)$, $C_b({\bf Z}^n)$ simply consist of the functions on ${\bf
Z}^n$ which have finite support or which are bounded, respectively.
The convolution $\lambda * f$ can be defined as before when $\lambda$
is a nonnegative linear functional on $C_{00}({\bf Z}^n)$ and $f$ is a
function on ${\bf Z}^n$ with finite support, or when $\lambda$ is a
bounded nonnegative linear functional on $C_{00}({\bf Z}^n)$ and $f$
is a bounded function on ${\bf Z}^n$, or when $\lambda$ is a
nonnegative linear functional on $C_{00}({\bf Z}^n)$ with finite
support and $f$ is any function on ${\bf Z}^n$.  When $\lambda$ and
$f$ are bounded, so is the convolution $\lambda * f$.  In this case we
can view the complex exponential $E^w(x)$ as defining a function of
$x$ on ${\bf Z}^n$ for all $w \in {\bf T}^n$.  This is because
$E^{w_1}$ and $E^{w_2}$ define the same function on ${\bf Z}^n$ when
$w_1 - w_2 \in {\bf Z}^n$, so that we can think of $E^w(x)$ for $x \in
{\bf Z}^n$ as being defined for each $w \in {\bf T}^n = {\bf R}^n /
{\bf Z}^n$.  If $\lambda$ is a bounded nonnegative linear functional
on $C_{00}({\bf Z}^n)$, then the Fourier transform
$\widehat{\lambda}(w)$ is now viewed as a function on ${\bf T}^n$,
defined as before by $\widehat{\lambda}(w) = \lambda(E^{-w})$, and one
can check that it is continuous.  Once again we have that $\lambda *
E^w = \widehat{\lambda}(w) \, E^w$ when $\lambda$ is a bounded
nonnegative linear functional on $C_{00}({\bf Z}^n)$ and $w \in {\bf
T}^n$.  If $\lambda_1$, $\lambda_2$ are nonnegative linear functionals
on $C_{00}({\bf Z}^n)$, then we can define the convolution $\lambda_1
* \lambda_2$ as a nonnegative linear functional on ${\bf Z}^n$ if at
least one of $\lambda_1$, $\lambda_2$ has finite support, or if both
$\lambda_1$, $\lambda_2$ are bounded.  If both $\lambda_1$,
$\lambda_2$ are bounded nonnegative linear functionals on $C_{00}({\bf
Z}^n)$, then the convolution $\lambda_1 * \lambda_2$ is a bounded
nonnegative linear functional on $C_{00}({\bf Z}^n)$ too.  If
$\lambda_1$, $\lambda_2$ both have finite support, then the
convolution $\lambda_1 * \lambda_2$ has finite support.  We also have
that $\lambda_1 * \lambda_2 = \lambda_2 * \lambda_1$ and $(\lambda_1 *
\lambda_2) * f = \lambda_1 * (\lambda_2 * f)$ in the cases where these
are defined.  If $\lambda_1$, $\lambda_2$ are bounded nonnegative
linear functionals on $C_{00}({\bf Z}^n)$, then $(\lambda_1 *
\lambda_2)^{\widehat{}}(w) = \widehat{\lambda}_1(w) \,
\widehat{\lambda}_2(w)$ for all $w \in {\bf T}^n$.

\section{Riemann--Stieltjes integrals}
\label{section on Riemann--Stieltjes integrals}
\index{Riemann--Stieltjes integrals}

	Let $a$, $b$ be real numbers with $a < b$, and let $\mu(x)$
be a monotone increasing function on the interval $[a, b]$.  Thus
we have
\begin{equation}
\label{mu(a) le mu(x) le mu(b)}
	\mu(a) \le \mu(x) \le \mu(b) 
\end{equation}
for all $x \in [a, b]$.

	Suppose that $f(x)$ is a continuous real-valued function on
$[a, b]$.  If $P = \{t_j\}_{j=0}^n$ is a partition of $[a, b]$, so
that
\begin{equation}
	a = t_0 < \cdots < t_n = b,
\end{equation}
then we define the \emph{mesh}\index{mesh of a partition of an interval}
of $P$ to be the maximum of $t_j - t_{j-1}$, $1 \le j \le n$.  To
such a partition $P$ and the functions $f$, $\mu$, we associate the
sum
\begin{equation}
	\sum_{j=1}^n f(t_j) \, (\mu(t_j) - \mu(t_{j-1})).
\end{equation}
Using the uniform continuity of $f$ on $[a, b]$ and the completeness
of the real numbers one can show that there is a real number $A$ such
that for every $\epsilon > 0$ there is a $\delta > 0$ so that
the absolute value of the sum above minus $A$ is less than $\epsilon$
whenever the mesh of the partition $P$ is less than $\delta$.
This number $A$ is called the Riemann--Stieltjes integral of $f$
with respect to $\mu$ on $[a, b]$, and it is denoted
\begin{equation}
	\int_a^b f(x) \, d\mu(x).
\end{equation}

	Now suppose that $\mu(x)$ is a monotone increasing real-valued
function on the whole real line.  Of course if one starts with $\mu(x)$
defined on a closed interval $[a, b]$ as before, then we can easily
extend $\mu(x)$ to a monotone increasing function on the whole real
line by setting $\mu(x) = \mu(a)$ when $x < a$, $\mu(x) = \mu(b)$ when
$x > b$.  If $f(x)$ is an element of $C_{00}({\bf R})$, then we
obtain the Riemann--Stieltjes integral
\begin{equation}
	\int_{\bf R} f(x) \, d\mu(x)
\end{equation}
of $f$ with respect to $\mu$ on the whole real line simply by taking
the integral $\int_a^b f(x) \, d\mu(x)$ for some $a$, $b$ such that
the support of $f$ is contained in the interval $[a, b]$.

	The integral of $f(x)$ with respect to $\mu$ on ${\bf R}$
defines a nonnegative linear functional on $C_{00}(M)$.  It is not
difficult to see that this nonnegative linear functional is bounded if
and only if $\mu$ is bounded, and more precisely it is bounded with
constant equal to 
\begin{equation}
	\sup_{x \in R} \mu(x) - \inf_{x \in R} \mu(x)
\end{equation}
when $\mu$ is bounded.

	It can be shown that every nonnegative linear functional on
$C_{00}({\bf R})$ arises in this manner, i.e., as a Riemann--Stieltjes
integral with respect to a monotone increasing function $\mu(x)$ on
${\bf R}$.

\section{Some geometry and analysis on metric spaces}
\label{section on some geometry and analysis on metric spaces}

	Let $(M, d(x,y))$ be a metric space.  Suppose that $A$ is a
nonempty finite set, and that for each $a \in A$ we have a
nonempty bounded subset $E_a$ of $M$.  Define
$\widehat{E}_a$ by
\begin{equation}
	\widehat{E}_a 
		= \{x \in M : \dist(x, E_a) \le \diam E_a\}.
\end{equation}
Thus
\begin{equation}
	\diam \widehat{E}_a \le 3 \diam E_a.
\end{equation}
	
	Let us define a subset $A_1$ of $A$ as follows.  First choose
an $a_1 \in A$ such that $\diam E_{a_1}$ is as large as
possible.  There may be more than one such $a_1$, and one simply
picks one.  Next, choose an $a_2 \in A$ such that $E_{a_2}$
is disjoint from $E_{a_1}$ and $\diam E_{a_2}$ is as large
as possible.  In general, if $a_1, \ldots, a_l$ have been
chosen, then one chooses $a_{l+1} \in A$ so that
$E_{a_{l+1}}$ is disjoint from $E_{a_j}$, $1 \le j \le l$,
and $\diam E_{a_{l+1}}$ is as large as possible.  If there are no
$a$'s in $A$ such that $E_a$ is disjoint from
$E_{a_j}$, $1 \le j \le l$, then we simply stop the process.  We
take $A_1$ to be the set of $a_j$'s in $A$ chosen in this manner.

	By construction,
\begin{equation}
\label{disjointness property}
	E_a \cap E_b = \emptyset 
		\quad\hbox{when } a, b \in A_1, a \ne b.
\end{equation}
Furthermore,
\begin{equation}
\label{covering property}
	\bigcup_{a \in A} E_a 
		\subseteq \bigcup_{b \in A_1} \widehat{E}_b.
\end{equation}
To see this, let $a \in A$ be given.  If $a \in A_1$, then
$E_a$ is certainly contained in the right side of (\ref{covering
property}), and so we assume that $a \in A \backslash A_1$.
From the construction of $A_1$ it follows that there is a $b \in A_1$
such that
\begin{equation}
	\diam E_a \le \diam E_b \quad\hbox{and}\quad
		E_a \cap E_b \ne \emptyset.
\end{equation}
This implies that $E_a \subseteq \widehat{E}_b$ for this choice
of $b$.

	Now let us assume that closed and bounded subsets of $M$ are
compact, and that $\lambda$ is a nonnegative linear functional
on $M$.  Let us assume also that $\lambda$ is bounded with constant
$C$, so that
\begin{equation}
\label{lambda bounded with constant C}
	|\lambda(f)| \le C \, \|f\|_{sup}
\end{equation}
for all $f$ in $C_{00}(M)$.  Note that if one starts with a
nonnegative linear functional $\lambda_1$ on $C_{00}(M)$ which is not
bounded, then one can get a bounded nonnegative linear functional
$\lambda$ on $C_{00}(M)$ by setting $\lambda(f) = \lambda_1(\phi \,
f)$, where $\phi$ is any nonnegative function in $C_{00}(M)$.

	For each $x \in M$, define $\mathcal{A}_x$ to be the set of
ordered pairs $(U, \phi)$, where $U$ is an open subset of $M$ such
that $x \in U$ and $\diam U > 0$, and $\phi$ is a function in
$C_{00}(M)$ such that $\supp \phi \subseteq U$ and $0 \le \phi(x) \le
1$ for all $x \in M$.  Fix a positive real number $\alpha$, and define
$\lambda_\alpha^*(x)$ for $x \in M$ by
\begin{equation}
	\lambda_\alpha^*(x) 
		= \sup \biggl\{\frac{\lambda(\phi)}{(\diam U)^\alpha} :
					(U, \phi) \in \mathcal{A}_x\}.
\end{equation}
Thus $\lambda_\alpha^*(x) \ge 0$, and $\lambda_\alpha^*(x) = +\infty$
is possible.

	Let $t$ be a positive real number, and consider the set
\begin{equation}
\label{{x in M : lambda_alpha^*(x) > t }}
	\{x \in M : \lambda_\alpha^*(x) > t \}.
\end{equation}
If $x$ is in this set, then there is a pair $(U, \phi) \in
\mathcal{A}_x$ such that
\begin{equation}
\label{frac{lambda(phi)}{(diam U)^alpha} > t}
	\frac{\lambda(\phi)}{(\diam U)^\alpha} > t.
\end{equation}
This implies that $\lambda_\alpha^*(y) > t$ for all $y \in U$, and
hence that (\ref{{x in M : lambda_alpha^*(x) > t }}) is an open subset
of $M$.

	Suppose that $K$ is a compact subset of the open set
(\ref{{x in M : lambda_alpha^*(x) > t }}).  For each $x \in K$,
there is a $(U_x, \phi_x) \in \mathcal{A}_x$ such that
\begin{equation}
\label{frac{lambda(phi_x)}{(diam U_x)^alpha} > t}
	\frac{\lambda(\phi_x)}{(\diam U_x)^\alpha} > t.
\end{equation}
By compactness there is a finite subset $F$ of $K$ such that
\begin{equation}
	K \subseteq \bigcup_{x \in F} U_x.
\end{equation}
The remarks at the beginning of the section imply that there is a
subset $F_1$ of $F$ such that the $U_x$'s for $x \in F_1$ are
pairwise-disjoint and
\begin{equation}
	K \subseteq \bigcup_{x \in F} U_x 
		\subseteq \bigcup_{y \in F_1} \widehat{U}_y,
\end{equation}
where $\widehat{U}_y$ is as defined before.  Hence
\begin{eqnarray}
	HF^\alpha_{con}(K) & \le & 
			\sum_{y \in F_1} (\diam \widehat{U}_y)^\alpha	\\
		& \le & \sum_{y \in F_1} 3^\alpha \, (\diam U_y)^\alpha
								\nonumber \\
		& \le & \sum_{y \in F_1} 3^\alpha \, t^{-1} \lambda(\phi_y)
								\nonumber \\
		& \le & 3^\alpha \, t^{-1} \, C,		\nonumber
\end{eqnarray}
where $C$ is as in (\ref{lambda bounded with constant C}).  This
uses (\ref{frac{lambda(phi_x)}{(diam U_x)^alpha} > t}) and the fact
that
\begin{equation}
	0 \le \sum_{y \in F_1} \phi_y(w) \le 1
\end{equation}
for all $z \in M$, since $0 \le \phi_y(w) \le 1$ for all $w \in M$, and
$\supp \phi_y \subseteq U_y$ and the $U_y$'s, $y \in F_1$, are
pairwise-disjoint.

\backmatter

\newpage

\addcontentsline{toc}{chapter}{Index}
\printindex

\end{document}